%% file: reflection.tex
\documentclass[10pt,leqno]{amsart}

\newtheorem{lemma}{Lemma} 
\newtheorem{theorem}{Theorem}
\newtheorem{corollary}{Corollary}

\newtheorem{conjecture}{Conjecture} 

\input{amssymb.sty}
\usepackage{graphics} 
\usepackage{epsfig}
\usepackage{enumerate}
\usepackage[french]{babel}
\usepackage[T1]{fontenc} 

\def\C{{\mathbb C}}
\def\N{{\mathbb N}}
\def\R{{\mathbb R}}
\def\V{{\mathbb V}}
\def\dl{[\![}
\def\dr{]\!]}

\begin{document}


\title[Envelopes of holomorphy and the reflection principle]
{On envelopes of holomorphy of domains covered by Levi-flat hats
and the reflection principle}

\author{Jo\"el Merker}

\address{CNRS, Universit\'e de Provence, CMI, 
39 rue Joliot-Curie, 13453 Marseille Cedex 13, France} 

\email{merker@cmi.univ-mrs.fr 00 33 / (0)4 91 11 36 72 / (0)4 91 53 99 05}



\date{\number\year-\number\month-\number\day}

\maketitle

\begin{abstract}
In the present paper, we associate the techniques of the Lewy-Pinchuk
reflection principle with the Behnke-Sommer continuity principle.
Extending a so-called {\it reflection function} to a parameterized
congruence of Segre varieties, we are led to studying the envelope of
holomorphy of a certain domain covered by a smooth Levi-flat ``hat''.
In our main theorem, we show that every $\mathcal{C}^\infty$-smooth CR
diffeomorphism $h: M\to M'$ between two globally minimal
real analytic hypersurfaces in $\C^n$ ($n\geq 2$) is real analytic at
every point of $M$ if $M'$ is holomorphically nondegenerate. More
generally, we establish that the reflection function $\mathcal{R}_h'$
associated to such a $\mathcal{C}^\infty$-smooth CR diffeomorphism
between two globally minimal hypersurfaces in $\C^n$ ($n\geq 1$)
always extends holomorphically to a neighborhood of the graph of $\bar
h$ in $M\times \overline{M}'$, without any nondegeneracy condition on
$M'$. This gives a new version of the Schwarz symmetry
principle to several complex variables. Finally, we show that every
$\mathcal{C}^\infty$-smooth CR mapping $h: M\to M'$ between
two real analytic hypersurfaces containing no complex
curves is real analytic at every point of $M$, without
any rank condition on $h$.

\medskip
\noindent
{\sc R\'esum\'e}.
Dans cet article, nous associons les techniques du principe de
r\'eflexion de Lewy-Pinchuk avec celles du principe de continuit\'e de
Behnke-Sommer.  Apr\`es avoir prolong\'e holomorphiquement une
fonction dite ``de r\'eflexion'' \`a une congruence de
sous-vari\'et\'es de Segre, nous sommes conduits \`a l'\'etude de
l'enveloppe d'holomorphie d'un domaine recouvert d'un ``chapeau''
Levi-plat lisse. D'apr\`es notre r\'esultat principal, tout
CR-diff\'eomorphisme $h: M\to M'$ de classe $\mathcal{C}^\infty$ entre
deux hypersurfaces analytiques r\'eelles globalement minimales de
$\C^n$ ($n\geq 2$) est analytique r\'eel en tout point de $M$ si $M'$
est holomorphiquement non-d\'eg\'en\'er\'ee. Plus g\'en\'eralement,
nous \'etablissons que la fonction de r\'eflexion $\mathcal{R}_h'$
associ\'ee \`a un tel diff\'eomorphisme CR de classe
$\mathcal{C}^\infty$ entre deux hypersurfaces analytiques r\'eelles
globalement minimales se prolonge toujours holomorphiquement \`a un
voisinage du graphe de $\bar h$ dans $M\times \overline{M}'$, sans
aucune condition de non-d\'eg\'en\'erescence sur $M'$. Cet \'enonc\'e
fournit une nouvelle version du principe de r\'eflexion de Schwarz en
plusieurs variables complexes. Enfin, nous d\'emontrons que toute
application $h:M\to M'$ de classe $\mathcal{C}^\infty$ et CR entre
deux hypersurfaces analytiques r\'eelles ne contenant pas de courbes
holomorphes est analytique r\'eelle en tout point de $M$, sans
aucune condition de rang sur $h$.
\end{abstract}

\smallskip

\begin{center}
\begin{minipage}[t]{10cm}
\baselineskip =0.35cm
{\scriptsize

\centerline{\bf Table of contents~:}

\smallskip

{\bf \S1.~Introduction and presentation of the results \dotfill 2.}

{\bf \S2.~Description of the proof of Theorem~1.2 \dotfill 8.}

{\bf \S3.~Biholomorphic invariance of the reflection function \dotfill 13.}

{\bf \S4.~Extension across a Zariski dense open subset of $M$ \dotfill 22.}

{\bf \S5.~Situation at a typical point of non-analyticity \dotfill 26.}

{\bf \S6.~Envelopes of holomorphy of domains with Levi-flat hats \dotfill 30.}

{\bf \S7.~Holomorphic extension to a Levi-flat union of Segre varieties \dotfill 34.}

{\bf \S8.~Relative position of the neighbouring Segre varieties \dotfill 38.}

{\bf \S9.~Analyticity of some degenerate $\mathcal{C}^\infty$-smooth CR mappings 
\dotfill 39.}

{\bf \S10.~Open problems and conjectures 
\dotfill 46.}

}\end{minipage}
\end{center}

\bigskip

\section*{\S1.~Introduction and presentation of the results}

\subsection*{1.1.~Main theorem}
Let $h:M\to M'$ be a $\mathcal{C}^\infty$-smooth CR diffeomorphism between
two geometrically smooth real analytic hypersurfaces in $\C^n$ $(n\geq
2)$. Call $M$ {\it globally minimal} (in the sense of
Tr\'epreau-Tumanov) if it consists of a single CR orbit ([Tr1,2],
[Tu1,2], [Me1], [MP1]). Call $M'$ {\it holomorphically nondegenerate}
(in the sense of Stanton) if there does not exist any nonzero $(1,0)$
vector field with holomorphic coefficients which is tangent to a
nonempty open subset of $M$ ([St1,2]).  Our principal result is as
follows.

\def\thetheorem{1.2}\begin{theorem}
If $M$ is globally minimal and if $M'$ is holomorphically
nondegenerate, then the $\mathcal{C}^\infty$-smooth CR diffeomorphism
$h$ is real analytic at every point of $M$.
\end{theorem}

Compared to classical results of the literature, in this theorem, no
pointwise, local or not propagating nondegeneracy condition is imposed
on $M'$, like for instance $M'$ be Levi nondegenerate, finitely
nondegenerate or essentially finite at every point. With respect to
the contemporary state of the art, the novelty in Theorem~1.2 lies in
the treatment of the locus of non-essentially finite points, which is
a proper real analytic subvariety of $M'$, provided $M'$ is
holomorphically nondegenerate.  There is also an interesting invariant
to study, more general than $h$, namely the {\it reflection function}
$\mathcal{R}_h'$.  Because the precise definition of $\mathcal{R}_h'$
involves a concrete defining equation of $M'$, it must be localized
around various points $p'\in M'$, so we refer to \S1.7 below for a
complete presentation. Generalizing Theorem~1.2, we show that
$\mathcal{R}_h'$ extends holomorphically to a neighborhood of each
point $(p,\overline{h(p)})\in M\times \overline{M}'$, assuming only
that $M$ is globally minimal and without any nondegeneracy condition
on $M'$ (Theorem~1.9). We deduce in fact Theorem~1.2 from the
extendability of $\mathcal{R}_h'$. This strategy of proof is inspired
from the deep works of Diederich-Pinchuk [DP1,2] ({\it see} also [V],
[Sha], [PV]) where the extension as a mapping is derived from the
extension as a correspondence.

In the sequel, we shall by convention sometimes denote by $(M,p)$ a
{\it small connected piece} of $M$ localized around a ``center'' point
$p\in M$. However, since all our considerations are semi-local and of
geometric nature, we shall {\it never} use the language
of germs.

\subsection*{1.3.~Development of the classical results and brief history}
The earliest extension result like Theorem~1.2 was found independently
by Pinchuk [P3] and after by Lewy [L]: if $(M,p)$ and
$(M',p')$ are strongly pseudoconvex, then $h$ is real analytic at
$p$. The classical proof in [P3] and [L] makes use of the so-called
{\it reflection principle} which consists to {\it solve}\, first the 
mapping $h$
with respect to the jets of $\bar h$ (by this, 
we mean a relation like $h(q)=\Omega(q,\bar q,
j^k\bar h(\bar q))$ where $\Omega$ is holomorphic
in its arguments and $q\in M$, {\it cf.}~(4.10) below) and to apply
afterwards the one-dimensional Schwarz symmetry principle in a
foliated union of transverse holomorphic discs. In 1978 and in
1982, Webster [W2,3] extended this result to Levi nondegenerate
CR manifolds of higher codimension.  Generalizing this principle,
Diederich-Webster proved in 1980 that a sufficiently smooth CR
diffeomorphism is analytic at $p\in M$ if $M$ is generically
Levi-nondegenerate and the morphism of jets of Segre varieties of $M'$
is injective ({\it see} \S2 of the fundamental article [DW] and~(1.11)
below for a definition of the Segre morphism). In 1983, Han [Ha]
generalized the reflection principle for CR diffeomorphisms between
what is today called {\it finitely nondegenerate} hypersurfaces ({\it
see} [BER2]).  In 1985, Derridj [De] studied the reflection principle
for proper mappings between some model classes of weakly pseudoconvex
boundaries in $\C^2$.  In 1985, Baouendi-Jacobowitz-Treves [BJT]
proved that every $\mathcal{C}^\infty$-smooth CR diffeomorphism $h:
(M,p)\to (M',p')$ between two real analytic CR-generic manifolds in
$\C^n$ which extends holomorphically to a fixed wedge of edge $M$, is
real analytic, provided $(M',p')$ is essentially finite. After the
work of Rea [R], in which holomorphic extension to one side of CR
functions on a minimal real analytic hypersurface was proved (the
weakly pseudoconvex case, which is not very different, was treated
long before in a short note by Bedford-Forn{\ae}ss [BeFo]; {\it see}\,
also [BT2]), after the work of Tumanov [Tu1], who proved wedge
extendability in general codimension, and after the work of
Baouendi-Rothschild [BR3], who proved the necessity of minimality for
wedge extension (in the meanwhile, Treves provided a simpler argument
of necessity), it was known that the automatic holomorphic extension
to a fixed wedge of the components of $h$ holds {\it if and only if}\,
$(M,p)$ is minimal in the sense of Tumanov. Thus, the optimal
extendability result in [Tu1] strengthened considerably the main
theorem of [BJT]. In the late eighties, the research on the
analyticity of CR mappings has been pursued by many authors
intensively. In 1987--88, Diederich-Forn{\ae}ss [DF2] and
afterwards (not independently)
Baouendi-Rothschild [BR1] extended this kind of reflection principle
to the non diffeomorphic case, namely for a
$\mathcal{C}^\infty$-smooth CR mapping $h$ between two essentially finite
hypersurfaces which is locally finite to one, or locally proper. This
result was generalized in [BR2] to $\mathcal{C}^\infty$-smooth
mappings $h\,:(M,p)\to (M',p')$ whose formal Jacobian determinant at
$p$ does not vanish identically, again with $(M',p')$ essentially
finite. In 1993-6, Sukhov [Su1,2] and Sharipov-Sukhov [SS] generalized the
reflection principle of Webster in [W2,3] by introducing a global
condition on the mapping, called Levi-transversality.  Following this
circle of ideas, Coupet-Pinchuk-Sukhov have pointed out in their
recent works [CPS1,2] that almost all the above-mentioned variations
on the reflection principle find a unified explanation in the fact
that a certain complex analytic variety $\V_p'$ is zero-dimensional,
which intuitively speaking means that $h$ is {\it finitely determined
by the jets of $\bar h$}, {\it i.e.} more precisely that each
components $h_j$ of $h$ satisfies a monic Weierstrass polynomial having
analytic functions depending on a finite jet of $\bar h$ as
coefficients (this observation appears also in [Me3]). They stated
thus a general result in the hypersurface case whose extension to a
higher codimensional minimal CR-generic source $(M,p)$ was achieved
recently by Damour in [Da2]. In sum, this last clarified unification
closes up what is attainable in the spirit of the so-called {\it
polynomial identities} introduced in [BJT], yielding a quite general
sufficient condition for the analyticity of $h$.  In the arbitrary
codimensional case, this general sufficient condition can be expressed
simply as follows. Let
$\overline{L}_1,\cdots, \overline{L}_m$ be a basis of $T^{0,1}M$,
denote $\overline{L}^\beta:= \overline{L}_1^{\beta_1}\dots
\overline{L}_{m}^{\beta_{m}}$ for $\beta\in\N^{m}$ and let
$\rho_{j'}'(t',\bar t')=0$, $1\leq j'\leq d'$, be a collection of real
analytic defining equations for a generic $(M',p')$ of codimension
$d'$. Then the complex analytic variety, called the (first) {\it
characteristic variety} in [CPS1,2], [Da1,2]:
\def\theequation{1.4}\begin{equation}
\V_p':=\{t'\in \C^n: \overline{L}^\beta [\rho'(t',\bar h(\bar t))]
\vert_{\bar t=\bar p}=0, \ 
\forall \, \beta \in \N^{m}\}.
\end{equation}
is always zero-dimensional at $p'\in \V_p'$ in [L], [P3], [W1], [W2],
[W3], [DW], [Ha], [De], [BJT], [DF2], [BR1], [BR2], [BR4], [Su1,2],
[BHR], [Su1], [Su2], [SS], [BER1], [BER2], [CPS1], [CPS2], [Da] (in
[P4], [DFY], [DP1,2], [V], [Sha], [PV], the variety $\V_p'$ is not
defined because these authors tackle the much more difficult problem
where {\it no initial regularity assumption}\, is supposed on the
mapping; in [DF2], some cases of non-essentially finite hypersurfaces
are admitted). Importantly, the condition $\dim_{p'} \V_p'=0$ requires
$(M',p')$ to be essentially finite.

\subsection*{1.5.~Non-essentially finite hypersurfaces}
However, it is known that the finest CR-regularity phenomena come down
to the consideration of a class of much more general hypersurfaces
which are called {\it holomorphically nondegenerate}\, by Stanton
[St1,2] and which are in general {\it not essentially finite}. In
1995, Baouendi-Rothschild [BR3] exhibited this condition as a {\it
necessary and sufficient} condition for the algebraicity of a local
biholomorphism between two real algebraic hypersurfaces. Thanks to the
nonlocality of algebraic objects, they could assume that $(M',p')$ is
essentially finite after a small shift of $p'$, which entails again
$\dim_{p'} \V_p'=0$, thus reducing the work
to the application of known techniques (even in
fact simpler, in the generalization to the higher codimensional case,
Baouendi-Ebenfelt-Rothschild came down to a direct application of the
algebraic implicit function theorem by solving algebraically $h$ with
respect to the jets of $\bar h$ [BER1];
Since then however, few
works have been devoted to the study of the analytic regularity of
smooth CR mapping between {\it non-essentially finite}\, hypersurfaces
in $\C^n$. It is well known that the main technical difficulties in
the subject happen to occur in $\C^n$ for $n\geq 3$ and that a great
deal of the obstacles which one naturally encounters can be avoided by
assuming that the target hypersurface $M$' is algebraic (with $M$
algebraic or real analytic), see {\it e.g.} the works [MM2],
[Mi1,2,3], [CPS1] (in case $M'$ is algebraic, its Segre varieties are
defined all over the compactification $P_n(\C)$ of $\C^n$, which
helps much). Finally, we would like to mention the papers of Meylan
[Mey], Maire and Meylan [MaMe], Meylan and the author [MM1], 
Huang, the author and Meylan [HMM] in this
respect (nevertheless, after division by a suitable holomorphic
function, the situation under study in these works is again reduced to
polynomial identities).

\subsection*{1.6.~Schwarz's reflection principle in higher dimension}
In late 1996, seeking a natural generalization of Schwarz's reflection
principle to higher dimension and inspired by the article
[DP1], the author ({\it see} [MM2], [Me3]) pointed out the interest of
the so-called {\it reflection function} $\mathcal{R}_h'$ associated
with $h$. This terminology is introduced {\it passim}\, in
[Hu,~p.~1802]; a different definition involving one more
variable is given in [Me3,5,6,7,8]; the biholomorphic invariance of
$\mathcal{R}_h'$ and the important observation that $\mathcal{R}_h'$
should extend holomorphically {\it without any nondegeneracy
condition} on $(M',p')$ appeared for the first time in the preprint versions
of [MM2], [Me3], which inspired the papers [Mi1,2].

 Indeed, the explicit expression of this function depends on a
local defining equation for $M'$, but its holomorphic extendability is
independent of coordinates and there are canonical rules of
transformation between two reflection functions ({\it see} \S3 below).
As the author believes, in the diffeomorphic case and provided
$M$ is at least globally minimal, {\it this function should extend
without assuming any nondegeneracy condition on $M'$, in pure analogy
with the Schwarzian case $n=1$}. It is easy to convince oneself that
the reflection function is the right invariant to study. In fact,
since then, it has been already studied thoroughly in the algebraic
and in the formal CR-regularity problems, {\it see} [Me3,5,6,7,8],
[Mi2,3,4]. For instance, the formal reflection mapping associated with
a formal CR equivalence between two real analytic CR-generic manifolds
in $\C^n$ which are minimal in the sense of Tumanov is convergent
({\it see} [Mi3,4] for partial results in this direction and [Me6,7,8]
for the complete statement). If $h$ is a holomorphic equivalence
between two real algebraic CR-generic manifolds in $\C^n$ which are
minimal at a Zariski-generic point, then the reflection mapping
$\mathcal{R}_h'$ is algebraic ({\it see} [Mi2] for the hypersurface
case and [Me5] for arbitrary codimension).  So we expect that totally
similar statements hold for smooth mappings between real analytic CR
manifolds.

\subsection*{1.7.~Analyticity of the reflection function}
For our part, we deal in this paper with smooth CR mappings
between {\it hypersurfaces}. Thus, as above, let $h: M\to M'$ be a
$\mathcal{C}^\infty$-smooth CR mapping between two connected
real analytic hypersurfaces in $\C^n$ with $n\geq 2$.  We
shall constantly assume that $M$ is globally minimal.  Equivalently,
$M$ is locally minimal (in the sense of Tr\'epreau-Tumanov) 
at every point, since
$M$ is real analytic (however, there exist $\mathcal{C}^2$-smooth or
$\mathcal{C}^\infty$-smooth hypersurfaces in $\C^n$, $n\geq 2$, which
are globally minimal but not locally minimal at many point, {\it see} [J],
[MP1]). Postponing generalizations and refinements to further
investigation, we shall assume here for simplicity that $h$ is a CR
diffeomorphism.  Of course, in this case, the assumption of global
minimality of $(M,p)$ can then be switched to $(M',p')$. The
associated reflection function $\mathcal{R}_h'$ is a complex function
which is defined in a neighborhood of the graph of $\bar h$
in $\C^n\times\C^n$ as
follows. Localizing $M$ and $M'$ at points $p\in M$ and $p'\in M'$
with $p'=h(p)$, we choose a complex analytic defining equation for
$M'$ in the form $\bar w'=\Theta'(\bar z',t')$, where
$t'=(z',w')\in\C^{n-1}\times \C$ are holomorphic coordinates vanishing
at $p'$ and where the power series $\Theta'(\bar z',t'):=
\sum_{\beta\in\N^{n-1}} (\bar z')^\beta \, \Theta_\beta'(t')$ vanishes
at the origin and converges normally in a small polydisc
$\Delta_{2n-1}(0,\rho')=\{(\bar z',t')\,: \vert \bar z'\vert, \vert
t'\vert < \rho'\}$, where $\rho'>0$ and where
$\vert t'\vert:=\max (\vert t_1'\vert,\dots,
\vert t_n'\vert)$ is the polydisc norm. Here, by reality of $M'$, the
holomorphic function $\Theta'$ is not arbitrary, it must satisfy the
power series identity $\Theta'(\bar z',z',\overline{\Theta}'(z',\bar
z', \bar w'))\equiv \bar w'$. Conversely, such a power series
satisfying this identity does define a real analytic hypersurface
$\bar w'=\Theta'(\bar z',t')$ of $\C^n$ as can 
be verified easily ([BER2, Remark~4.2.30]).
It is important to notice that once the coordinate system $t'$ is
fixed, with the $w'$-axis not complex tangent to $M'$ at $0$, then
there is {\it only one}\, complex defining equation for $M'$ of the form
$\bar w'=\Theta'(\bar z',t')$. 

By definition,
the {\it reflection function} $\mathcal{R}_h'$ associated with $h$ and
with such a local defining function for $(M',p')$ is the following
function of $2n$ complex variables:
\def\theequation{1.8}\begin{equation}
(t,\bar\nu')\mapsto
\bar \mu'-\sum_{\beta\in \N^{n-1}} (\bar{\lambda'})^\beta \,
\Theta_\beta'(h(t))=:\mathcal{R}_h'(t,\bar \nu'),
\end{equation}
where $\bar\nu'=(\bar\lambda', \bar\mu')\in\C^{n-1}\times \C$.  It can
be checked rigorously that this function is CR and of class
$\mathcal{C}^\infty$ with respect to the variable $t\in M$ in a
neighborhood of $p$ and that it is holomorphic with respect to the variable
$\bar\nu'$ in the polydisc neighborhood $\{\vert z'\vert < \rho'\}$ of
$\bar p'$ in $\C^n$ ({\it see}\, Lemma~3.8 below). Let us call the
functions $\Theta_\beta'(h(t))$ the {\it components}\, of the
reflection function. Since $M$ is in particular minimal at the point
$p\in M$, the components $h_j$ of the mapping $h$ and hence also the
components $\Theta_\beta'(h(t))$ of $\mathcal{R}_h'$ extend
holomorphically to a one-sided neighborhood $D_p$ of $M$ at $p$,
obtained by gluing Bishop discs to $(M,p)$. Our first main result is
as follows.

\def\thetheorem{1.9}\begin{theorem}
If $h: M\to M'$ is a $\mathcal{C}^\infty$-smooth CR diffeomorphism
between two globally minimal real analytic hypersurfaces in $\C^n$,
then for every point $p\in M$ and for every choice of a coordinate
system vanishing at $p':=h(p)$ as above in which $(M',p')$ is
represented by $\bar w'=\Theta'(\bar z', t')$, the associated
reflection function $\mathcal{R}_h'(t,\bar\nu')=
\bar\mu'-\Theta'(\bar\lambda',h(t))$ centered at $p\times
\overline{p'}$ extends holomorphically to a neighborhood of $p\times
\overline{p'}$ in $\C^n\times \C^n$.
\end{theorem}

\smallskip
\noindent
In \S3 below, we provide some fundational material about the
reflection function. Especially, we prove that the
holomorphic extendability to a neighborhood of $p\times \bar p'$
does not depend on the choice of
a holomorphic coordinate system vanishing at $p'$.  By
differentiating~(1.8) with respect to $\bar \nu'$, we may observe that
the holomorphic extendability of $\mathcal{R}_h'$ to a neighborhood of
$p$ is equivalent to the following statement: {\it all the component
functions $\Theta_\beta'(h(t))=: \theta_\beta'(t)$ $($an infinite
number$)$ extend holomorphically to a fixed neighborhood of $p$ and
there exist constants $C,\rho,\rho'>0$ such that $\vert t\vert < \rho
\Rightarrow \vert \theta_\beta'(t)\vert < C\,(\rho')^{-\vert \beta\vert}$}
({\it see}\, Lemma~3.16 below). So
Theorem~1.9 may be interpreted as follows: instead of asserting that
the mapping $h$ extends holomorphically to a neighborhood of $p$, we
state that a certain invariant infinite collection of {\it
holomorphic functions of the components $h_j$ of the mapping} (which
depends directly on $M'$) {\it do extend holomorphically to a neighborhood
of $p$}. The important fact here is that we do not put any extra
nondegeneracy condition on $M'$ at $p'$ (except minimality). Another
geometric interpretation is as follows. Let $\underline{S}_{t'}':=\{
(\bar \lambda',\bar\mu')\in\C^n\,: \bar
\mu'=\Theta'(\bar\lambda',t')\}$ denote the {\it conjugate Segre variety}
associated with the fixed point having coordinates $t'$ (usually, to
define Segre varieties, one fixes instead the point $\bar \nu'$; nevertheless
conjugate Segre varieties are equally interesting, as argued in
[Me4]). Then Theorem~1.9 can be interpreted as saying that the not
rigorously defined intuitive ``Segre mapping'' $t\mapsto
\underline{S}_{h(t)}'$ extends holomorphically at $p$. In fact,
the target value of this mapping should be thought to be represented
concretely by the defining function of $\underline{S}_{h(t)}'$, namely
this intuitive ``Segre mapping'' {\it must}\, (and can only) be represented by
the {\it rigorous}\, reflection function $(t,\bar\nu')\mapsto \bar
\mu'-\Theta'(\bar\lambda',h(t))$. In sum, Theorem~1.9 precisely
asserts that the ``Segre mapping'' extends holomorphically to a
neighborhood of $p\times \bar p'$, without any nondegeneracy 
condition on $(M',p')$.  In certain circumstances, {\it
e.g.}  when $(M',p')$ is moreover assumed to be Levi-nondegenerate,
finitely nondegenerate or essentially finite, one may deduce
afterwards, thanks to the holomorphic extendability of the components
$\Theta_\beta'(h(t))$, that $h$ itself extends holomorphically at $p$
({\it cf.}  [DF2], [BR1], [DFY], [DP1,2], [V], [Sha],
[PV]). Analogously, in Theorem~1.14 below, we shall derive from
Theorem~1.9 above an important expected {\it necessary and sufficient
condition} for $h$ to be holomorphic at $p$.

\subsection*{1.10.~Applications} We give essentially two 
important applications.
Firstly, associated with $M$', there is an invariant integer
$\kappa_{M'}'$ with $0\leq \kappa_{M'}'\leq
n-1$, called the {\it holomorphic degeneracy degree} of
$M'$, which counts the maximal number of $(1,0)$ vector fields with
holomorphic coefficients defined in a neighborhood of $M'$ which are tangent
to $M'$ and which are linearly independent at a Zariski-generic point.
In particular, $M'$ is holomorphically nondegenerate if and only if
$\kappa_{M'}'=0$. Inspired by the geometric reflection principle
developed in [DW], [DF4], [F], we can provide another (equivalent)
definition of the integer $\kappa_{M'}'$ in terms of the morphism of
jets of Segre varieties as follows ({\it see} also [Me6,7,8]; historically,
finite order jets of $\mathcal{C}^\infty$-smooth CR mappings
together with finite order jets of the Segre morphism were first studied
in the reflection principle by Diederich-Forn{ae}ss in [DF4]). By
complexifying the variable $\bar t'$ as $(\bar t')^c=:\tau'$ and by
fixing $\tau'$, we may consider the {\it complexified Segre variety}\, which
is defined by $\mathcal{ S}_{\tau'}':=\{(w',z'):
w'=\overline{\Theta}'(z',\tau')\}$. For some supplementary information about
the canonical geometric correspondence between complexified Segre
varieties and complexified CR vector fields, we refer the interested
reader to [Me4,5]. Let $j_{t'}^k\mathcal{S}_{\tau'}'$ denote the
$k$-jet at the point $t'$ of $\mathcal{S}_{\tau'}'$.  This $k$-jet is
in fact defined by differentiating the defining equation of
$\mathcal{S}_{\tau'}'$ with respect to $z'$ as follows. For
$\beta\in\N^{n-1}$, we denote $\vert\beta\vert:=
\beta_1+\cdots+\beta_{n-1}$ and
$\partial_{z'}^\beta:=\partial_{z_1'}^{\beta_1}\cdots
\partial_{z_{n-1}'}^{\beta_{n-1}}$. Then the $k$-jet provides in fact
a holomorphic mapping which is defined over the extrinsic
complexification $\mathcal{M}':=\{(t',\tau'): w'-\overline{\Theta}'(
z',\tau')=0\}$ of $M'$ as shown in the following definition:
\def\theequation{1.11}\begin{equation}
j_k' :
{\mathcal{M}}' \ni (t', \tau') \mapsto j_{t'}^k\mathcal{
S}_{\tau'}':=(t',\{\partial_{z'}^\beta [w'-\overline{\Theta}'(z',\tau')]
\}_{\vert\beta\vert\leq k}) \in \C^{n+{(n-1+k)!\over (n-1)! \, k!}}.
\end{equation}
For $k$ large enough, 
the analytic properties of these jet mappings $j_k'$ govern the geometry
of $M'$, as was pointed out in [DW] for the first time. For instance,
Levi nondegeneracy, finite nondegeneracy and essential finiteness of
$(M',p')$ may be characterized in terms of the mappings $j_k'$ ([DW], [DF4],
[Me6,7,8]). In our case, it is clear that there exists an integer
$\chi_{M'}'$ with $1\leq \chi_{M'}'\leq n$ such that the generic rank
of $j_k'$ equals $n-1+\chi_{M'}'$ for all $k$ large enough, since the
generic ranks increase and are bounded by $2n-1$. Then the holomorphic
degeneracy degree can also be defined equivalently by
$\kappa_{M'}':=n-\chi_{M'}'$. We may notice in particular that $M'$ is
Levi-flat if and only if $\chi_{M'}'=1$, since
$\overline{\Theta}'(z',\tau')\equiv\tau_n'$ in this case. Consequently, 
we always
have $\chi_{M'}'\geq 2$ in this paper 
since we constantly assume that $M'$ is globally
minimal. The biholomorphic invariance of Segre varieties makes it easy
to precise in which sense the jet mapping $j_k'$ is invariantly attached
to $M'$, namely how it changes when one varies the coordinate
system. Then the fact that $\chi_{M'}'$ is defined in terms of the
generic rank of an invariant holomorphic mapping together with the
connectedness of $M'$ explains well that the integers $\chi_{M'}'$ and
$\kappa_{M'}'$ do not depend on the center point $p'\in M'$ in 
a neighborhood of which we define the mappings $j_k'$ (we prove
this in \S3). In particular, this explains why $M'$ is holomorphically
degenerate at one point if and only if it is holomorphically
degenerate at every point ([BR4]). On the contrary, the direct
definition of $\kappa_{M'}'$ in terms of locally defined tangent
holomorphic vector fields provided in [BR4], [BER2] makes this point
less transparent, even if the two definitions are equivalent. So, we
believe that the definitionof $\kappa_{M'}'$ 
in terms of $j_k'$ is more adequate.
Furthermore, to be even more concrete, let us add that the behavior of
the map~(1.11) depends mostly upon the infinite collection of
holomorphic mappings $(\overline{\Theta}_\beta'(\tau'))_{\beta\in\N^{n-1}}$,
since we essentially get rid of $z'$ by differentiating
$w'-\sum_{\beta\in\N^{n-1}} (z')^\beta\,\overline{\Theta}_\beta'(\tau')$
with respect to $z'$ in~(1.11). Equivalently, after conjugating, we
may consider instead the simpler holomorphic mappings $\mathcal{Q}_k':
t'\ni\C^n\mapsto (\Theta_\beta'(t'))_{\vert\beta\vert \leq k}\in
\C^{{(n-1+k)!\over (n-1)! \, k!}}$. Then the generic rank of
$\mathcal{Q}_k'$ is equal to the same integer $\chi_{M'}'$, for all
$k$ large enough.  This again supports the thesis that the components
$\Theta_\beta'(t')$ occuring in the defining function of $(M',p')$ and
in the reflection function are over all important.  In \S3 below, some
more explanations about the mappings $\mathcal{Q}_k'$ are provided.

Let $\chi_{M'}'$ be as above and let
$\Delta$ be the unit disc in $\C$. It is known that there exists a proper
real analytic subset $E_{M'}'$ of $M'$ such that for each point $q'\in
M'\backslash E_{M'}'$, there exists a neighborhood of $q'$ in $\C^n$
in which $(M',q')$ is biholomorphically equivalent to a product
$\underline{M}_{q'}'\times \Delta^{n-\chi_{M'}'}$ of a small real
analytic hypersurface $\underline{M}_{q'}'$ contained in the smaller
complex space $\C^{\chi_{M'}'}$ by a $(n-\chi_{M'}')$-dimensional
polydisc. As expected of course, the hypersurface
$\underline{M}_{q'}'$ is a {\it holomorphically nondegenerate}
hypersurface (Lemma~3.54), namely $\kappa_{\underline{M}_{q'}'}'=0$.
Now, granted Theorem~1.9, we observe that the local graph
$\{(t,h(t)): t\in (M,p)\}$ of $h$ is clearly contained in
the following local complex analytic set passing through $p\times p'$:
\def\theequation{1.12}\begin{equation}
\mathcal{C}_h':=\{(t,t')\in \C^n\times \C^n:
\Theta_\beta'(t')=\theta_\beta'(t), \ \forall
\, \beta\in\N^{n-1}\}.
\end{equation} 
It follows from the considerations of \S3 below that the various local
complex analytic sets $\mathcal{C}_h'$ centered at points
$(p, h(p))$ stick together in a well defined complex analytic set,
independent of coordinates.  Furthermore, since the generic rank of
$\mathcal{Q}_\infty'$ is equal to $\chi_{M'}'$, there exists a well defined
irreducible component $\mathcal{C}_h''$ of $\mathcal{ C}_h'$ of
dimension $(2n-\chi_{M'}')$ containing the local graph of $h$. We
deduce:

\def\thecorollary{1.13}\begin{corollary}
Let $(n-\chi_{M'}')$ be the holomorphic degeneracy degree of $M'$. Then
there exists a semi-global closed complex analytic subset
$\mathcal{C}_h''$ defined in a neighborhood of the graph of $h$ in
$\C^n\times \C^n$ which is of dimension $(2n-\chi_{M'}')$ and which
contains the graph of $h$ over $M$. In particular, $h$ extends as a
complex analytic set to a neighborhood of $M$ if $\chi_{M'}'=n$,
{\it i.e.} if $M'$ is holomorphically nondegenerate.
\end{corollary}

\smallskip
Of course, the most interesting case of Corollary~1.13 is when
$\chi_{M'}'=n$. Extendability of $h$ as an analytic set can be
improved.  Using the approximation theorem of Artin ([Ar]) we shall
deduce the following expected result ({\it see}\, Lemma~4.14), which
is identical with Theorem~1.2:

\def\thetheorem{1.14}\begin{theorem}
Let $h: M\to M'$ be a $\mathcal{C}^\infty$-smooth CR diffeomorphism
between two connected globally minimal real analytic hypersurfaces in
$\C^n$. If $M'$ is holomorphically nondegenerate, then $h$ is real
analytic at every point of $M$.
\end{theorem}

\noindent
Of course, real analyticity of $h$ is equivalent to its holomorphic
extendability to a neighborhood of $M$ in $\C^n$, by a classical
theorem due to Severi and generalized to higher codimension by
Tomassini. In particular, Theorem~1.14 entails that a pair of globally
minimal holomorphically nondegenerate real analytic hypersurfaces in
$\C^n$ are $\mathcal{C}^\infty$-smoothly CR equivalent {\it if and only if}\,
they are biholomorphically equivalent.

\subsection*{1.15.~Necessity}
Since 1995-6 ({\it see} [BR4], [BHR]), it was known that Theorem~1.14
above might provide an expected {\it necessary and sufficient
condition} for $h$ be analytic (provided of course that the local
CR-envelope of holomorphy of $M$, which already contains one side
$D_p$ of $M$ at $p$, does not contain the other side). Indeed,
considering self-mappings of $M'$, we have:

\def\thelemma{1.16}\begin{lemma}
\text{\rm ([BHR])}
Conversely, if $(M',p')$ is holomorphically {\rm degenerate} and if
there exists a $\mathcal{C}^\infty$-smooth CR function 
defined in a neighborhood of $p'\in M'$ which does {\rm not}
extend holomorphically to a neighborhood of $p'$, then there exists a
$\mathcal{C}^\infty$-smooth CR-automorphism of $(M',p')$ fixing $p'$
which is {\rm not} real analytic at $p'$.
\end{lemma}

\subsection*{1.17.~Organization of the paper}
To be brief, in \S2 we present first a thorough intuitive description
(in words) of our strategy for the proof of Theorems~1.2 and~1.9. This
presentation is really important, since it helps to understand the
general point of view without entering 
excessively technical considerations. Then \S3,
\S4, \S5, \S6, \S7 and \S8 are devoted to complete all the proofs.  We
would like to mention that in the last \S9, we provide a proof of the
following assertion, which might be interesting in itself, because it
holds without any rank assumption on $h$.  We refer the reader to the
beginning of \S9 for comments, generalizations and applications.

\def\thetheorem{1.18}\begin{theorem}
Let $h: M\to M'$ be a $\mathcal{C}^\infty$-smooth CR mapping between
two connected real analytic hypersurfaces in $\C^n$ $(n\geq 2)$.  If
$M$ and $M'$ do not contain any complex curve, then $h$ is
real analytic at {\rm every}\, point of $M$.
\end{theorem}

\subsection*{1.19.~Acknowledgement}
The author is very grateful to Egmont Porten, who poin\-ted out to him
the interest of gluing half-discs to the Levi flat hypersurfaces
$\Sigma_\gamma$ below. Also, the author wishes to thank Herv\'e
Gaussier and the referee for clever and helpful suggestions concerning
this paper.

\section*{\S2. Description of the proof of Theorem~1.2}

\subsection*{2.1.~Continuity principle and reflection principle}
According to the extendability theorem proved in [R], [BT2] and generalized
to only $\mathcal{C}^2$-smooth hypersurfaces by Tr\'epreau [Tr1], for
every point $p\in M$, the mapping $h$ in Theorems~1.9 and~1.14 
already extends
holomorphically to a one-sided neighborhood $D_p$ of $M$ at $p$ in
$\C^n$. This extension is performed by using small Bishop discs
attached to $M$ and by applying the approximation theorem proved in
[BT1]. These $D_p$ may be glued to yield a domain $D$ attached to $M$
which contains at least one side of $M$ at every point. In this
concern, we would like to remind the reader of the well known and
somewhat ``paradoxical'' phenomenon of {\it automatic holomorphic
extension of CR functions on $M$ to both sides}, which can render the
above Theorem~1.9 surprisingly trivial. Indeed, let $U_M$ denote the
(open) set of points $q$ in $M$ such that the envelope of holomorphy
of $D$ contains a neighborhood of $q$ in $\C^n$ (as is well known, if,
for instance, the Levi form of $M$ has one positive and one negative
eigenvalue at $q$, then $q\in U_M$; more generally, the {\it local
envelope of holomorphy} of $M$ or of the one-sided neighborhood $D$ of
$M$ at an arbitrary point $q\in M$ is always {\it one-sheeted}, as can
be established using the approximation theorem proved in [BT1]). Then
clearly, the $n$ components $h_1,\dots,h_n$ of our CR diffeomorphism
extend holomorphically to a neighborhood of $U_M$ in $\C^n$, as does
any arbitrary CR function on $M$. But it remains to extend $h$
holomorphically across $M\backslash U_M$ and the techniques of the
reflection principle are then unavoidable. {\it Here lies the
``paradox''}: sometimes the envelope of holomorphy trivializes the
problem, sometimes near some pseudoconvex points of finite D'Angelo
type (but not all) it helps to control the behavior of the mapping
thanks to local peak functions, sometimes it does not help at all,
especially at every point of the ``border'' between the pseudoconvex
and the pseudoconcave parts of $M$. In the interesting articles
[DF2,3], Diederich-Forn{\ae}ss succeeded in constructing the local
envelope of holomorphy at many points of a real analytic non-pseudoconvex
bounded boundary in $\C^2$ for which the border consists of a compact
maximally real submanifold and they deduced that any biholomorphic
mapping between two such domains extends continuously up to the
boundary as a CR homeomorphism.  In general, it is desirable to
describe constructively the local envelope of holomorphy at every
point of the border of $M$.  However, this general problem seems to be out of
the reach of the presently known techniques of study of envelopes of
holomorphy by means of analytic discs.  Fortunately, in the study of
the smooth reflection principle, the classical techniques usually do
not make any difference between the two sets $U_M$ and $M\backslash
U_M$ and these techniques provide a {\it uniform method} of extending
$h$ across $M$, no matter the reference point $p$ belongs to $U_M$ or
to $M\backslash U_M$ ({\it see} [L], [P3,4], [W1], [W2], [DW], [W3],
[BJT], [BR1], [BR2], [DF2], [Su1], [Su2], [SS], [BHR], [BER1], [BER2],
[CPS1], [CPS2]). Such a uniform method seems to be quite
satisfactory. On the other hand, the recent far reaching works of
Diederich-Pinchuk in the study of the {\it geometric reflection
principle} show up an accurate analysis of the relative
pseudo-convex(-concave) loci of $M$. Such an analysis originated in
the works of Diederich-Forn{\ae}ss [DF2,3] and in the work of
Diederich-Forn{\ae}ss-Ye [DFY]. In [P4], [DP1,2], [Hu], [Sha], the
authors achieve the propagation of holomorphic extension of a ``germ''
along the Segre varieties of $M$ (or the Segre sets), taking into
account their relative position with respect to $M$ and its local
convexity. In such reasonings, various discussions concerning
envelopes of holomorphy come down naturally in the proofs (which
involve many sub-cases). However, comparing these two trends of
thought, it seems to remain still really paradoxical that both
phenomena contribute to the reflection principle, without an
appropriate understanding of the general links between these two
techniques.  Guided by this observation, we have devised a new
two-sided technique.  In this article, we shall indeed perform the
proof of Theorem~1.9 by {\it combining the technique of the reflection
principle together with the consideration of envelopes of
holomorphy}. Further, we have been guided by a deep analogy between
the various reflection principles and the results on {\it propagation
of analyticity for CR functions along CR curves}, in the spirit of the
Russian school in the sixties, of Treves' school, in the spirit of the
works of Tr\'epreau, of Tumanov, of J\"oricke, of Porten and others:
the vector fields of the complex tangent bundle $T^cM$ being the {\it
directions} of propagation for the one-sided holomorphic extension of
CR functions, and the Segre varieties giving these directions (because
$T_q^cM=T_qS_{\bar q}$ for all $q\in M$), one can expect that Segre varieties
{\it also propagate the analyticity of CR mappings}.  Of
course, such a propagation property is already well known and
intensively studied since the historical works of Pinchuk [P1,2,3,4]
and since the important more recent articles of
Diederich-Forn{\ae}ss-Ye [DFY] and of Diederich-Pinchuk
[DP1,2]. However, in the classical works, one propagates along a
single Segre variety $S_{\bar p}$ and perhaps afterwards along the
subsequent ``Segre sets'' if necessary ([BER1,2], [Me4,5,6,7,8],
[Mi3,4]). As argued in [Me4,5], this terminology ``Segre sets'' is not 
the best one.
But in the present article we will propagate the analytic
properties along a {\it bundle} of Segre varieties of $M$, namely {\it
along a Levi-flat union of Segre varieties} $\Sigma_\gamma:=\cup_{q\in
\gamma} S_{\bar q}$, parametrized by a smooth curve $\gamma$
transversal to $T^cM$, in total analogy with the propagation of
analyticity of CR functions, where one uses a {\it bundle of attached
analytic discs}, parametrized by a curve transversal to $T^cM$ ({\it
cf.}~Tumanov's version of propagation [Tu2]; in this concern, we would
like to mention that recently, Porten [Po] has discovered a simple
strategy of proof using only CR orbits, deformations of bundles of
analytic discs and Levi forms on manifolds with boundary which treats
in an unified way the local ([Tu1]) and the global ([Tr2], [Tu2],
[Me1], [J]) wedge extension theorem). Let us now explain our strategy
in full details and describe our proof. To avoid excessive
technicalities in this presentation, we shall discuss the proof of
Theorem~1.2 instead of Theorem~1.9.

\subsection*{2.2.~Description of the proof of Theorem~1.2}
To begin with, recall from \S1 that the generic rank of the locally
defined holomorphic mapping $\mathcal{Q}_\infty'\,: t'\mapsto
(\Theta_\beta'(t'))_{\beta\in\N^{n-1}}$ is equal to the integer
$\chi_{M'}'$. The generic rank of an infinite collection
of holomorphic functions can always be interpreted in terms of finite
subcollections $\mathcal{Q}_k'(t')=
(\Theta_\beta'(t'))_{\vert\beta\vert\leq k}$. Of course, using the CR
diffeomorphism assumption, we may prove carefully that
$\chi_M=\chi_{M'}'$ ({\it see} Lemma~4.3).  It is known that $M'$ is
holomorphically nondegenerate if and only if $\chi_{M'}'=n$.  In the
remainder of \S2, we shall assume that $M'$ is holomorphically
nondegenerate. Let $q'\in M'$ be a point where the rank of
$\mathcal{Q}_k'(t')$ is equal to
$n$, hence locally constant. In our first step, we
will show that $h$ is real analytic at the reciprocal image of each
such point $h^{-1}(q')\in M$.  In fact, these points $q'$ are the {\it
finitely nondegenerate}\, points of $M'$, in the sense of [BER2,~\S11.2].  In
this case, it will appear that our proof of the first step is a
reminiscence of the Lewy-Pinchuk reflection principle and in fact, it
is a mild easy generalization of it, just by differentiating more than one
time. Afterwards, during the second (crucial and much more delicate)
step, to which \S5--8 below are devoted, we shall extend $h$ at
each point $h^{-1}(q')$, where $q'$ belongs to the real analytic subset
$E_{M'}'\subset M'$ where the mapping $\mathcal{Q}_\infty'$ is not of
rank $n$.  This is where we use envelopes of
holomorphy.  We shall start as follows.  By \S3.47, there exists a
proper real analytic subset $E_{M'}'$ of $M'$ such that the rank of
the mapping $\mathcal{Q}_\infty'$ localized around points $p'\in M'$
equals $n$ at each point $q'$ close to $p'$ not belonging to
$E_{M'}'$. Let $E_{\text{\rm na}}'\subset E_{M'}'\subset M'$
(``$\text{\rm na}$'' for ``non-analytic'') denote the closed set of
points $q'\in M'$ such that $h$ is {\it not}\, real 
analytic in a neighborhood of
$h^{-1}(q')$. By the first step, $E_{\rm na}'$ is necessarily
contained in $E_{M'}'$.  If $E_{\text{\rm na}}'=\emptyset$,
Theorem~1.9 would be proved, gratuitously. We shall therefore assume
that $E_{\text{\rm na}}'\neq \emptyset$ and we shall endeavour to
derive a contradiction in several nontrivial steps as
follows. Assuming that $E_{\rm na}'$ is nonempty, in order to come to
an absurd, it suffices to exhibit at least one point $p'$ of $E_{\rm
na}'$ such that $h$ is in fact real analytic in a
neighborhood of $h^{-1}(p')$. This is what we
shall achieve and the proof is long.  In analogy with what is done in
[MP1,2], we shall first show that we can choose a particular point
$p_1'\in E_{\text{\rm na}}'$ which is nicely disposed as follows ({\it
see} {\sc Figure~1}).

\def\thelemma{2.3}\begin{lemma}
\text{\rm ({\it cf.} [MP1, Lemma~2.3])} Let $E'\subset M'$ be an
arbitrary closed subset of an everyhere locally minimal real analytic
hypersurface $M'\subset \C^n$, with $n\geq 2$. If $E'$ and
$M'\backslash E'$ are nonempty, then there exists a point $p_1'\in E'$
and a real analytic one-codimensional submanifold $M_1'$ of $M'$ with
$p_1'\in M_1'\subset M'$ which is \text{\rm generic} in $\C^n$ and
which divides $M'$ near $p_1'$ in two open parts ${M_1'}^-$ and
${M_1'}^+$ such that $E'\backslash\{p_1'\}$ is contained in the open side
${M_1'}^+$ near $p_1'$.
\end{lemma}

\noindent
To reach the desired contradiction, it will suffice to prove that $h$
is analytic at the point $h^{-1}(p_1')$, where $p_1'\in E_{\text{\rm
na}}'\cap M_1'$ is such a special 
point as in Lemma~2.3 above. To this aim, we
shall pick a long embedded real analytic arc $\gamma'$ contained in
${M_1'}^-$ transverse to the complex tangential directions of $M'$,
with the ``center'' $q_1'$ of $\gamma'$ very close to $p_1'$ ({\it
see} {\sc Figure~1}). Next, using the inverse mapping $h^{-1}$, we can
copy back these objects on $M$, namely we set $E_{\text{\rm
na}}:=h^{-1}(E_{\text{\rm na}}')$, $\gamma:=h^{-1}(\gamma')$,
$p_1:=h^{-1}(p_1')$, $q_1:=h^{-1}(q_1')$, whence $M_1:=h^{-1}(M_1')$,
$M_1^-=h^{-1}({M_1'}^-)$ and $M_1^+=h^{-1}({M_1'}^+)$.

\bigskip
\begin{center}
\input figure6.pstex_t
\end{center}
\bigskip

To the analytic arc $\gamma'$, we shall associate holomorphic coordinates
$t'=(z',w')\in\C^{n-1}\times \C$, $w'=u'+iv'$, such that $p_1'=0$ and
$\gamma'$ is the $u'$-axis (in particular, some ``normal'' coordinates
in the sense of [BJT] would be appropriate, but not indispensable) and
we shall consider the reflection function $\mathcal{R}_h'(t,\bar\nu')=
\bar\mu'-\sum_{\beta\in \N^{n-1}} \bar{\lambda'}^\beta \
\Theta_\beta'(h(t))$ {\it in these coordinates $(z',w')$}. The
functions $\Theta_\beta'(h(t))$ will be called the {\it components of
the reflection function $\mathcal{R}_h'$}. Next, we choose coordinates
$t\in \C^n$ near $(M,p_1)$ vanishing at $p_1$. To the
$\mathcal{C}^\infty$-smooth arc $\gamma$, we shall associate the
following $\mathcal{C}^\infty$-smooth Levi-flat hypersurface:
$\Sigma_\gamma:= \bigcup_{q\in \gamma} S_{\bar q}$, where $S_{\bar q}$
denotes the Segre variety of $M$ associated to various points $q\in M$
({\it see} {\sc Figure~2}). Let $\Delta_n(0,\rho):=\{t\in \C^n: \vert
t\vert < \rho\}$ be the polydisc with center $0$ of polyradius
$(\rho,\dots,\rho)$, where $\rho>0$. Using the tangential Cauchy-Riemann
operators to differentiate the fundamental identity which reflects the
assumption $h(M)\subset M'$, we shall establish the following crucial
observation.

\def\thelemma{2.4}\begin{lemma}
There exists a positive real number $\rho>0$ independent of
$\gamma'$ such that all the components $\Theta_\beta'(h(t))$
of the reflection function extend as CR functions
of class $\mathcal{C}^\infty$ over $\Sigma_\gamma\cap
\Delta_n(0,\rho)$.
\end{lemma}

\noindent 
Furthermore, by global minimality of $M$, there exists a global
one-sided neighborhood $D$ of $M$ to which all CR functions (hence the
components of $h$) extend holomorphically ({\it see}\, the details in
\S3.6). We now recall that, by construction of
$M_1'$, the CR mapping $h$ is already holomorphic in a small
neighborhood of $h^{-1}(q')$ for every point $q'\in {M_1'}^-$.  It
follows that the components $\Theta_\beta'(h(t))$ of the reflection
function are already holomorphic in a fixed neighborhood, say
$\Omega$, of $M_1^-$ in $\C^n$. Also, they are already holomorphic at
each point of the global one-sided neighborhood $D$. In particular,
they are holomorphic in a neighborhood $\omega_\gamma\subset \Omega$
in $\C^n$ of $\gamma\subset M_1^-$. Then according to the
Hanges-Treves extension theorem [HaTr], we deduce that all the
components $\Theta_\beta'(h(t))$ of the reflection function extend
holomorphically to a neighborhood $\omega(\Sigma_\gamma)$ of
$\Sigma_\gamma$ in $\C^n$, which is a (very thin) neighborhood whose
size depends of course on the size of $\omega_\gamma$ (and the size of
$\omega_\gamma$ goes to zero without any explicit control as the
center point $q_1$ of $\gamma$ tends towards $p_1\in E_{\text{\rm na}}$).

\bigskip
\begin{center}
\input figure3.pstex_t
\end{center}
\bigskip

To achieve the final step, we shall consider the envelope of
holomorphy of $D\cup \Omega \cup \omega(\Sigma_\gamma)$ (in fact, to
prevent from poly-dromy phenomena, we shall instead consider a certain
subdomain of $D\cup \Omega\cup \omega(\Sigma_\gamma)$, see the details
in \S6 below), which is a kind of round domain $D\cup \Omega$ covered
by a thin Levi-flat almost horizontal ``hat-domain''
$\omega(\Sigma_\gamma)$ touching the ``top of the head'' $M$ along the
one-dimensional arc~$\gamma$ ({\it see} {\sc Figure~3}).

Our purpose will be to show that, if the
arc~$\gamma'$ is sufficiently close to $M_1'$ (whence $\gamma$ is also
very close to $M_1$), then the envelope of holomorphy of $D\cup \Omega
\cup \omega(\Sigma_\gamma)$ contains the point $p_1$, {\it
even if $\omega(\Sigma_\gamma)$ is arbitrarily thin}. We will
therefore deduce that all the components of the reflection function
extend holomorphically at $p_1$, thereby deriving the desired
contradiction. By exhibiting a special curved Hartogs domain, we shall
in fact prove that holomorphic functions in $D\cup \Omega \cup
\omega(\Sigma_\gamma)$ extend holomorphically to the lower one sided
neighborhood $\Sigma_\gamma^-$ (the ``same'' side as
$D=M^-$, {\it see}\,
{\sc Figure}~3); we explain below why this analysis gives analyticity at
$p_1$, even in the (in fact simpler) case where $p_1$ belongs to the other
side $\Sigma_\gamma^+$. 

\bigskip
\begin{center}
\input figure4.pstex_t
\end{center}
\bigskip

Notice that, because the order of contact
between $\Sigma_\gamma$ and $M$ is at least equal to two (because
$T_qM=T_q\Sigma_\gamma$ for every point $q\in \gamma$), we cannot
apply directly any version of the edge of the wedge theorem to this
situation. Another possibility (which, on the contrary, might well succeed)
would be to apply repeatedly the Hanges-Treves theorem, in the disc
version given in [Tu2] ({\it see} also [MP1]) to deduce that
holomorphic functions in $D\cup \Omega \cup \omega(\Sigma_\gamma)$
extend holomorphically to the lower side $\Sigma_\gamma^-$, just by
sinking progressively $\Sigma_\gamma$ into $D$. But this would require
a too complicated analysis for the desired statement.  Instead, by
performing what seems to be the simplest strategy, we shall use some
deformations (``translations'') of the following half analytic disc
attached to $\Sigma_\gamma$ along $\gamma$. We shall consider the
inverse image by $h$ of the half-disc $(\gamma')^c\cap D'$ obtained by
complexifying $\gamma'$ ({\it see} {\sc Figure~2} and
{\sc Figure~3}). Rounding off the corners and reparametrizing
the disc, we get an analytic disc $A\in \mathcal{O}(\Delta)\cap
\mathcal{C}^\infty(\overline{\Delta})$ with $A(b^+\Delta)\subset 
\gamma\subset \Sigma_\gamma$, where $b^+\Delta:= b\Delta \cap \{\text{\rm 
Re} \, \zeta \geq 0\}$, $b\Delta=\{\vert z\vert =1\}$ and $A(1)=q_1$. It is
this half-attached disc that we shall ``translate'' along the complex
tangential directions to~$\Sigma_\gamma$ as follows.

\def\thelemma{2.5}\begin{lemma}
There exists a $\mathcal{C}^\infty$-smooth 
$(2n-2)$-parameter family of analytic discs $A_\sigma: 
\Delta\to \C^n$, $\sigma\in\R^{2n-2}$, 
$\vert \sigma\vert <\varepsilon$, satisfying
\begin{itemize}
\item[{\bf (1)}]
The disc $A_\sigma\vert_{\sigma=0}$ coincides with the above disc $A$.
\item[{\bf (2)}]
The discs $A_\sigma$ are half-attached to $\Sigma_\gamma$, namely
$A_\sigma(b^+\Delta)\subset \Sigma_\gamma$.
\item[{\bf (3)}]
The boundaries $A_\sigma(b\Delta)$ of the discs $A_\sigma$
are contained in $D\cup \Omega \cup \omega(\Sigma_\gamma)$.
\item[{\bf (4)}]
The mapping $(\zeta,\sigma)\mapsto A_\sigma(\zeta)\in \Sigma_\gamma$ is a 
$\mathcal{C}^\infty$-smooth diffeomorphism from a neighborhood
of $(1,0)\in b\Delta\times\R^{2n-2}$ onto a neighborhood of
$q_1$ in $\Sigma_\gamma$.
\item[{\bf (5)}]
As $\gamma=h^{-1}(\gamma')$ varies and as $q_1$ tends to $p_1$, these 
discs depend $\mathcal{C}^\infty$-smoothly upon $\gamma'$ and properties
{\bf (1-4)} are stable under perturbations of $\gamma'$.
\item[{\bf (6)}]
If $\gamma(0)=q_1$ is sufficiently close to $M_1$, and if $p_1\in
\Sigma_\gamma^-$ is under $\Sigma_\gamma$
$($as in {\sc Figure~3}$)$, then the envelope of
holomorphy of $($an appropriate subdomain of$)$ 
$D\cup \Omega \cup \omega(\Sigma_\gamma)$ contains $p_1$.
\end{itemize}
\end{lemma}

Consequently, using these properties {\bf (1-6)} and applying the
continuity principle to the family $A_\sigma$, we shall obtain that
the envelope of holomorphy of $D\cup \Omega\cup \omega(\Sigma_\gamma)$
(in fact of a good subdomain of it, in order to assure monodromy)
contains a large part of the side $\Sigma_\gamma^-$ of $\Sigma_\gamma$
in which $D (=:M^-)$ lies. In the case where $p_1$ lies in this side
$\Sigma_\gamma^-$, and provided that the center point $q_1$ of
$\gamma$ is sufficiently close to $p_1$, we are done: the components
of the reflection function extend holomorphically at $p_1$ (this case
is drawn in {\sc Figure~3}). Of course,
it can happen that $p_1$ lies in the other side $\Sigma_\gamma^+$ or in
$\Sigma_\gamma$ itself. In fact, the following tri-chotomy is in order to
treat the problem. To apply Lemma~2.5 correctly, and to complete the study
of our situation, we shall indeed distinguish three cases.  
\begin{itemize}
\item[{\bf Case~I.}] \ 
The Segre variety $S_{\bar p_1}$ cuts $M_1^-$ along an
infinite sequence of points $(q_k)_{k\in\N}$ tending towards $p_1$.
\smallskip
\item[{\bf Case~II.}] \ 
The Segre variety $S_{\bar p_1}$ does not intersect
$M_1^-$ in a neighborhood of $p_1$ and it goes under $M_1^-$, 
namely inside $D$.
\smallskip
\item[{\bf Case~III.}] \ 
The Segre variety $S_{\bar p_1}$ does not intersect
$M_1^-$ in a neighborhood of $p_1$ and it goes over $M_1^-$, namely 
over $D\cup M_1^-$.
\end{itemize}
In the first case, choosing the point $q_1$ above to be one of the
points $q_k$ which is sufficiently close to $p_1$, and using the fact
that $p_1$ {\it belongs to} $S_{\bar q_1}$ (because $q_1\in S_{\bar
p_1}$), we have in this case $p_1\in \Sigma_\gamma$ and the
holomorphic extension to a neighborhood $\omega(\Sigma_\gamma)$
already yields analyticity at $p_1$ (in this case, we have
nevertheless to use Lemma~2.5 to insure monodromy of the
extension). In the second case, we have $S_{\bar p_1}\cap
D\neq\emptyset$. We then choose the center point $q_1$ of $\gamma$
very close to $p_1$. Because we have in this case a uniform
control of the size of $\omega(\Sigma_\gamma)$, we again get that
$p_1$ always belongs to $\omega(\Sigma_\gamma)$ and Lemma~2.5 is again used
to insure monodromy. In the third (a priori more delicate) case, by a
simple calculation, we shall observe that $p_1$ always belong to the
lower side $\Sigma_\gamma^-$ (as in {\sc Figure~3}) and Lemma~2.5
applies to yield holomorphic extension and monodromy of the extension.
In sum, we are done in all the three cases: we have shown that the
components $\Theta_\beta'(h(t))$ all extend holomorphically at
$p_1$. Finally, using a complex analytic set similar to
$\mathcal{C}_h'$ defined in~(1.12) and Lemma~4.14 below, we deduce
that $h$ is real analytic at $p_1$.

In conclusion to this presentation, we would like to say that some
unavoidable technicalities that we have not mentioned here will render
the proof a little bit more complicated (especially about the choice
of $q_1$ sufficiently close to $p_1$, about the choice of $\gamma$ and
about the smooth dependence with respect to $\gamma$ of
$\Sigma_\gamma$ and of $A_\sigma$). The remainder of the paper is
devoted to complete these technical features thoroughly. At first, we
provide some necessary background material about the reflection function.

\section*{\S3.~Biholomorphic invariance of the reflection function}

\subsection*{3.1.~Preliminary and notation} Let $p'\in M'$, let
$t'=(z_1',\dots,z_{n-1}',w')=(z',w')$ be holomorphic coordinates
vanishing at $p'$ such that the projection $T_{p'}^cM'\to
\C_{z'}^{n-1}$ is submersive. As in \S1, we can represent $M'$ by a
complex analytic defining equation of the form $\bar w'=\Theta'(\bar
z',t')$, where the right hand side function converges normally in the
polydisc $\Delta_{2n-1}(0,\rho')$ for some $\rho'>0$. Here, by normal
convergence we mean precisely that there exists a constant $C>0$ such
that if we develope $\Theta'(\bar
z',t')=\sum_{\beta\in\N^{n-1}}\sum_{\alpha\in\N^n}\, (\bar z')^\beta
(t')^\alpha\,\Theta_{\beta,\alpha}'$, with $\Theta_{\beta,
\alpha}'\in\C$, then we have
\def\theequation{3.2}\begin{equation}
\vert\Theta_{\beta,\alpha}'\vert \leq
C\, (\rho')^{-\vert\alpha\vert-\vert\beta\vert},
\end{equation} 
for all multi-indices $\alpha$ and $\beta$.  Furthermore, by the
reality of $M'$ the function $\Theta'$ satisfies the power series
identity $\Theta'(\bar z',z',\overline{\Theta}'(z',\bar z',\bar
w'))\equiv \bar w'$. It follows from this identity that
$\Theta_0'(t')$ does not vanish identically, and in fact contains the
monomial $w'\equiv \Theta_0'(0,w')$. 
We set $p:=h^{-1}(p')$ and similarly, we represent a
local defining equation of $M$ near $p$ as $\bar w=\Theta(\bar z,t)$,
where $\Theta$ converges normally in $\Delta_{2n-1}(0,\rho)$ for some
$\rho>0$. We denote the mapping by $h:=(f,g):=
(f_1,\dots,f_{n-1},g)$. Then the assumption that $h$ maps $M$
into $M'$ yields that
\def\theequation{3.3}\begin{equation} 
\overline{g(t)}=\Theta'(\overline{f(t)},h(t)),
\end{equation}
for all $t\in M$ near $p$. For this relation to hold locally, it is convenient
to assume that $\vert h(t)\vert < \rho'$ for every $t\in M$ with
$\vert t \vert < \rho$.  

Since by assumption the $h_j$ are of class $\mathcal{C}^\infty$ and CR
over $M$, we can extend them to a neighborhood of $M$ in $\C^n$ as
functions $\widetilde{h}_j$ of class $\mathcal{C}^\infty$ with
antiholomorphic derivatives
$\partial_{\bar t_l}\widetilde{h}_j$ vanishing to infinite order on
$M$, $l=1,\dots,n$.  
So, if we develope these extensions in real Taylor series at
each point $q\in M$ as follows: 
\def\theequation{3.4}\begin{equation}
T_q^\infty \widetilde{h}_j=
\widetilde{h}_j(q)+\sum_{\alpha\in\N^n\backslash\{0\}}\,\partial_t^\alpha
\widetilde{h}_j(q)\,
(t-q)^\alpha/\alpha\,! \ \in \ \C\dl t\dr,
\end{equation}
there are no antiholomorphic term.  

The reflection function associated
with such a coordinate system and with such a defining equation, 
namely 
\def\theequation{3.5}\begin{equation}
\mathcal{R}_h'(t,\bar\nu'):=\bar\mu'-\sum_{\beta\in\N^{n-1}}\,
(\bar\lambda')^\beta\,\Theta_\beta'(h(t)),
\end{equation}
where $\bar\nu'=(\bar\lambda',\bar\mu')$,
converges normally with respect to $t\in M$ with $\vert t\vert < \rho$
and $\bar\nu'\in\C^n$ with $\vert \bar\nu'\vert < \rho'$, hence
defines a function which is CR of class $\mathcal{C}^\infty$ on $M$ near
$p$ and holomorphic with respect to $\bar\nu'$. The main goal of this
paragraph is to study its invariance with respect to changes of
coordinates. 

\subsection*{3.6.~Holomorphic extension to a one-sided neighborhood attached
to $M$} Before treating invariance, recall that thanks to the local
minimality at every point, all CR functions on $M$ and in particular
the $h_j$ extend holomorphically to one side of $M$ at every point of
$M$ (the simplest proof of this result can be found in [R]; {\it
see}\, also the excellent survey [Tr3] for a proof using Bishop
discs). Of course, the side may vary. We do not require that $M$ be
orientable, but anyway the small pieces $(M,p)$ always divide locally
$\C^n$ in two components $(M,p)^\pm$. By shrinking these one-sided
neighborhoods covered by attached analytic discs, we may assume that
for every point $p\in M$, all CR functions on $M$ extend
holomorphically to the intersection of a small nonempty open ball
$B_p$ centered at $p$ with one of the two local open components
$(M,p)^\pm$.  Let $D_p$ denote the resulting open side of $M$ at $p$,
namely $D_p=B_p\cap (M,p)^+$ or $D_p=B_p\cap (M,p)^-$. Since the union
of the various open sets $D_p$ does not necessarily make a domain, we
introduce the following definition. By a {\it global one-sided
neighborhood}\, of $M$ in $\C^n$, we mean a {\it domain}\, $D$ such
that for every point $p\in M$, $D$ contains a local one-sided
neighborhood of $M$ at $p$. In particular, $D$ necessarily contains a
neighborhood of a point $q\in M$ if it contains the two local sides of
$M$ at $q$. To construct a global one-sided neighborhood to which all
$\mathcal{C}^\infty$-smooth and even $\mathcal{C}^0$-smooth CR
functions on $M$ extend holomorphically, it suffices to set
\def\theequation{3.7}\begin{equation}
D:=\bigcup_{q\in M} \, D_q \ 
\bigcup_{D_p\cap D_q=\emptyset} \
\left(
B_p\cap \overline{D_p}\cap 
\overline{D_q}\cap B_q
\right).
\end{equation}
The second part of this union consists of an open subset of $M$ which
connects every meeting pair of local one-sided neighborhoods in the
case where their respective sides differ. If the radii of the $B_p$ are
sufficiently small compared to the geometric distortion of $M$, then
the open set defined by~(3.6) is a domain in $\C^n$. Moreover, using
the uniqueness principle for CR functions, it is elementary to see
that every CR function $\phi$ on $M$ extends as a unique
holomorphic function globally defined over $D$. In this concern, we
would like to mention that a more general construction in arbitrary
codimension in terms of attached wedges is provided in [Me2],
[MP1,2] and in [Da2]. 

Since each $D_p$ is contained in some union of small Bishop discs with
boundaries contained in $(M,p)$, it follows that the maximum modulus
of the holomorphic extension of $\phi$ to $D_p$ is less than or equal
to the maximum modulus of the CR function $\phi$ over the piece
$(M,p)$, which is a little bit larger than $\overline{D_p}\cap M$. To
be precise, after shrinking $B_p$ if necessary, we can assume that the
Bishop discs covering $D_p$ have their boundaries attached to $M\cap
\{\vert t\vert < \rho\}$.  Since $\vert h(t)\vert < \rho'$ for $t\in
M$ with $\vert t \vert < \rho$, the same majoration holds for $t\in
D_p$ (maximum principle), so it follows that the series defined
by~(3.5) also converges normally with respect to $t$ inside $D_p$. In
conclusion, we have established the following.

\def\thelemma{3.8}\begin{lemma}
With the above notation, $\mathcal{R}_h'$ is
defined in the set
\def\theequation{3.9}\begin{equation}
\left[D_p\cup (M\cap \Delta_n(0,\rho))\right]\times
\Delta_n(0,\rho') \ \subset \ 
\Delta_n(0,\rho)\times\Delta_n(0,\rho').
\end{equation}
Precisely, $\mathcal{R}_h'$ is holomorphic with respect to
$(t,\bar\nu')$ in $D_p\times \Delta_n(0,\rho')$ and it is CR
of class $\mathcal{C}^\infty$ over the real analytic hypersurface
\def\theequation{3.10}\begin{equation}
\left[M\cap \Delta_n(0,\rho)\right]\times \Delta_n(0,\rho') \ \subset \
\Delta_n(0,\rho)\times \Delta_n(0,\rho').
\end{equation}
\end{lemma}

\subsection*{3.11.~Characterization of the holomorphic extendability
of $\mathcal{R}_h'$} Let $x\in\C^m$, $x'\in \C^{m'}$ and consider a
power series of the form 
\def\theequation{3.12}\begin{equation}
R(x,x'):=\sum_{\alpha\in\N^m,
\alpha'\in\N^{m'}} \,R_{\alpha,\alpha'}\, x^\alpha\, (x')^{\alpha'},
\end{equation}
where the $R_{\alpha,\alpha'}$ are complex coefficients. Let us assume
that $R$ converges normally in some polydisc
$\Delta_m(0,\sigma)\times\Delta_{m'}(0,\sigma')$, for some two
$\sigma, \sigma'>0$. By normal convergence, we mean that there exists
a constant $C>0$ such that the Cauchy inequalities $\vert
R_{\alpha,\alpha'}\vert \leq C\, (\sigma)^{-\vert\alpha\vert}\,
(\sigma')^{-\vert\alpha'\vert}$ hold. Let us define $R_{\alpha'}(x):=
\sum_{\alpha\in\N^m}\, R_{\alpha,\alpha'}\, x^\alpha=
\left[
{1\over \alpha'!}\,\partial_{x'}^{\alpha'}R(x,x')
\right]_{x':=0}$.  Classically in
the basic theory of converging power series, it follows 
that for every positive $\widetilde{\sigma}<\sigma$, there
exists a constant $C_{\widetilde{\sigma}}$ which depends on
$\widetilde{\sigma}$ such that for all $x$ satisfying $\vert x\vert <
\widetilde{\sigma}$, the estimate $\vert R_{\alpha'}(x)\vert\leq
C_{\widetilde{\sigma}}\, (\sigma')^{-\vert\alpha'\vert}$ holds. Indeed,
we simply compute for $\vert x\vert < \widetilde{\sigma}$ the
elementary series:
\def\theequation{3.13}\begin{equation}
\left\{
\aligned
{} &
\vert R_{\alpha'}(x)\vert \leq \sum_{\alpha\in\N^m}\,
\vert R_{\alpha,\alpha'}\vert\,\vert x\vert^\alpha\leq\\
& \ \ \ \ \ \ \ \ \ \ \ \ \ 
\leq C \ \sum_{\alpha\in\N^m}\,
\sigma^{-\vert\alpha\vert}\,(\sigma')^{-\vert\alpha'\vert}\, 
\widetilde{\sigma}^{\vert\alpha\vert}=
C\, \left(
{\sigma\over \sigma-\widetilde{\sigma}}
\right)^m\, 
(\sigma')^{-\vert\alpha'\vert}.
\endaligned\right.
\end{equation}
As an application, such an inequality applies to the defining function of $M'$:
for every positive $\widetilde{\rho'}<\rho'$, there exists a constant
$C_{\widetilde{\rho'}}$ such that for all $\vert t'\vert < 
\widetilde{\rho'}$ we have
\def\theequation{3.14}\begin{equation}
\vert\Theta_\beta'(t')\vert \leq
C_{\widetilde{\rho'}}\
(\rho')^{-\vert\beta\vert}.
\end{equation}
The estimation~(3.13) also exhibits an interesting basic
property. Suppose for a while that the reflection function
$\mathcal{R}_h'$ defined by~(3.5) extends holomorphically to the
polydisc $\Delta_n(0,\sigma)\times \Delta_n(0,\sigma')$ for some
positive $\sigma, \sigma'>0$ with $\sigma<\rho$ and
$\sigma'<\rho'$. Then the functions $\theta_\beta'(t)$ defined by
\def\theequation{3.15}\begin{equation}
\theta_\beta'(t):=\left[{1\over \beta!}\, \partial_{\bar\lambda'}^\beta
\mathcal{R}_h'(t,\bar\nu')
\right]_{\bar\lambda':=0}
\end{equation}
satisfy a Cauchy estimate, namely $\vert \theta_\beta'(t)\vert \leq
C_{\widetilde{\sigma}} \,(\sigma')^{-\vert\beta\vert}$ for all $\vert
t\vert < \widetilde{\sigma}<\sigma$.  By~(3.5), notice that
$\theta_\beta'(t)\equiv \Theta_\beta'(h(t))$ over $M\cap
\Delta_n(0,\rho)$ and inside $D_p$, so the holomorphic
extendability of $\mathcal{R}_h'$ implies that all
the components $\Theta_\beta'(h(t))$ extend holomorphically 
to $\Delta_n(0,\sigma)$.  These preliminary observations
are appropriate to obtain the following useful characterization of the
holomorphic extendability of $\mathcal{R}_h'$ which says in substance
that it suffices that all its components $\Theta_\beta'(h(t))$ extend
at $p$ and then afterwards the Cauchy estimate holds automatically.

\def\thelemma{3.16}\begin{lemma}
The following three properties are equivalent:
\begin{itemize}
\item[{\bf (i)}]
There exists $\sigma>0$ with $\sigma<\rho$ and $\sigma<\rho'$ such
that $\mathcal{R}_h'$ extends holomorphically to the polydisc
$\Delta_n(0,\sigma)\times \Delta_n(0,\sigma)$.
\item[{\bf (ii)}]
There exists $\sigma>0$ with $\sigma<\rho$ such that all
$\mathcal{C}^\infty$-smooth CR functions $\Theta_\beta'(h(t))$ defined
on $M\cap \Delta_n(0,\rho)$ extend holomorphically to the polydisc
$\Delta_n(0,\sigma)$ as holomorphic 
functions $\theta_\beta'(t)$ which satisfy the
inequality $\vert \theta_\beta'(t)\vert \leq C
\,(\sigma')^{-\vert\beta\vert}$ for some two positive constants $C>0$,
$\sigma'<\rho'$ and for all $\vert t\vert <\sigma$.
\item[{\bf (iii)}]
There exists $\sigma>0$ with $\sigma<\rho$ such that all
$\mathcal{C}^\infty$-smooth CR functions $\Theta_\beta'(h(t))$ defined
on $M\cap \Delta_n(0,\rho)$ extend holomorphically to the polydisc
$\Delta_n(0,\sigma)$ as holomorphic functions $\theta_\beta'(t)$.
\end{itemize}
\end{lemma}

\proof
Of course, {\bf (i)} implies {\bf (ii)} which in turn implies {\bf
(iii)} trivially. Conversely, let us show that {\bf (iii)} implies
{\bf (ii)}. By~(3.4) with $q=0$, the Taylor series of $h_j$ at the
origin $H_j(t):=T_0^\infty \widetilde{h}_j(t)$ involves only
holomorphic monomials $t^\alpha$ and {\it no}\, antiholomorphic
monomial. We notice that the Taylor series at the origin of
$\Theta_\beta'(h(t))$ coincides with the composition of formal power
series $\Theta_\beta'(H(t))$. Consequently, by the assumption {\bf
(iii)}, the formal power series mapping $H(t)$ is a formal solution of
some evident complex analytic equations. Indeed, we have
\def\theequation{3.17}\begin{equation}
R_\beta'(t,H(t)):=\Theta_\beta'(H(t))-\theta_\beta'(t)\equiv 0 \ \ \ \ \ 
{\rm in} \ \C\dl t\dr,
\end{equation}
for all $\beta\in\N^{n-1}$. By the Artin approximation 
theorem ({\it see}\, [Ar]), there exists an analytic power series 
$\widetilde{H}(t)$ with $\widetilde{H}(0)=0$,
which converges normally in some polydisc, say 
$\Delta_n(0,\sigma)$ with $\sigma>0$, and which satisfies
\def\theequation{3.18}\begin{equation}
R_\beta'(t,\widetilde{H}(t))
:=\Theta_\beta'(\widetilde{H}(t))-\theta_\beta'(t)\equiv 0,
\end{equation}
for all $t\in\Delta_n(0,\sigma)$.  Shrinking $\sigma$ if necessary, we
may assume that for $\vert t\vert < \sigma$, we have $\vert
\widetilde{H}(t)\vert < \sigma'<\rho'$. Then the Cauchy
estimate~(3.14) valuable for the $\Theta_\beta'(t')$ yields by
composition a Cauchy estimate for $\Theta_\beta'(\widetilde{H}(t))$
which in turn yields the desired Cauchy estimate for the $\theta_\beta'(t)$
as stated in the end of {\bf (ii)}, thanks to the
relations~(3.18). This completes the proof.
\endproof

\subsection*{3.19.~Invariance of the reflection function} Our definition
of the reflection function $\mathcal{R}_h'$ seems to be
unsatisfactory, because it heavily depends on the choice of
coordinates and on the choice of a local defining function for
$(M',p')$. Our purpose is now to show that Theorem~1.9 holds true for
{\it every}\, system of coordinates provided it holds for {\it one}\,
such system. This requires to analyze how the components
$\Theta_\beta'(h(t))$ behave under the action of biholomorphisms.  Let
$t''=\Lambda(t')$ be a local biholomorphic mapping such that
$\Lambda(0)=0$, denote $t''=(z'',w'')=(z_1'',\dots,z_{n-1}'',w'')$
and denote $\Lambda=(\Phi_1,\dots,\Phi_{n-1},\Psi)$ accordingly. By
the implicit function theorem, if we assume that the linear mapping
$\pi''\circ d\Lambda: T_0^c M'\to \C_{z''}^{n-1}$ is bijective, where
$\pi'':\C_{z'',w''}^n\to \C_{z''}^{n-1}$ is the projection parallel to
the $w''$ axis, then the image $\Lambda(M')$ can also be defined locally
in a neighborhood of the origin by a defining equation of the
form $\bar w''= \Theta''(\bar z'',t'')$ similar to that
of $M'$. Equivalently, this differential geometric
condition can be expressed by the nonvanishing
\def\theequation{3.20}\begin{equation}
{\rm det}\, \left(
\overline{L}_j'\,\overline{\Phi}_k(0)\right)_{1\leq j,k\leq n-1}
\neq 0, 
\end{equation}
where the $\overline{L}_j'$ constitute a basis for the CR vector fields
on $M'$, namely $\overline{L}_j'=\partial_{\bar z_j'}+
\Theta_{\bar z_j'}'(\bar z',t')\,\partial_{\bar w'}$ for $j=1,\dots,n-1$.
Thus, we aim to compare the two reflection functions
\def\theequation{3.21}\begin{equation}
\left\{
\aligned
{} &
\mathcal{R}_h'(t,\bar\nu'):=\bar\mu'-\sum_{\beta\in\N^{n-1}}\,
(\bar\lambda')^\beta\,\Theta_\beta'(h(t)),\\
&
\mathcal{R}_{\Lambda\circ h}''(t,\bar\nu''):=
\bar\mu''-\sum_{\beta\in\N^{n-1}}\,(\bar\lambda'')^\beta\,
\Theta_\beta''(\Lambda\circ h(t)).
\endaligned\right.
\end{equation}
Without loss of generality, we can assume that $\Theta''$ converges
normally in $\Delta_{2n-1}(0,\rho'')$ and that
$\Lambda(\Delta_n(0,\rho'))$ is contained in $\Delta_n(0,\rho'')$. The
following lemma exhibits the desired invariance under biholomorphic
transformations fixing the center point $p'$ and Lemma~3.37 below will
show the invariance under local translations of the center point.

\def\thelemma{3.22}\begin{lemma}
The following two conditions are equivalent:
\begin{itemize}
\item[{\bf (i)}]
There exists $\sigma>0$ with
$\sigma<\rho$ and $\sigma<\rho'$ such that $\mathcal{R}_h'(t,\bar \nu')$
extends holomorphically to the polydisc
$\Delta_n(0,\sigma)\times\Delta_n(0,\sigma)$.
\item[{\bf (ii)}]
There exists $\sigma>0$ with $\sigma<\rho$ and $\sigma<\rho''$ such
that $\mathcal{R}_{\Lambda\circ h}''(t,\bar\nu'')$ extends
holomorphically to the polydisc
$\Delta_n(0,\sigma)\times\Delta_n(0,\sigma)$.
\end{itemize}
\end{lemma}

\proof
Of course, it suffices to prove that {\bf (i)} implies {\bf (ii)},
because $\Lambda$ is invertible. The proof is a little bit long and
calculatory, but the principle is quite simple (in advance, the reader
may skip to equation~(3.35) and to the paragraph following which explain
well the relation between the components of the two reflection functions).
As $\Lambda$ maps $M'$ into $M''$, there exists a converging
power series $A(t',\bar t')$ such that the following identity 
holds for all $t'$ with $\vert t'\vert < \rho'$:
\def\theequation{3.23}\begin{equation}
\overline{\Psi}(\bar t')-\Theta''(\overline{\Phi}(\bar t'),
\Lambda(t'))\equiv
A(t',\bar t')\,
\left[\bar w'-\Theta'(\bar z',t')\right]
\end{equation}
Replacing $\bar w'$ by $\Theta'(\bar z',t')$ on the left hand side, 
we get an interesting formal power series identity at the origin in
$\C^{2n-1}$
\def\theequation{3.24}\begin{equation}
\overline{\Psi}(\bar z',\Theta'(\bar z',t'))\equiv
\Theta''(\overline{\Phi}(\bar z',\Theta'(\bar z',t')),
\Lambda(t')),
\end{equation}
which converges for all $\vert \bar z'\vert <\rho'$ and $\vert t'\vert
< \rho'$. Putting $\bar z'=0$, we see first that 
\def\theequation{3.25}\begin{equation}
\overline{\Psi}(0,\Theta'(0,t'))\equiv \Theta''
(\overline{\Phi}(0,\Theta'(0,t')),\Lambda(t')).
\end{equation}
Next, we differentiate the relation~(3.24) with respect to $\bar z_j'$
for $j=1,\dots,n-1$. Remembering that $\overline{L}_j'=
\partial_{\bar z_j'}+\Theta_{\bar z_j'}'(\bar z',t')\,\partial_{\bar
w'}$, we see that differentiation with respect to $\bar z_j'$ is the
same as applying the operator $\overline{L}_j'$ and we get by the
chain rule
\def\theequation{3.26}\begin{equation}
\overline{L}_j'\overline{\Psi}(\bar z',\Theta'(\bar z',t'))\equiv
\sum_{k=1}^{n-1}\,\overline{L}_j'\overline{\Phi}_k(\bar z',
\Theta'(\bar z',t'))\,{\partial\Theta''\over\partial \bar z_k''}(
\overline{\Phi}(\bar z',\Theta'(\bar z',t')),\Lambda(t')).
\end{equation}
Consider the following determinant, which, by the assumption~(3.20) does
not vanish at the origin:
\def\theequation{3.27}\begin{equation}
\mathcal{D}'(\bar z',t'):=
{\rm det}\,\left(
\overline{L}_j'\overline{\Phi}_k(\bar z',\Theta'(\bar z',t'))
\right)_{1\leq j,k\leq n-1}.
\end{equation}
Shrinking $\rho'$ if necessary, we can assume that $\mathcal{D}'$
is nonzero at every point of $\Delta_{2n-1}(0,\rho')$. Then using
the rule of Cramer, we can solve in~(3.26) the first order
partial derivatives of $\Theta''$ with respect to the rest. 
We obtain an expression of the form
\def\theequation{3.28}\begin{equation}
{\partial\Theta''\over\partial \bar z_k''}(
\overline{\Phi}(\bar z',\Theta'(\bar z',t')),\Lambda(t'))\equiv
{R_k(\{(\overline{L}')^\gamma\overline{\Lambda}_i(
\bar z',\Theta'(\bar z',t'))\}_{\vert\gamma\vert=1,1\leq i\leq n})\over
\mathcal{D}'(\bar z',t')}.
\end{equation}
Here, for every multi-index $\gamma\in\N^{n-1}$, we denote by
$(\overline{L}')^\gamma$ the antiholomorphic derivation of order
$\vert\gamma\vert$ defined by $(\overline{L}_1')^{\gamma_1}\cdots
(\overline{L}_{n-1}')^{\gamma_{n-1}}$. Moreover, in~(3.28), it is a fact
that the terms $R_k$ are certain universal polynomials in their 
$n(n-1)$ arguments.

By differentiating again~(3.28) with respect to the $\bar z_j'$, 
using Cramer's rule, and making an inductive argument, it follows
that for every multi-index $\beta\in\N^{n-1}$, there exists a certain 
complicated but universal polynomial $R_\beta$ such that the
following relation holds:
\def\theequation{3.29}\begin{equation}
{1\over \beta!}\,
{\partial^{\vert\beta\vert}\Theta''\over
\partial(\bar z'')^\beta}(
\overline{\Phi}(\bar z',\Theta'(\bar z',t')),\Lambda(t'))\equiv
{R_\beta(\{(\overline{L}')^\gamma
\overline{\Lambda}_i(\bar z',\Theta'(\bar z',t'))
\}_{1\leq \vert\gamma\vert\leq
\vert\beta\vert,\,1\leq i\leq n})\over
[\mathcal{D}'(\bar z',t')]^{2\vert \beta\vert -1}}.
\end{equation}
Now, we put $\bar z':=0$ in these identities. An important observation
is in order. The composed derivations $(\overline{L}')^\gamma$ are
certain differential operators with nonconstant coefficients. Using
the explicit expression of the $\overline{L}_j'$, we see that all
these coefficients are certain universal polynomials of the collection
of partial derivatives $\{\partial^{\vert\delta\vert} \Theta' (\bar
z',t')/\partial (\bar z')^\delta\}_{1\leq \vert \delta\vert \leq
\vert\gamma\vert}$. Thus the numerator of~(3.29), after putting $\bar
z':=0$, becomes a certain holomorphic function of the collection
$\{\Theta_\gamma'(t')\}_{0\leq \vert\gamma\vert\leq\vert\beta\vert}$
(recall $\Theta_\gamma'(t')= \left[ {1\over
\gamma!}\,\partial^{\vert\gamma\vert} \Theta' (\bar z',t')/\partial
(\bar z')^\gamma\right]_{\bar z':=0}$). A similar property holds for
the denominator. In summary, we have shown that there exists
an infinite collection of holomorphic functions $S_\beta$ of their 
arguments such that 
\def\theequation{3.30}\begin{equation}
{1\over \beta!}\,
{\partial^{\vert\beta\vert}\Theta''\over
\partial(\bar z'')^\beta}(
\overline{\Phi}(0,\Theta'(0,t')),\Lambda(t'))\equiv
S_\beta(\{\Theta_\gamma'(t')\}_{\vert\gamma\vert\leq\vert\beta\vert})=:
s_\beta(t'),
\end{equation}
where the left and right hand sides are holomorphic functions of $t'$
running in the polydisc $\Delta_n(0,\rho')$. Furthermore, by Cauchy's
integral formula, there exists a positive constant $C$ such that
for all $\vert \bar z''\vert, \vert t''\vert<\rho''/2$, we have the majoration
\def\theequation{3.31}\begin{equation}
\left\vert
{1\over\beta!}\,
{\partial^{\vert\beta\vert}\Theta''\over
\partial(\bar z'')^\beta}(\bar z'',t'')
\right\vert\leq C\,(\rho''/2)^{-\vert\beta\vert}.
\end{equation}
Consequently we
get the estimate $\vert s_\beta(t')\vert \leq 
C\,(\rho''/2)^{-\vert\beta\vert}$. Now, 
let us rewrite the relations~(3.30) in a more explicit form, 
taking into account that $\Theta'(0,t')=\Theta_0'(t')$ by definition:
\def\theequation{3.32}\begin{equation}
\left\{
\aligned
{} &
\Theta_\beta''(\Lambda(t'))+\sum_{\gamma\in\N_*^{n-1}}\,
(\overline{\Phi}(0,\Theta_0'(t')))^\gamma\
{(\beta+\gamma)!\over\beta!\ \gamma!}\
\Theta_{\beta+\gamma}''(\Lambda(t'))\equiv\\
& \ \ \ \ \ \ \ \ \ \ \ \ \ \ \ \
\equiv S_\beta(\{\Theta_\delta'(t')\}_{\vert\delta\vert\leq\vert\beta\vert})=:
s_\beta(t'),
\endaligned\right.
\end{equation}
where we denote $\N_*^{n-1}:=\N^{n-1}\backslash\{0\}$.
This collection of equalities may be considered as an infinite upper
triangular linear system with unknowns being the functions
$\Theta_\beta''(\Lambda(t'))$.  This system can be readily inverted.
Indeed, using Taylor's formula in the convergent case or proceeding
directly at the formal level, it is easy to see that if we are given
an infinite collection of equalities with complex coefficients and
with $\zeta\in\N^{n-1}$ which is of the form
\def\theequation{3.33}\begin{equation}
\Theta_\beta''+\sum_{\gamma\in\N_*^{n-1}}\,
\zeta^\gamma \ {(\beta+\gamma)!\over\beta!\ \gamma!} \
\Theta_{\beta+\gamma}''=S_\beta,
\end{equation}
for all multi-indices $\beta\in\N^{n-1}$, then we can solve the 
unknowns $\Theta_\beta''$ in terms of the right hand side terms $S_\beta$
by means of a totally similar formula, except for signs:
\def\theequation{3.34}\begin{equation}
S_\beta+\sum_{\gamma\in\N^{n-1}_*}\,
\zeta^\gamma \ (-1)^\gamma \ {(\beta+\gamma)!\over\beta!\ \gamma!} \
S_{\beta+\gamma}=\Theta_\beta'',
\end{equation}
for all $\beta\in\N^{n-1}$. Applying this observation to~(3.32) and
using the above Cauchy estimates on $s_\beta(t')$, 
we deduce the convergent representation
\def\theequation{3.35}\begin{equation}
\left\{
\aligned
{} &
\Theta_\beta''(\Lambda(t'))\equiv 
S_\beta(\{\Theta_\delta'(t')\}_{\vert\delta\vert\leq\vert\beta\vert})+\\
& \ \ \ \ \ 
+\sum_{\gamma\in\N^{n-1}_*}\,
(\overline{\Phi}(0,\Theta_0'(t')))^\gamma \
(-1)^\gamma \
{(\beta+\gamma)!\over \beta! \ \gamma!} \
S_{\beta+\gamma}(\{\Theta_\delta'(t')\}_{\vert\delta\vert\leq
\vert\beta\vert+\vert\gamma\vert}),
\endaligned\right.
\end{equation}
which is valuable for $\vert t'\vert < \rho'$. Here, we recall that
the functions $S_\beta$ only depend on the biholomorphism $\Lambda$
and that they are holomorphic with respect to their arguments. Now, we
can prove that {\bf (i)} implies {\bf (ii)} in Lemma~3.22. By the
equivalence between {\bf (i)} and {\bf (ii)} of Lemma~3.16, it
suffices to show that all component functions $\Theta_\beta''(\Lambda
(h(t)))$ extend holomorphically to a neighborhood of the origin
provided all component functions $\Theta_\beta'(h(t))$ extend
holomorphically (by construction, the Cauchy estimates are already at
hand).  But this is evident by reading~(3.35) after replacing $t'$ by
$h(t)$.  This completes the proof of Lemma~3.22.
\endproof

\subsection*{3.36.~Translation of the center point} 
We have shown that the holomorphic extendability of the reflection
function $\mathcal{R}_h'$ centered at {\it one}\, point $p\times
\overline{h(p)}$ is an invariant property. On the other hand, suppose
that $\mathcal{R}_h'$ is holomorphic in the product polydisc
$\Delta_n(0,\sigma)\times \Delta_n(0,\sigma')$, for $0<\sigma<\rho$ and
$0<\sigma'<\rho'$. Does it follow that the reflection functions
centered at points $q\times \overline{h(q)}\in\Delta_n(0,\sigma)\times
\Delta_n(0,\sigma')$ also extends holomorphically at these points~?
Without loss of generality, we can assume that $h(M\cap
\Delta_n(0,\rho))\subset \Delta_n(0,\rho')$ and that $h(M\cap
\Delta_n(0,\sigma))\subset \Delta_n(0,\sigma')$. Let $q\in
\Delta_n(0,\sigma)$ be an arbitrary point and set $q':=h(q)$. Recall
that as in \S3.1 above, we are given coordinates $t$ and $t'$ centered
at the origin in which the equations of $M$ and of $M'$ are of the
form $\bar w=\Theta(\bar z,t)$ and $\bar w'=\Theta'(\bar z',t')$, with
$\Theta$ converging normally in the polydisc $\Delta_{2n-1}(0,\rho)$
and similarly for $\Theta'$.  We can center new holomorphic
coordinates at $q$ and at $q'$ simply by setting $t_*:=t-q$ and
$t_*':=t'-q'$. We shall denote $\vert q\vert=:\varepsilon$ and $\vert
q'\vert=:\varepsilon'$. Let $M_*:=M-q$ and $M_*':=M'-q'$ be the two
new hypersurfaces obtained by such geometric translations. In the new
coordinates, we naturally have two new defining equations
$w_*=\Theta_*(\bar z_*,t_*)$ and $\bar w_*'=\Theta_*'(\bar z_*',t_*')$
for $M_*$ and for $M_*'$ with $\Theta_*$ converging (at least) in
$\Delta_{2n-1}(0,\rho-\varepsilon)$ and with $\Theta_*'$ converging
(at least) in $\Delta_{2n-1}(0,\rho'-\varepsilon')$.  The explicit
expression of $\Theta_*'$ will be computed in a while.
Let $h_*(t_*):=h(q+t_*)$. Let 
$\bar\nu_*':=(\bar\lambda_*',\bar\mu_*'):=\bar\nu'-\bar q'$.
Define the transformed reflection 
function $\mathcal{R}_{*h_*}'(t_*,\bar\nu_*')$ accordingly.

\def\thelemma{3.37}\begin{lemma}
If $\varepsilon<\sigma$ and $\varepsilon'<\sigma'$, then the
reflection function $\mathcal{R}_{*h_*}'(t_*,\bar\nu_*'):=
\bar\mu_*'-\Theta_*'(\bar \lambda_*',h_*(t_*))$ extends
holomorphically to the polydisc $\Delta_n(0,\sigma-\varepsilon)\times
\Delta_n(0,\sigma'-\varepsilon')$.
\end{lemma}

\proof
At first, we compute the defining equation of $M_*'$.
To obtain the explicit expression of $\Theta_*'$, it suffices to transform
the equation
\def\theequation{3.38}\begin{equation}
\bar w'-\bar w_{q'}=\Theta'(\bar z',t')-\Theta'(\bar z_{q'}',t_{q'}')=
\sum_{\beta\in\N^{n-1}}\,
(\bar z')^\beta\,\Theta_\beta'(t')-\sum_{\beta\in\N^{n-1}}\,
(\bar z_{q'})^\beta\,\Theta_{\beta}'(t_{q'}')
\end{equation}
in the form 
\def\theequation{3.39}\begin{equation}
\bar w_*'=\Theta_*'(\bar z_*',t_*')=
\sum_{\beta\in\N^{n-1}}\,
(\bar z_*')^\beta\,
\Theta_{*\beta}'(t_*').
\end{equation}
Differentiating with respect to $\bar z'$ and 
setting $\bar z':=\bar z_{q'}'$, we obtain
\def\theequation{3.40}\begin{equation}
\left\{
\aligned
{}&
\Theta_{*0}'(t_*'):=
\sum_{\gamma\in\N^{n-1}}\,
(\bar z_{q'}')^\gamma\, \Theta_\gamma'(q'+t_*')-
\sum_{\gamma\in\N^{n-1}}\,(\bar z_{q'}')^\gamma\,
\Theta_\gamma'(q'),\\
&
\Theta_{*\beta}'(t_*'):=
\Theta_\beta'(q'+t_*')+
\sum_{\gamma\in\N_*^{n-1}}\,
(\bar z_{q'}')^\gamma\,
\Theta_{\beta+\gamma}'(q'+t_*')\,
{(\beta+\gamma)!\over \beta! \ \gamma!},
\endaligned\right.
\end{equation}
for all $\beta\in\N_*^{n-1}$. Now, suppose that the reflection 
function $\mathcal{R}_h'(t,\bar \nu')$ in the old system of coordinates 
extends holomorphically to the product polydisc
$\Delta_n(0,\sigma)\times \Delta_n(0,\sigma')$ as
a function that we shall denote by
\def\theequation{3.41}\begin{equation}
R'(t,\bar\nu'):=\bar\mu'-\sum_{\beta\in\N^{n-1}}\,
(\bar\lambda')^\beta\, \theta_\beta'(t).
\end{equation}
By Lemma~3.16, the functions $\theta_\beta'(t)$ are holomorphic in
$\Delta_n(0,\sigma)$ and they extend holomorphically the
$\mathcal{C}^\infty$-smooth CR functions $\Theta_\beta'(h(t))$ defined
on $M\cap \Delta_n(0,\rho)$. Immediately, $R'$ is holomorphic in 
an obvious product polydisc centered at $q\times \bar q'$, namely 
in $\Delta_n(q,\sigma-\varepsilon)\times\Delta_n(\bar q',
\sigma'-\varepsilon')$.
Let $t_*:=t-q$ and $\bar\nu_*':=\bar\nu'-\bar q'$. The unique 
function $R_*'(t_*,\bar\nu_*')$ satisfying 
\def\theequation{3.42}\begin{equation}
R'(t,\bar\nu')=R_*'(t_*,\bar\nu_*')=
\bar\mu_*'-\sum_{\beta\in\N^{n-1}}\,
(\bar\lambda_*')^\beta\,\theta_{*\beta}'(t_*)
\end{equation}
possesses coefficients necessarily given by
\def\theequation{3.43}\begin{equation}
\left\{
\aligned
{} &
\theta_{*0}'(t_*):=\sum_{\gamma\in\N^{n-1}}\,
(\bar z_{q'}')^\gamma\,\theta_\gamma'(q+t_*)-\sum_{\gamma\in\N^{n-1}}\,
(\bar z_{q'}')^\gamma\,\theta_\gamma'(q),\\
& 
\theta_{*\beta}'(t_*):=\theta_\beta'(q+t_*)+
\sum_{\gamma\in\N_*^{n-1}}\,(\bar z_{q'}')^\gamma\,
\theta_{\beta+\gamma}'(q+t_*)\,
{(\beta+\gamma)!\over \beta! \ \gamma!},
\endaligned\right.
\end{equation}
for all $\beta\in\N_*^{n-1}$. 
In the new coordinate system, the reflection function centered at $q\times
\overline{h(q)}$ can be defined as
\def\theequation{3.44}\begin{equation}
\mathcal{R}_{*h_*}'(t_*,\bar\nu_*'):=
\bar\mu_*'-\sum_{\beta\in\N^{n-1}}\,
(\bar\lambda_*')^\beta\,
\Theta_{*\beta}'(h_*(t_*)),
\end{equation}
for $t_*\in M_*$ with $\vert t_*\vert<\sigma-\varepsilon$.
Substituting $t_*'$ by $h_*(t_*)$ in the equations~(3.40) and
using afterwards that the $\theta_\beta'(t)$ extend the
$\Theta_\beta'(h(t))$, we deduce that the functions
\def\theequation{3.45}\begin{equation}
\left\{
\aligned
{}&
\Theta_{*0}'(h_*(t_*))=
\sum_{\gamma\in\N^{n-1}}\,
(\bar z_{q'}')^\gamma\, \Theta_\gamma'(h(q+t_*))-
\sum_{\gamma\in\N^{n-1}}\,(\bar z_{q'}')^\gamma\,
\Theta_\gamma'(q'),\\
&
\Theta_{*\beta}'(h_*(t_*))=
\Theta_\beta'(h(q+t_*))+
\sum_{\gamma\in\N_*^{n-1}}\,
(\bar z_{q'}')^\gamma\,
{(\beta+\gamma)!\over \beta! \ \gamma!}\,
\Theta_{\beta+\gamma}'(h(q+t_*))
\endaligned\right.
\end{equation}
extend holomorphically to the polydisc
$\Delta_n(0,\sigma-\varepsilon)$ as functions of $t_*$ given by the
right hand sides of~(3.43). The convergence of these series follows from
the Cauchy estimates on the $\theta_\beta'(t)$.  This completes the
proof of Lemma~3.37.
\endproof

\subsection*{3.46.~Delocalization and propagation} At this stage, 
we can summarize what the term ``reflection function'' really
means. Let $h: M\to M'$ be a (not necessarily local)
$\mathcal{C}^\infty$-smooth CR mapping between two connected real
analytic CR manifolds.  For any product of points $p\times
\overline{h(p)}$ lying in the graph of $\bar h$ in $M\times
\overline{M'}$ and for any system of coordinates $t'$ vanishing at
$p':=h(p)$ in which the complex defining equation of $M'$ is an
uniquely defined {\it graph}\, of the form $\bar w'=\Theta'(\bar
z',t')$, we define the associated reflection {\it centered}\, at
$p\times \overline{p'}$ by $\mathcal{R}_h(t,\bar\nu'):=\bar\mu'-
\Theta'(\bar\lambda',h(t))$. If it exists, its holomorphic extension
at $p\times \overline{p'}$ is unique, thanks to the uniqueness
principle on the boundary ([P1]). Also, its holomorphic extension does
not depend on the system of coordinates $t'$ vanishing at $p'$. And
finally, its holomorphic extension propagates at nearby points.
Although for some real analytic hypersurface $M'$ there does not exist
a {\it global}\, defining equation of the form $\bar w'=\Theta'(\bar
z',t')$, we believe that the transformation rules explained in
Lemmas~3.22 and~3.37 justify that we speak of ``the'' reflection
function.

The two analytic relations~(3.35) and~(3.40) are extremely important.
In \S3.47 just below, we shall see that they permit to establish that
certain CR geometric concepts defined in terms of the collection
$(\Theta_\beta'(t'))_{\beta\in\N^{n-1}}$ are biholomorphically
invariant.

\subsection*{3.47.~The exceptional locus of $M'$}
As above, let $p'\in M'$ and assume that the defining equation of $M'$
converges normally in the polydisc $\Delta_{2n-1}(0,\rho')$. Let us
consider the infinite Jacobian matrix of the infinite 
holomorphic mapping $\mathcal{Q}_\infty'(t')=
(\Theta_\beta'(t'))_{\beta\in\N^{n-1}}$ introduced in \S1.10:
\def\theequation{3.48}\begin{equation}
\mathcal{J}_\infty(t'):=
(\partial\Theta_\beta'(t')/\partial t_j')_{\beta\in\N^{n-1}, 
1\leq j\leq n}.
\end{equation}
Concretely, by ordering the multi-indices $\beta$, we may think of
$\mathcal{J}_\infty(t')$ as a horizontally infinite $\infty\times n$
complex matrix.  Also, it is convenient to truncate this matrix by limiting
the multi-indices to run over $\vert \beta\vert \leq k$. Let us denote
such finite matrices by
\def\theequation{3.49}\begin{equation}
\mathcal{J}_k(t'):=
(\partial\Theta_\beta'(t')/\partial t_j')_{\vert\beta\vert\leq k, 
1\leq j\leq n}.
\end{equation}
As a holomorphic mapping of $t'$, the generic rank of
$\mathcal{J}_k(t')$ increases with $k$. Let $\chi_{M'}'$ denote the
maximal generic rank of these finite matrices. Equivalently, there
exists a minor of size $\chi_{M'}'$ of the matrix $\mathcal{J}_\infty$
which does not vanish identically as a holomorphic function of $t'$,
but all minors of size $(\chi_{M'}'+1)$ of $\mathcal{J}_\infty(t')$ do
vanish identically. We call this integer the {\it generic rank}\, of
the infinite matrix $\mathcal{J}_\infty(t')$. Of course, $\chi_{M'}'$
is at least equal to $1$, because the term $\Theta_0'(t')$ does not
vanish identically and is nonconstant ({\it see}~\S3.1). So we have
$1\leq \chi_{M'}'\leq n$. Apparently, the integer $\chi_{M'}'$ seems
to depend on $p'$ and on the choice of coordinates centered at $p'$,
but in fact it is a biholomorphic invariant of the hypersurface $M'$
itself, which explains in advance the notation. Recall that $M'$ is
connected, which is important. We shall check this invariance in two
steps.

\def\thelemma{3.50}\begin{lemma}
Let $p'\in M'$, let $t'$ be a system of coordinates vanishing at $p'$
and let $t''$ be another system of coordinates vanishing at $p'$
defined by $t''=\Lambda(t')$ as in Lemma~3.22. Then the two generic
ranks of the associated infinite Jacobian matrices are identical.
\end{lemma}

\proof
Looking at the family of relations~(3.35) and applying the rank
inequality for composed holomorphic mappings, we see that the generic
rank of $\mathcal{J}_\infty(t'')$ is certainly less than or equal to
the generic rank of $\mathcal{J}_\infty(t')$.  As the mapping
$\Lambda$ is invertible, a relation similar to~(3.35) holds if we
reverse the r\^oles of $t'$ and $t''$, and we get the opposite
inequality between generic ranks.
\endproof

\def\thelemma{3.51}\begin{lemma}
Let $p'\in M'$, let $q'\in M'$ be close to $p'$ as Lemma~3.37 and
consider the infinite Jacobian matrix $\mathcal{J}_\infty(t_*')$ 
associated with the functions $\Theta_{*\beta}'(t_*')$ defined by~(3.40).
Then the generic ranks of $\mathcal{J}_\infty(t')$ and
of $\mathcal{J}_\infty(t_*')$ coincide.
\end{lemma}

\proof
This is immediate, because the relation~(3.40) between the
two collections $(\Theta_{*\beta}'(t_*'))_{\beta\in\N^{n-1}}$
and $(\Theta_\beta'(t'))_{\beta\in\N^{n-1}}$ is linear, 
upper triangular and invertible.
\endproof

So we may prove that $\chi_{M'}'$ is a global biholomorphic
invariant of the connected hypersurface $M'$.  
Indeed, any two points $p_1'\in M'$ and $p_2'\in
M'$ can be connected by a finite chain of intermediate points which
are contained in pairs of overlapping coordinate system for which
Lemmas~3.50 and~3.51 apply directly.

Here is an interesting and useful application. Locally in a
neighborhood of an arbitrary point $p'\in M'$, we may define a {\it
proper complex analytic}\, subset of $\Delta_n(0,\rho')$ denoted by
$\mathcal{E}'$ which is obtained as the vanishing locus of all the
minors of size $\chi_{M'}'$ of $\mathcal{J}_\infty(t')$. As in the
proofs of Lemmas~3.50 and~3.51, by looking more closely at the two
families of infinite relations~(3.35) and~(3.40), we observe that the
set of points $t'$ close to $p'$ at which the {\it rank}\, of
$\mathcal{J}_\infty(t')$ is maximal equal to $\chi_{M'}'$ is
independent of coordinates.  Consequently, the complex analytic set
$\mathcal{E}'$, which we shall call the {\it extrinsic exceptional
locus of $M'$}, is an invariant complex analytic subset defined in a
neighborhood of $M'$ in $\C^n$. Moreover, $\mathcal{E}'$ is {\it
proper} ({\it i.e.}  of dimension $\leq n-1$), because
$\chi_{M'}'\geq 1$, so there is at least one not identically zero
minor in the definition of $\mathcal{E}'$. The {\it intrinsic
exceptional locus of $M'$} denoted by $E_{M'}'$ is defined to be the
intersection of $\mathcal{E}'$ with $M'$. This is also a {\it
proper}\, real analytic subset of $M'$ (maybe empty).

\def\thelemma{3.52}\begin{lemma}
If $M'$ is globally minimal, then the real dimension of
$E_{M'}'$ is less than or equal to $(2n-3)$.
\end{lemma}

\proof
Suppose on the contrary that there exists a stratum $S$ of real 
dimension $(2n-2)$. This stratum cannot be generic at any point, 
because otherwise $\mathcal{E}'$ which contains $S$ would be
of complex dimension $n$. So $S$ is a complex hypersurface
contained in $M'$, contradicting local minimality at every point.
\endproof

This dimension estimate should be compared to that of the Levi
degeneracy locus: unless $M'$ is everywhere Levi degenerate, the set
of points at which $M'$ is Levi degenerate is a proper real analytic
subvariety, but in general of dimension less than or equal to
$(2n-2)$, with this bound attained. This is so because the Levi
degeneracy locus is {\it not}\, contained in a complex analytic subset
of a neighborhood of $M'$. The fact that the real codimension of
$E_{M'}'$ is at least two will be crucial for the proof of
Theorems~9.2 and~9.3 below.
 
\subsection*{3.53.~Local product structure at a Zariski-generic point}
In the beginning of \S4 below we shall need the following
geometric straightening statement.

\def\thelemma{3.54}\begin{lemma}
In a small neighborhood of an arbitrary point $q'\in M'\backslash
E_{M'}'$, the hypersurface $M'$ is biholomorphic to a product
$\underline{M}_{q'}'\times\Delta^{n-\chi_{M'}'}$ by a polydisc of
dimension $(n-\chi_{M'}')$, where $\underline{M}_{q'}'$ is a real
analytic hypersurface in $\C^{\chi_{M'}'}$. Furthermore, at the
point $\underline{q'}$, the rank of an associated infinite matrix
$\underline{\mathcal{J}}_\infty(\underline{t'})$, where
$\underline{t}'\in\C^{\chi_{M'}'}$ are holomorphic coordinates
vanishing at $\underline{q}'$, is maximal equal to $\chi_{M'}'$.
\end{lemma}

\proof
Choose coordinates $t'$ vanishing at $q'$. By assumption, the mapping
$t'\mapsto (\Theta_\beta'(t'))_{\vert \beta\vert \leq k}$ is of
constant rank $\chi_{M'}'$ for all $t'$ near the origin and for all
$k$ large enough. By the rank theorem, it follows that the union of
level sets $\mathcal{F}_{r'}:=\{t':
\Theta_\beta'(t')=\Theta_\beta'(r'), \forall\,\beta\in\N^{n-1}\}$ for
$r'$ running in a neighborhood of $q'$ do constitute a local
holomorphic foliation by complex leaves of dimension
$n-\chi_{M'}'$. We can straighten this foliation in a neighborhood
of $q'$ so that (after an eventual dilatation) $\C^n$ decomposes as
the product $\Delta^{\chi_{M'}'}\times \Delta^{n-\chi_{M'}'}$,
where the second term corresponds to the leaves of this foliation.  In
these new straightening coordinates, which we will denote by $t''$, we
claim that the leaves of this foliation are again defined by the level
sets of the functions $\Theta_\beta''(t'')$, namely
$\mathcal{F}_{r''}:=\{t'': \Theta_\beta''(t'')=\Theta_\beta''(r''),
\forall\,\beta\in\N^{n-1}\}$. This is so, thanks to the important
relations~(3.35). For simplicity, let us denote these coordinates
again by $t'$ instead of $t''$. We claim that if the point $r'$
belongs to $M'$, then its leaf $\mathcal{F}_{r'}$ is entirely contained
in $M'$ in a neighborhood of $q'$. Indeed, let $s'\in \mathcal{F}_{r'}$, 
so we have $\Theta_\beta'(s')=\Theta_\beta'(r')$ for all $\beta\in\N^{n-1}$
by definition. It follows first that 
\def\theequation{3.55}\begin{equation}
0=\bar w_{r'}'-\Theta'(\bar z_{r'}',t_{r'}')=
\bar w_{r'}'-\Theta'(\bar z_{r'}',t_{s'}').
\end{equation}
Next, thanks to the reality of $M'$, there exists a
nonzero holomorphic function $a'(t',\tau')$, where 
$\tau'=(\zeta',\xi')\in\C^{n-1}\times \C$, such that
\def\theequation{3.56}\begin{equation}
\xi'-\Theta'(\zeta',t')\equiv a'(t',\tau')\,\left[
w'-\overline{\Theta}'(z',\tau')\right],
\end{equation}
for all $t',\tau'$ running in a neighborhood of the origin.
Using crucially this identity, we can transform~(3.55) as follows
\def\theequation{3.57}\begin{equation}
0=w_{s'}'-\overline{\Theta}'(z_{s'}',\bar t_{r'}').
\end{equation}
Now, conjugating this new identity, we get $\bar w_{s'}'-\Theta'(\bar
z_{s'}', t_{r'}')=0$ and finally, using a second time
$\Theta_\beta'(s')=\Theta_\beta'(r')$ for all $\beta\in\N^{n-1}$, we
obtain
\def\theequation{3.58}\begin{equation}
\bar w_{s'}'-\Theta'(\bar z_{s'}',t_{s'}')=0,
\end{equation}
which shows that $s'\in M'$, as claimed. In summary, in the straightened
coordinates $(\underline{t}',\widetilde{t}')\in
\C^{\chi_{M'}'}\times \C^{n-\chi_{M'}'}$, those leaves
$\{\widetilde{t}'=ct.\}$ intersecting $M'$ are entirely contained in $M'$.
It follows that there exists a defining equation for $M'$ in a
neighborhood of the origin which is of the form
\def\theequation{3.59}\begin{equation}
\overline{\underline{w}}'=\underline{\Theta}'(\overline{\underline{z}}',
\underline{t}'),
\end{equation}
namely it is independent of the coordinates $\widetilde{t}'$.  We
define $\underline{M}_{q'}'$ to be the hypersurface of
$\C^{\chi_{M'}'}$ defined by the equation~(3.59). The infinite
Jacobian matrix $\mathcal{J}_\infty(t')$ of $M'$ therefore coincides
with the infinite Jacobian matrix of $\underline{M}_{q'}'$. By
assumption, $\mathcal{J}_\infty(t')$ is of rank $\chi_{M'}'$ at the
origin (this means that all finite submatrices $\mathcal{J}_k(t')$ are
of rank $\chi_{M'}'$ for all large enough $k$). So the rank at the
origin of $\underline{\mathcal{J}}_\infty(\underline{t}')$ is also
equal to $\chi_{M'}'$. The proof of Lemma~3.54 is complete.
\endproof

\subsection*{3.60.~Pointwise nondegeneracy conditions on $M'$}
We shall call the (always connected) hypersurface $M'$ {\it
holomorphically nondegenerate}\, if $\chi_{M'}'=n$. By examinating the
proof of Lemma~3.54, one can see that this definition coincides with
the original definition of Stanton [St1,2] in terms of tangent
holomorphic vector fields ({\it cf.}~also [Me5,\S9]).  By \S3.47
above, holomorphic nondegeneracy is a global property of $M'$.
Furthermore, we shall say that $M'$ is {\it finitely nondegenerate}\,
at the point $p'$ if for one (hence for all) system(s) of coordinates
vanishing at $p'$, the rank of $\mathcal{J}_\infty(t')$ is equal to
$n$ at the origin. By the above definitions, a connected real analytic
hypersurface $M'$ is holomorphically nondegenerate if and only if
there exists a proper complex analytic subset of a
neighborhood of $M'$ in $\C^n$, namely the extrinsic exceptional locus
$\mathcal{E}'$, such that $M'$ is finitely nondegenerate at every
point of $M'$ not belonging to $\mathcal{E}'$. Also, Lemma~3.54 above
may be interpreted as a sort of geometric quotient procedure: locally
in a neighborhood of a Zariski-generic point $q'\in M'$, {\it
i.e.}~for $q'\not\in E_{M'}'$, after dropping the innocuous polydisc
$\Delta^{n-\chi_{M'}'}$, we are left with a finitely nondegenerate
real analytic hypersurface $\underline{M}_{q'}'$ in a smaller complex
affine space. Finally, we shall say that $M'$ is {\it essentially
finite}\, at the point $p'$ if for one (hence for all) system(s) of
coordinates vanishing at $p'$, the local holomorphic mappings $t'\mapsto
(\Theta_{\beta}'(t'))_{\vert\beta\vert\leq k}$ are finite-to-one
in a neighborhood of the origin for all $k$ large enough.  It can be
checked that this definition coincides with the one introduced in [DW]
and subsequently studied by many authors. We shall consider
essentially finite hypersurfaces in \S9 below.

\subsection*{3.61.~Conclusion} All the considerations of this paragraph
support well the thesis that {\it the collection of holomorphic
functions $(\Theta_\beta'(t'))_{\beta\in\N^{n-1}}$ is the most
important analytic object attached to a real analytic hypersurface
$M'$ localized at one of its points}.

\section*{\S4.~Extension across a Zariski dense open subset of $M$}

\subsection*{4.1.~Holomorphic extension at a Zariski-generic point} Let 
$h:M\to M'$ be a $\mathcal{C}^\infty$-smooth CR diffeomorphism
between two connected real analytic hypersurfaces in $\C^n$.

\def\thelemma{4.2}\begin{lemma}
If $M$ is globally minimal, then $M'$ is also globally minimal.
\end{lemma} 

\proof
Indeed, as $h$ is CR, it sends every $\mathcal{C}^\infty$-smooth curve
$\gamma$ of $M$ running into complex tangential directions
diffeomorphically onto a curve $\gamma':=h(\gamma)$ also running in
complex tangential directions.  Then Lemma~4.2 is a direct consequence
of the definition of CR orbits. We do not enter the details.
\endproof

The starting point of the proof of Theorem~1.9 is to show that the various
reflection functions already extend holomorphically to a neighborhood
of $q\times \overline{h(q)}$ for all points $q$ running in the Zariski
open subset $M\backslash E_M$ of $M$, where $E_M$ is the intrinsic
exceptional locus of $M$ defined in the end of \S3.47 above. It is
convenient to observe first that $h$ maps $E_M$ bijectively onto
$E_{M'}'$.

\def\thelemma{4.3}\begin{lemma}
A point $q\in M$ belongs to $M\backslash E_M$ if and only if
its image $h(q)$ belongs to $M'\backslash E_{M'}'$. Furthermore,
$\chi_M=\chi_{M'}'$.
\end{lemma}

\proof
Let $q\in M$ be arbitrary, let $t$ be coordinates
vanishing at $q$ and let $t'$ be coordinates vanishing at $q':=h(q)$
in which we have
\def\theequation{4.4}\begin{equation}
\overline{g(t)}-\Theta'(\overline{f(t)},
h(t))=a(t,\bar t)\,\left[\bar w-\Theta(\bar z,t)\right],
\end{equation}
for all $t\in M$ close to the origin and for some nonvanishing
function $a(t,\bar t)$ of class $\mathcal{C}^\infty$.  By developping
the Taylor series of all $\mathcal{C}^\infty$-smooth functions
in~(4.4) and by polarizing, we see that the Taylor series $H$ of $h$
at the origin induces a formal mapping between $(M,q)$ and $(M',q')$,
namely there exists a formal power series $A(t,\tau)$ with nonzero
constant term such that the following identity holds between formal
power series in the $2n$ variables $(t,\tau)$:
\def\theequation{4.5}\begin{equation}
\overline{G}(\tau)-\Theta'(\overline{F}(\tau),H(t))\equiv
A(t,\tau)\,\left[ \xi-\Theta(\zeta,t)\right].
\end{equation}
Now, the computations of Lemma~3.22 can be performed at a purely
formal level, replacing the mapping $\Lambda$ there by the formal
mapping $H$. We obtain a relation similar to~(3.35), interpreted at
the formal level, with
$\Lambda$ replaced by $H$. Using the invertibility of $H$ to get a
second relation like~(3.35) with $\Lambda$ replaced by $H^{-1}$, it
then follows that the rank of the mapping $t\mapsto
(\Theta_\beta(t))_{\beta\in\N^{n-1}}$ at $q$ is the same as the rank
of the mapping $t'\mapsto (\Theta_\beta'(t'))_{\beta\in\N^{n-1}}$ at
$h(q)$. This property yields the desired conclusion.
\endproof

Thus, the starting point of the proof of Theorem~1.9 is the following
Zariski dense holomorphic extension result.

\def\thelemma{4.6}\begin{lemma}
If $h:M\to M'$ is a $\mathcal{C}^\infty$-smooth CR diffeomorphism
between two globally minimal real analytic hypersurfaces in $\C^n$,
then for every point $q\in M\backslash E_M$ lying outside the
intrinsic exceptional locus of $M$ and for every choice of a
coordinate system vanishing at $q':=h(q)$ in which $(M',q')$ is
represented by $\bar w'=\Theta'(\bar z',t')$, the associated reflection
function $\mathcal{R}_h'(t,\bar\nu')=\bar\mu'-\Theta'(
\bar\lambda',h(t))$ extends holomorphically to a neighborhood of
$q\times \overline{q'}$ in $\C^n\times\C^n$.
\end{lemma}

\proof
First, by Lemma~4.3, we already know that $q'$ does not belong to
$E_{M'}'$ and that $\chi_{M}= \chi_{M'}'$. For short, let us denote
this integer by $\chi$.  By Lemma~3.22, the holomorphic extendability
of the reflection function is invariant, so let us choose adapted
convenient coordinates. Using Lemma~3.54, we can find coordinates
near $q'\in M'$ of the form $t'=(z',v',w')\in\C^{\chi-1}\times
\C^{n-\chi}\times\C^1$ in which the equation of $M'$ near the origin
is given by $\bar w'=\Theta'(\bar z',z',w')$.
Notice that the $(v',\bar v')$ coordinates do not
appear in the defining equation, because of the product structure. We
do the same straightening near $q\in M$, so that we can split the
coordinates as $t=(z,v,w)\in\C^{\chi-1}\times \C^{n-\chi}\times\C^1$
in which the equation of $M$ near the origin is also given in the form
$\bar w=\Theta(\bar z,z,w)$. Finally, we split the mapping accordingly
as $h=(f,l,g)\in\C^{\chi-1}\times\C^{n-\chi}\times \C^1$. It is
important to notice that in these coordinates, the reflection function
$\mathcal{R}_h'(t,\bar\lambda',\bar\upsilon',
\bar\nu')=\bar\mu'-\Theta'(\bar\lambda',f(t),g(t))$, where
$(\bar\lambda',\bar\upsilon',\bar\nu')\in\C^{\chi-1}\times
\C^{n-\chi}\times\C^1$, {\it neither}\, depends on the $(n-\chi)$
middle components $(l_1,\dots, l_{n-\chi})=(h_\chi,\dots,h_{n-1})$
of $h$ {\it nor}\, on $\bar\upsilon'$. Clearly, to show that this reflection
function extends holomorphically at $q\times \overline{q'}$, it would
suffice to show that the $\chi$ components $(f_1,\dots,f_{\chi-1},
g)=(h_1,\dots,h_{\chi-1},h_n)$ of $h$ extend holomorphically to 
a neighborhood of the
origin. We need some notation. Let $\underline{h}$ denote these $\chi$
special components $(f,g)$, let $\underline{M}$ denote the hypersurface $\bar
w=\Theta(\bar z,z,w)$ of $\C^\chi$ and similarly let $\underline{M}'$
denote the hypersurface $\bar w'=\Theta'(\bar z',z',w')$ of
$\C^\chi$. {\it A priori}, it is not clear whether $\underline{h}$ induces a
$\mathcal{C}^\infty$-smooth CR mapping between
$(\underline{M},\underline{q})$ and $(\underline{M}',\underline{q}')$,
since $\underline{h}$ might well depend on the variables 
$(v_1,\dots,v_{n-\chi})$.

\def\thelemma{4.7}\begin{lemma}
The $\chi$ components $(f_1,\dots,f_{\chi-1},g)$ of $h$ are independent
of the $(n-\chi)$ coordinates $v$. Consequently, the mapping $h$ induces a 
well defined CR mapping $\underline{h}:(\underline{M},q)\to 
(\underline{M}',q')$ of class $\mathcal{C}^\infty$.
\end{lemma}

\proof
Let $\overline{L}_1,\dots, \overline{L}_{n-1}$ be a commuting basis
of $T^{0,1}M$ with real analytic coefficients, for instance
$\overline{L}_j={\partial \over \partial \bar z_j}+\Theta_{\bar
z_j}(\bar z,z,w){\partial \over \partial \bar w}$ for
$j=1,\dots,\chi-1$ and also $\overline{L}_i={\partial \over\partial
\bar v_i}$, for $i=1,\dots,n-\chi$. Notice that
the $(1,0)$ vector field $L_i$, $i=1,\dots,n-\chi$ commute with the
$(0,1)$ vector fields $\overline{L}_j$, $j=1,\dots,\chi-1$. Since $h$ is a
$\mathcal{C}^\infty$-smooth CR diffeomorphism, after a possible
linear change of coordinates, we can assume that the determinant
$\mathcal{D}(t,\bar t):= \hbox{det} \, \left(\overline{L}_j \bar
f_k(\bar t)\right)_{1\leq j,k\leq \chi-1}$ is nonzero at the origin.
Applying the derivations $\overline{L}_1,\dots,
\overline{L}_{\chi-1}$ to the fundamental identity $\overline{g(t)}=
\Theta'(\overline{f(t)},\underline{h}(t))$ for $t$ on $(M,q)$, we get
first
\def\theequation{4.8}\begin{equation}
\overline{L}_j\overline{g(t)}=
\sum_{k=1}^{\chi-1}\,\overline{L}_j\overline{f_k(t)}\
{\partial\Theta'\over\partial \bar z_k'}
(\overline{f(t)},\underline{h}(t)).
\end{equation}
Shrinking $\sigma>0$ if necessary, we can assume that the
determinant $\mathcal{D}(t,\bar t)$ does not vanish for
all $\vert t\vert<\sigma$. By Cramer's rule, we can solve 
in~(4.8) the first order partial derivatives $\partial_{\bar z_k'}
\Theta'$ with respect
to the other terms. As in the proof of Lemma~3.22, by induction, it follows
that for every multi-index $\beta\in\N^{\chi-1}$, there exists a certain 
universal polynomial $R_\beta$ such that the following
relation holds for all $t\in M$ with $\vert t\vert < \sigma$:
\def\theequation{4.9}\begin{equation}
{1\over \beta!}\
{\partial^{\vert\beta\vert}\Theta'\over
\partial (\bar z')^\beta}
(\overline{f(t)},\underline{h}(t))=
{R_\beta(\{(\overline{L})^\gamma\,
\overline{\underline{h}(t)}\}_{\vert\gamma\vert\leq\vert\beta\vert})\over
\left[\mathcal{D}(t,\bar t)\right]^{2\vert\beta\vert-1}}.
\end{equation}
Next, since by assumption the point $q'$ does not belong to $E_{M'}'$,
the second sentence of Lemma~3.54 tells us that there exists a positive
integer $k$ such that the rank of the mapping $\C^\chi\ni
(z',w')\mapsto (\Theta_\beta'(z',w'))_{\vert\beta\vert\leq k}$ is maximal
equal to $\chi=\chi_{M'}'$. Writing the equalities~(4.9) only for
$\vert\beta\vert\leq k$ and applying the implicit function theorem, it
follows finally that we can solve $\underline{h}(t)$ with respect to
the derivatives of $\overline{\underline{h}(t)}$, namely there exist
$\chi$ {\it holomorphic}\, functions $\Omega_j$ 
in their variables such that for
$j=1,\dots,\chi$ and $t\in M$ with $\vert t\vert < \sigma$ (shrinking
$\sigma$ if necessary), we have:
\def\theequation{4.10}\begin{equation}
\underline{h}_j(t)=\Omega_j(\{(\overline{L})^\gamma\,
\overline{\underline{h}(t)}\}_{\vert\gamma\vert\leq k})=
\Omega_j(t,\bar t,\{\partial_{\bar t}^\gamma
\overline{\underline{h}(t)}\}_{\vert\gamma\vert\leq k}).
\end{equation}
Applying now the $n-\chi$ vector fields $L_i={\partial\over \partial
v_i}$, $i=1,\dots,n-\chi$, to these identities, using the fact that
these $\partial_{v_i}$ do commute with the antiholomorphic derivations
$\overline{L}_1^{\gamma_1}\cdots\overline{L}_{\chi-1}^{\gamma_{\chi-1}}$,
and noticing that the $\overline{\underline{h}(t)}$ are anti-CR, we
obtain that the $\partial_{v_i}\underline{h}_j(t)$ do vanish
identically on $M$ near the origin. Since the $\underline{h}_j(t)$ are
CR and of class $\mathcal{C}^\infty$, we already know that the
derivatives $\partial_{\bar v_i} \underline{h}_j(t)$ also vanish
identically on $M$ near the origin.  This proves that the
$\underline{h}_j$ are independent of the coordinates $(v,\bar v)$, as
desired.
\endproof

Finally, the following lemma achieves to prove that 
the reflection function, which only depends on $\underline{h}$, does
extend holomorphically to a neighborhood of $q\times \overline{q'}$,
as claimed.

\def\thelemma{4.11}\begin{lemma}
The mapping $\underline{h}$ extends holomorphically to a 
neighborhood of the origin.
\end{lemma}

\proof
The proof of this lemma is an easy generalization of the Lewy-Pinchuk
reflection principle and in fact, it can be argued that Lemma~4.11 is
almost completely contained in [P3] (and also in [W2,3], [DW]). Formally
indeed, the calculations in the proof of Lemma~4.7 above are totally
similar to the ones in the Levi nondegenerate case except for the
order of derivations. Of course, the interest of derivating further
the equations~(4.8) does not lie in this (rather evident or
gratuitous) generalization of the reflection principle. Instead, the
interest lies in the fact that there are large classes of everywhere
Levi-degenerate hypersurfaces for which it is natural to introduce the
concept of finite nondegeneracy expressed in terms of the fundamental
functions $\Theta_\beta'(t')$. Indeed, finite nondegeneracy correspond
to the (not rigorous, in the folklore) intuitive notion of ``higher order
Levi-forms''. Furthermore, holomorphically nondegenerate hypersurfaces
are almost everywhere finitely nondegenerate. In sum, from the point
of view of local analytic CR geometry, higher order derivations are
very natural.

Although Lemma~4.11 is explicitely stated or covered by [DW], [Ha],
[BJT], {\it etc.}, we shall summarize its proof for completeness.
Recall that by \S3.6, the components of $h$ extend holomorphically to
a global one-sided neighborhood $D$ of $M$ which contains one side $D_q$ of
$M$ at $q$.  Let $M^-$ denote the side of $(M,q)$ containing $D_q$ and
let $M^+$ denote the other side. As in the Lewy-Pinchuk reflection
principle, using the real analyticity of the coefficients of the
$\overline{L}_j$ and using the one-dimensional Schwarz reflection
principle in the complex lines $\{w=ct.\}$ which are transverse to $M$
near the origin, we observe that the functions $\Omega_j$ in the right
hand side of~(4.10) extend $\mathcal{C}^\infty$-smoothly to $M^+$ as
functions $\omega_j$ which are partially holomorphic with respect to
the transverse variable $w$. Since by ~(4.10) the values of the
$\underline{h}_j$ coincide on $(M,q)$ with the values of the
$\omega_j$ and since the $\underline{h}_j$ are already holomorphic
inside $D_q$, it follows from a rather easy (because everything is
$\mathcal{C}^\infty$-smooth) separate analyticity principle that the
$\underline{h}_j$ and the $\omega_j$ stick together in holomorphic
functions defined in a neighborhood of $q$. This provides the desired
holomorphic extension of $\underline{h}$ and hence the holomorphic
extension of $\mathcal{R}_h'$. The proofs of Lemmas~4.6 and~4.11 are
complete.  
\endproof

\subsection*{4.12.~Holomorphic extension of the mapping}
We end up this paragraph by showing that Theorem~1.9 implies
Theorem~1.14 (or equivalently Theorem~1.2). Under the assumptions of
Theorem~1.14, let $p\in M$ be arbitrary and let $t$ be coordinates
vanishing at $p$. By the holomorphic extendability of the reflection
function, we know that all the $\mathcal{C}^\infty$-smooth CR
functions $\Theta_\beta'(h(t))$ extend as holomorphic functions
$\theta_\beta'(t)$ to a fixed neighborhood of $p$. Thanks to the
holomorphic nondegeneracy of $M'$, there exist $n$ different
multi-indices $\beta^1,\dots,\beta^n\in\N^{n-1}$ such that the
generic rank of the holomorphic 
mapping $t'\mapsto (\Theta_{\beta^k}'(t'))_{1\leq
k\leq n}$ equals $\chi_{M'}'=n$, or equivalently
\def\theequation{4.13}\begin{equation}
{\rm det}\, ([\partial\Theta_{\beta^k}'/
\partial t_l'](t'))_{1\leq k,l\leq n}\not\equiv 0 \ \ \ 
{\rm in} \ \C\{t'\}.
\end{equation}
Let $H(t)$ denote the formal Taylor series of $h$ at the origin. Since
the Jacobian determinant of $h$ at $0$ does not vanish, it follows
that~(4.13) holds in $\C\dl t\dr$ after $t'$ is replaced by $H(t)$. Then
the holomorphic extendability of $h$ at the origin is covered by the
following assertion.

\def\thelemma{4.14}\begin{lemma}
Let $p\in M$, let $t$ be coordinates vanishing at $p$,
let $h_1,\dots,h_n$ be CR functions of class $\mathcal{C}^\infty$ on
$(M,p)$ vanishing at the origin, let $H_j(t)$ denote the formal Taylor
series of $h_j$ at $0$, let $Q_1'(t'),\dots,Q_n'(t')$ be holomorphic
functions satisfying 
\def\theequation{4.15}\begin{equation}
{\rm det} \, ([\partial Q_k'/\partial
t_l'](H(t)))_{1\leq k,l\leq n}\not\equiv 0 \ \ \ {\rm in} \
\C\dl t\dr. 
\end{equation}
Assume that there exist holomorphic functions $q_1'(t),\dots,q_n'(t)$
defined in a neighborhood of the origin such that $Q_k'(h(t))\equiv
q_k'(t)$ for all $t\in (M,p)$ close to the origin. Then
$h_1(t),\dots,h_n(t)$ extend holomorphically to a neighborhood of the
origin.
\end{lemma}

\proof
Clearly, the holomorphic functions $S_{k,l}'(t',t_*')$ defined by the
relations
\def\theequation{4.16}\begin{equation}
Q_k'(t')-Q_k'(t_*')=\sum_{l=1}^n\,
S_{k,l}'(t',t_*')\,
(t_l'-t_{*l}')
\end{equation}
satisfy the relation $S_{k,l}'(t',t')= [\partial Q_k'(t')/\partial
t_l'](t')$. We first prove that the Taylor 
series $H_j(t)$ are convergent.  By the
Artin approximation theorem ([Ar]), for every integer $N\in \N_*$,
there exists a {\it converging}\, power series mapping
$\mathcal{H}(t)\in\C\{t\}^n$ with  $\mathcal{H}(t)\equiv H(t) \ {\rm
mod}\, \vert t\vert^N$ such that $Q_k'(\mathcal{H}(t))\equiv
q_k'(t)$. If $N$ is large enough, it follows from the main assumption
of Lemma~4.14 that the following formal determinant does not vanish
identically in $\C\dl t\dr$:
\def\theequation{4.17}\begin{equation}
{\rm det}\,(S_{k,l}'(H(t),\mathcal{H}(t))
)_{1\leq k,l\leq n}\not\equiv 0.
\end{equation}
Finally, by the relation 
\def\theequation{4.18}\begin{equation}
\left\{
\aligned
{}&
0=q_k'(t)-q_k'(t)=Q_k'(H(t))-Q_k'(\mathcal{H}(t))=\\
&  \ \ 
=\sum_{l=1}^n\,S_{k,l}'(H(t),\mathcal{H}(t))\,
[H_l(t)-\mathcal{H}_l(t)]
\endaligned\right.
\end{equation}
and thanks to the invertibility of the matrix $S_{k,l}'$ ({\it see}~(4.17)),
we deduce that $H_j(t)=\mathcal{H}_j(t)$ is convergent, as claimed.
Secondly, for $t\in (M,p)$ close to the origin, we again use~(4.16)
with $t':=h(t)$ and $t_*':=\mathcal{H}(t)$, which yields a relation 
like~(4.18) with $H(t)$ replaced by $h(t)$, namely
\def\theequation{4.19}\begin{equation}
0=Q_k'(h(t))-Q_k'(\mathcal{H}(t))=
\sum_{l=1}^n\,S_{k,l}'(h(t),\mathcal{H}(t))\,
[h_l(t)-\mathcal{H}_l(t)].
\end{equation}
Then the corresponding determinant~(4.17) (with $H(t)$ replaced by
$h(t)$) does not vanish identically on $(M,p)$, because it has a
nonvanishing formal Taylor series by~(4.17) and because $(M,p)$ is
generic. Consequently, relation~(4.19) implies that $h_l(t)\equiv
\mathcal{H}_l(t)$ for all $t\in (M,p)$ close to the origin. This completes
the proof of Lemma~4.14 (similar arguments are provided in [N]).
Also, the proof of Theorem~1.14 (taking Theorem~1.9 for granted)
is complete.
\endproof

\section*{\S5. Situation at a typical point of non-analyticity}

Thus, we already know that $\mathcal{R}_h'$ is analytic at every point
$q\times \overline{h(q)}$ for $q$ running in the open dense subset
$M\backslash E_M$ of $M$. It remains to show that $\mathcal{R}_h'$ is
analytic at all the points $p\times \overline{h(p)}$, where $p\in
E_M$, which entails $h(p)\in E_{M'}'$ by Lemma~4.3.  This objective
constitutes the principal task of the demonstration. In fact, we shall
prove a slightly more general semi-global statement which we summarize
as follows.

\def\thetheorem{5.1}\begin{theorem}
Let $h: M\to M'$ be a $\mathcal{C}^\infty$-smooth CR diffeomorphism
between two globally minimal real analytic hypersurfaces in $\C^n$. If
the local reflection mapping $\mathcal{R}_h'$ is analytic at {\rm
one}\, point $q\times \overline{h(q)}$ of $M\times M'$, then it is
analytic at {\rm every} point $p\times \overline{h(p)}$ of the graph
of $\bar h$ in $M\times M'$.
\end{theorem}

To prove Theorem~5.1, we shall proceed by contradiction. We define
the following subset of $M'$
\def\theequation{5.2}\begin{equation}
\mathcal{A}':=\{p'\in M':\mathcal{R}_h' \ \text{\rm is analytic in a 
neighborhood of} \ h^{(-1)}(p')\times \overline{p'}\}.
\end{equation}
A similar subset $\mathcal{A}$ of $M$ such that $h$ maps $\mathcal{A}$
bijectively onto $\mathcal{A}'$ can be defined, but in fact, it will
be more adapted to our purposes to work in $M'$ with $\mathcal{A}'$.
Recall that we already know that Theorem~1.9 implies
Theorem~1.2. However, for a direct proof of Theorem~1.2 ({\it
cf.}~\S2), it would have been convenient to define the set
$\mathcal{A}'$ above as the set of point $p'\in M'$ such that $h$ is
analytic in a neighborhood of $h^{-1}(p')$. Anyway, the set
$\mathcal{A}'$ defined by~(5.2) is nonempty, by the assumption of
Theorem~5.1. For the proof of Theorem~1.9, $\mathcal{A}'$ is also
nonempty, because it contains $M'\backslash E_{M'}'$ thanks to
Lemma~4.6 above. So let us start with~(5.2). If $\mathcal{A}'=M'$,
Theorem~5.1 would be proved, gratuitously.  As in \S2.2, we shall
therefore suppose that its complement $E_{\rm
na}':=M'\backslash\mathcal{A}'$ is nonempty and we shall endeavour to
derive a contradiction. In fact, to derive a contradiction, it clearly
suffices to prove that there exists at least one point $p'\in E_{\rm na}'$
such that $\mathcal{R}_h'$ is analytic at $h^{(-1)}(p')\times
\overline{p'}$. It is convenient to choose a ``good'' such point $p_1'$
which is geometrically well located, namely it belongs to $E_{\rm
na}'$ and in a neighborhood of $p_1'$, the closed set $E_{\rm na}'$
is not too pathological or wild: it lies {\it behind}\, a smooth
generic ``wall'' $M_1'$.

\subsection*{5.3.~Construction of a generic wall} 
As in Lemma~2.3, this point $p_1'$ will belong to a generic
one-codimensional submanifold $M_1'\subset M'$, a kind of ``wall'' in
$M'$ dividing $M'$ locally into two open sides, which will be disposed
conveniently in order that one open side of the ``wall'', say
${M_1'}^-$, will {\it contain only points where $\mathcal{R}_h'$ is
already real analytic}. To show the existence of such a point $p_1'\in
E_{\text{\rm na}}'$ and of such a manifold (``wall'') $M_1'$, we shall
proceed similarly as in [MP1,~Lemma~2.3]. The following picture summarizes
how we proceed intuitively speaking.

\bigskip
\begin{center}
\input figure1.pstex_t
\end{center}
\bigskip

\def\thelemma{5.4}\begin{lemma}
There is a point $p_1'\in E_{\text{\rm na}}'$ and a real analytic
generic hypersurface $M_1'\subset M'$ passing
through $p_1'$ so that $E_{\text{\rm na}}'\backslash\{p_1'\}$ lies near $p_1'$
in one side of $M_1'$ $(${\rm see}
{\sc Figure~4}$)$.
\end{lemma}
 
\proof
Let $q'\in E_{\rm na}'\neq \emptyset$ be an arbitrary point and let
$\gamma'$ be a piecewise real analytic curve running in complex
tangential directions to $M'$ (CR-curve) linking $q'$ with another
point $p'\in M'\backslash E_{\rm na}'$. Such a curve $\gamma'$ exists
because $M'$ and $M'\backslash E_{\rm na}'$ are globally minimal by
assumption (in fact, every open subset of $M'$ is globally minimal,
because $M'$ is locally minimal at every point). After shortening
$\gamma'$, we may suppose that
$\gamma'$ is a smoothly embedded segment, that $p'$ and
$q'$ belong to $\gamma'$ and
are close to each other. Therefore $\gamma'$ can
be described as a part of an integral curve of some nonvanishing real
analytic CR vector field ({\it i.e.} a section of $T^cM'$) $L'$ defined in a
neighborhood of $p'$.

Let $H'\subset M'$ be a small $(2n-2)$-dimensional
real analytic hypersurface passing through $p'$ and transverse to
$L'$. Integrating $L'$ with initial values in $H'$ we obtain
real analytic coordinates $(u',v')\in\R\times\R^{2n-2}$ so
that for fixed $v_0'$, the segments $(u',v_0')$ are contained in the
trajectories of $L'$. After a translation, we may assume that the origin $(0,0)$
corresponds to a point of $\gamma'$ close to $p'$ which is not contained in
$E_{\text{\rm na}}'$, again denoted by $p'$. 
Fix a small $\varepsilon>0$ and for real $\delta\geq 1$, 
define the ellipsoids ({\it see} again {\sc Figure~1} above)
$$
Q_\delta':=\{(u',v'):|u'|^2/\delta+|v'|^2<\varepsilon\}.
$$
There is a minimal $\delta_1>1$ with $\overline{Q_{\delta_1}'}\cap
E_{\text{\rm na}}'\not=\emptyset$. Then $ \overline{Q_{\delta_1}'}\cap
E_{\text{\rm na}}'=\partial Q_{\delta_1}'\cap E_{\text{\rm na}}'$ and
${Q_{\delta_1}'}\cap E_{\text{\rm na}}'=\emptyset$. Observe that every
boundary $\partial Q_{\delta}'$ is transverse to the trajectories of
$L'$ out off the equatorial set
$\Upsilon':=\{(0,v'):|v'|^2=\varepsilon\}$ which is contained in
$M'\backslash E_{\text{\rm na}}'$. Hence $\partial Q_{\delta_1}'$ is
transverse to $L'$ in all points of $\partial Q_{\delta_1}'\cap
E_{\text{\rm na}}'$. So $\partial Q_{\delta_1}'\backslash \Upsilon'$
is generic in $\C^{m+n}$, since $L'$ is a CR field.

To conclude, it suffices to choose a point $p_1'\in\partial Q_{\delta_1}'\cap
E_{\rm na}'$ and to take for $M_1'$ a small real analytic
hypersurface passing through $p_1'$ which is tangent to 
$\partial Q_{\delta_1}'$ at $p_1'$ and satisfies
$M_1'\backslash \{p_1'\} \subset Q_{\delta_1}'$.
\endproof

In summary, it suffices now for our purposes to establish the
following assertion.

\def\thetheorem{5.5}\begin{theorem}
Let $p_1'\in E_{\text{\rm na}}'$ and assume that there exists a
real analytic one-codimensional submanifold $M_1'$ with
$p_1'\in M_1'\subset M'$ which is generic in $\C^n$ such that
$E_{\text{\rm na}}' \backslash \{p_1'\}$ is completely contained in
one of the two open sides of $M'$ divided by $M_1'$ at $p_1'$, say in
${M_1'}^+$, and such that $\mathcal{R}_h'$ is analytic at the points 
$h^{-1}(q')\times \overline{q'}$, for every point $q'$ belonging 
to the other side ${M_1'}^-$. Then the reflection function
$\mathcal{R}_h'$ extends holomorphically at the point $h^{-1}(p_1')\times
\overline{p_1'}$.
\end{theorem}

\noindent
By the CR diffeomorphism
assumption, the formal Taylor series of $h$ at $p_1$
induces an invertible formal CR mapping between $(M,p_1)$ and
$(M',p_1')$. It is shown in [Me6,8] that the associated formal
reflection function converges at $p_1\times \overline{p_1'}$ and (as a
corollary) that there exists a local {\it biholomorphic equivalence}\,
from $(M,p_1)$ onto $(M',p_1')$. Consequently, 
it would be possible to suppose,
without loss of generality, that $(M',p_1')=(M,p_1)$ in
Theorem~5.5. However, since the proof would be completely the same
(except in notation), we shall maintain the general hypotheses. In
coordinates $t'$ vanishing at $p_1'$, we can assume that $M'$ is given
by the {\it real}\, equation ${\rm Im}\,w'=\varphi'(z',\bar z',{\rm
Re}\, w')$, {\it i.e.}  $v'=\varphi'(z',\bar z',u')$ if $w':=u'+iv'$,
or equivalently by the {\it complex}\, equation $\bar
w'=\Theta'(\bar z',t')$ with $\Theta'$ converging in the
polydisc $\Delta_{2n-1}(0,\rho')$ and satisfying $\bar w'\equiv
\Theta'(\bar z',z', \overline{\Theta}'(z',\bar z',\bar w'))$. In fact,
given $\varphi'$, the function $\Theta'$ is the unique solution of the
implicit functional equation $w'-\Theta'(\bar z',t')\equiv
2i\,\varphi'(z',\bar z', (w'+\Theta'(\bar z',t'))/2)$. It is convient
to choose the coordinates in order that $T_0M'=\{w'=\bar w'\}$.  Moreover, an
elementary reasoning using only linear changes of coordinates and
Taylor's formula shows that, after a possible deformation of the manifold
$M_1'$ in a new manifold still passing through $p_1'$ which is bent
quadratically in the left side ${M_1'}^-$, we can assume for simplicity
that $M_1'$ is given by the two equations $\bar w'=\Theta'(\bar
z',t')$ and $x_1'=-[{y_1'}^2+\vert z_\sharp'\vert^2+{u'}^2]$, where we
decompose $z_1'=x_1'+iy_1'$ in real and imaginary part and where we
denote $z_\sharp':=(z_2',\dots, z_{n-1}')$. In this notation, the new side
${M_1'}^-$ is given by:
\def\theequation{5.6}\begin{equation}
{M_1'}^-: \ \ \ \{(z',w')\in M': \ x_1'<-[
{y_1'}^2+\vert z_\sharp'\vert^2+{u'}^2]\}.
\end{equation} 
({\it Warning}: For ease of readability, in Figure~5 below, we have
drawn $M_1'$ as if the defining equation of $M_1'$ was equal
to $x_1'=+[ {y_1'}^2+\vert z_\sharp'\vert^2+{u'}^2]$, so
Figure~5 is slightly incorrect.) Shrinking $\rho'$ if necessary, by
Lemma~5.4, we know that $E_{\rm na}'\backslash \{p_1'\}$ is contained
in the right {\it open}\, part ${M_1'}^+\cap \Delta_n(0,\rho')$.  We set
$p_1:=h^{-1}(p_1')$ and $M_1:=h^{-1}(M_1')$. Then $M_1$ is
one-codimensional generic submanifold of $M$ which is only of class
$\mathcal{C}^\infty$, because the CR diffeomorphism $h$ is only of
class $\mathcal{C}^\infty$. The reader may observe that even if we
take the conclusion of the proof of Theorem~5.5 for granted, namely
even if we admit that $\mathcal{R}_h'$ is real analytic at
$h^{-1}(p_1')\times \overline{p_1'}$, it does {\it not}\, follow necessarily
(unless $M'$ is holomorphically nondegenerate) that $h$ is real
analytic ({\it cf.}  Lemma~1.16), so the hypersurface $h^{-1}(M_1')$
is {\it not}\, real analytic in general. Let $D$ be the global
one-sided neighborhood of automatic extendability of CR functions on
$M$ constructed in \S3.6. Let $D_{p_1}\subset D$ be a small local
one-sided neighborhood of $(M,p_1)$. Since we are working
at $p_1$, we shall identify the two notations
$D_{p_1}$ and $D$ in the sequel. By the considerations of \S3.6,
the reflection function $\mathcal{R}_h'$ associated with these
coordinates is already holomorphic in $D\times
\Delta_n(0,\rho')$, shrinking $\rho'>0$ if necessary. Moreover,
$\mathcal{R}_h'$ is also holomorphic at each point $h^{-1}(q')\times
\overline{q'}$, for all $q'$ belonging to ${M_1'}^{-}$ in a
neighborhood of the origin. Using the computation of \S3.36
(especially, equations~(3.40)), we can make this property more
explicit. Let $(\Psi_{q'}')_{q'\in M'}$ denote the family of
biholomorphisms sending $q'\in M'$ to $0$ simply obtained by
translation of coordinates $t'\mapsto t_*':=t'-q'$. The $\Psi_{q'}'$ are
holomorphically parametrized by $q'\in \Delta_n(0,\rho'/2)$.  Let
$\bar w_*'=\Theta_*'(\bar z_*',t_*')$ denote the corresponding equation
of $\Psi_{q'}'(M')$. Let $h_{q'}$ denote the mapping $h-q'$ obtained
by this translation of coordinates, namely 
$h_*(t):=h(t)-q'$. Let $q:=h^{-1}(q')$.  By assuming
that the reflection function extends holomorphically to a neighborhood of
$h^{-1}(q')\times \overline{q'}$ for every point $q'\in {M_1'}^-$, we
mean precisely that each translated reflection function
$\mathcal{R}_{*,h_*}'$ in these coordinates vanishing at $q'$
extends holomorphically to a neighborhood of $q\times 0$. By
Lemma~3.22, this property is invariant under changes of coordinates
fixing $q\times 0$.  However, we need to express this property in
terms of a {\it single}\, coordinate system, for instance in the
system $t'$ vanishing at $p_1'$, and this is not obvious.

\subsection*{5.7.~Holomorphic extendability in a fixed coordinate system}
This part is delicate and we begin with some heuristic explanations.
As presented in \S2 with a slightly different
definition of $E_{\rm na}'$, in the situation of Theorem~1.2 (where 
$M'$ is holomorphically nondegenerate) and
of the corresponding Theorem~5.5, the mapping
$h$ already extends holomorphically to a neighborhood of $M_1^-$ in 
$\C^n$. However, in the situation of Theorems~1.9 and~5.5, this is untrue in 
general. Consequently, we raise the following question: if we fix the
coordinate system $t'$ vanishing at the point $p_1'\in M_1'$ of Theorem~5.5,
is it also true that the components $\Theta_\beta'(h(t))$ extend
holomorphically to a neighborhood of $M_1^-$ in $\C^n$~?
Let $q'\in {M_1'}^-$ be close to the
point $p_1'$, which is the origin in the coordinates
$t'$. Let $t_*':=t'-t_*'$ as
in \S3.36.  Let $\bar w_*'=\Theta_*'(\bar z_*',t_*')$ be the translated
equation of $M'$. Also, denote $h(t)-q'$ by
$h_*(t)$. By assumption, the reflection function
$\bar\mu_*'-\Theta_*'(\bar\lambda_*',h_*(t))$ extends holomorphically to
$\Delta_n(q,\sigma_q)\times \Delta_n(0,\sigma_{q'}')$, for some two
positive real numbers $\sigma_q>0$ and $\sigma_{q'}'>0$. By
Lemma~3.16, we have a Cauchy estimate for the
holomorphic extensions $\theta_{*\beta}'(t)$ of the components 
$\Theta_{*\beta}'(h(t))$ of
the reflection function, say $\vert \theta_{*\beta}'(t)\vert
\leq C\,(\sigma_{q'}')^{-\vert\beta\vert}$ for all $\vert t-q\vert
< \sigma_q$. Possibly, $\sigma_{q'}'$ is smaller than 
$\vert q'\vert$. In the previous coordinate system $t'$, it would be
natural to deduce that the $\mathcal{C}^\infty$-smooth
CR functions $\Theta_\beta'(h(t))$ extend
holomorphically to a neighborhood of $q$ in $\C^n$. Unfortunately, 
by formulas~(3.43), we would necessarily have the
following representation for the desired holomorphic extensions
$\theta_\beta'(t)$
of the components $\Theta_\beta'(h(t))$ of the reflection function
(if the series would be convergent for $\vert t-q\vert
\leq \sigma_q$):
\def\theequation{5.8}\begin{equation}
\theta_\beta'(t):=\theta_{*\beta}'(t)+
\sum_{\gamma\in\N_*^{n-1}}\,(\bar z_{q'}')^\gamma\,
(-1)^\gamma\,\theta_{*\beta+\gamma}'(t)\,
{(\beta+\gamma)!\over \beta!\ \gamma!}.
\end{equation}
For this formulas to converge {\it normally}\, and to define a
holomorphic function of $t$, it would be necessary that the modulus of $\bar
z_{q'}'$ be smaller than $\sigma_{q'}'$, which is {\it not a
priori}\, true in general.  This difficulty is meaningful, unavoidable
and important.

At present, we may nevertheless 
observe a useful trick: {\it if $\bar
z_{q'}'$ vanishes, then formulas~(5.8) automatically yield holomorphic
functions $\theta_\beta'(t)$ in the polydisc $\{\vert t-q\vert <
\sigma_q\}$}. Indeed, if $\bar z_{q'}'=0$, there are no infinite series
at all~! Indeed, 
$\theta_\beta'(t)\equiv \theta_{*\beta}'(t)$.
So choosing a point $q'$ with vanishing coordinate $\bar
z_{q'}'$ is a crucial observation allowing to bypass the
nonconvergency of the series~(5.8). Moreover, we will crucially use
this trick in the proof of Lemma~7.7 (corresponding to Lemma~2.4). In
sum, we have observed that Cauchy estimates might be killed after complex
tangential displacement whereas they are trivially conserved after complex
transversal displacement.

\subsection*{5.9.~Holomorphic extension to a neighborhood of
$M_1^-$} Fortunately, thanks to Artin's approximation theorem, we can
bypass the general difficulty above and we can 
make $\sigma_{q'}'$ larger than $\bar
z_{q'}'$ at the cost of reducing $\sigma_q$. In advance, the following
Lemma~5.10 is adapted to its application in \S7 below. Let
$q_1'\in M'$ be close to $p_1'$, let $t'$ be a fixed system 
of coordinates centered at $q_1'$, so $q_1'$
is identified with the origin. Without loss of
generality, we can assume that $h(M\cap \Delta_n(0,\rho/2))\subset
M'\cap \Delta_n(0,\rho'/2)$. Let $E\subset M\cap
\Delta_n(0,\rho/2)$ be an arbitrary closed subset, not necessarily
passing through the origin.  Set $E':=h(E)$. As in Theorem~5.5, let us assume
that the reflection function centered at points $q\times
\overline{h(q)}$ is locally holomorphically extendable, for all $q\in 
(M\backslash E)\cap \Delta_n(0,\rho/2)$. Then the following holds.

\def\thelemma{5.10}\begin{lemma}
In the fixed system of coordinates $t'$ centered at $q_1'$, there
exists a neighborhood $\Omega$ of $(M\backslash E)\cap
\Delta_n(0,\rho/2)$ in $\C^n$ to which the components
$\Theta_\beta'(h(t))$ of the reflection function extend
holomorphically.
\end{lemma}

\proof
So, 
let $q\in (M\backslash E)\cap \Delta_n(0,\rho/2)$ be an 
arbitrary point and let $q':=h(q)$.  As in Lemma~3.37,
let $t_*':=t'-q'$, let $t_*:=t-q$ and let
$\bar\nu_*'-\sum_{\beta\in\N^{n-1}}\,
(\bar\lambda_*')^\beta\,\Theta_{*\beta}'(h_*(t_*))$ be the reflection
function centered at $q\times \overline{q'}$. Here, we have $\vert
q\vert < \rho/2$ and $\vert q'\vert <\rho'/2$. Let $\sigma_q>0$ and
$\sigma_{q'}'>0$ be such that $\mathcal{R}_h'(t, \bar\nu')$ extends as
a holomorphic function
$R_*'(t_*,\bar\nu_*'):=\bar\mu_*'-\sum_{\beta\in\N^{n-1}}\,
(\bar\lambda_*')^\beta\,\theta_{*\beta}'(t_*)$ for $\vert t_*\vert <
\sigma_q$ and $\vert \bar\nu_*'\vert < \sigma_{q'}'$. Of course, it
follows that the holomorphic functions $\theta_{*\beta}'(t_*)$
converge for $\vert t_*\vert < \sigma_q$ and that if $H_*(t_*)$
denotes the formal power series of $h_*(t_*)$ at $t_*=0$, then
$\Theta_{*\beta}'(H_*(t_*))\equiv \theta_{*\beta}'(t_*)$ in $\C\dl
t_*\dr$. By Artin's approximation theorem, there exists a holomorphic
mapping $\mathcal{H}_*(t_*)$ defined for $\vert t_*\vert<
\sigma_*<\sigma_q$ such that
$\Theta_{*\beta}'(\mathcal{H}_*(t_*))\equiv \theta_{*\beta}'(t_*)$ in
$\C\{t_*\}$. By the Cauchy estimates for $\Theta_\beta'(t')$, since
$\vert q'\vert < \rho'/2$, there exists a constant $C>0$ such that we
have $\vert \Theta_{*\beta}'(t_*')\vert \leq C\,
(\rho'/2)^{-\vert\beta\vert}$ for all $\vert t_*'\vert <
\rho'/4$. Shrinking $\sigma_*$ if necessary, we can assume that $\vert
\mathcal{H}_*(t_*)\vert < \rho'/4$ for all $\vert t_*\vert
<\sigma_*$. It follows that
\def\theequation{5.11}\begin{equation}
\vert\Theta_{*\beta}'(\mathcal{H}_*(t_*))\vert
=\vert\theta_{*\beta}'(t_*)\vert < C\, (\rho'/2)^{-\vert\beta\vert},
\end{equation}
for all $\beta\in\N^{n-1}$. Finally, this Cauchy estimate is
appropriate to deduce that the series defined in equations~(5.8) do
converge normally and do define holomorphic extension to the polydisc
$\Delta_n(q,\sigma_*)$ of the components of the reflection function
centered at $q_1\times\overline{q_1'}$.  For all $q$ running in
$(M\backslash E)\cap \Delta_n(0,\rho/2)$, the various obtained
extensions of course stick together thanks to the uniqueness principle
at the boundary ([P1]). The proof of Lemma~5.10 is complete.
\endproof

In particular, in the situation of Theorem~5.5, it follows from
Lemma~5.10 (with $q_1:=p_1$) that we can assume that the components of
the reflection function centered at $p_1\times \overline{p_1'}$ extend
holomorphically to a neighborhood of $(M\cap E_{\rm na})\cap
\Delta_n(0,\rho/2)$.  Now we can begin our principal geometric
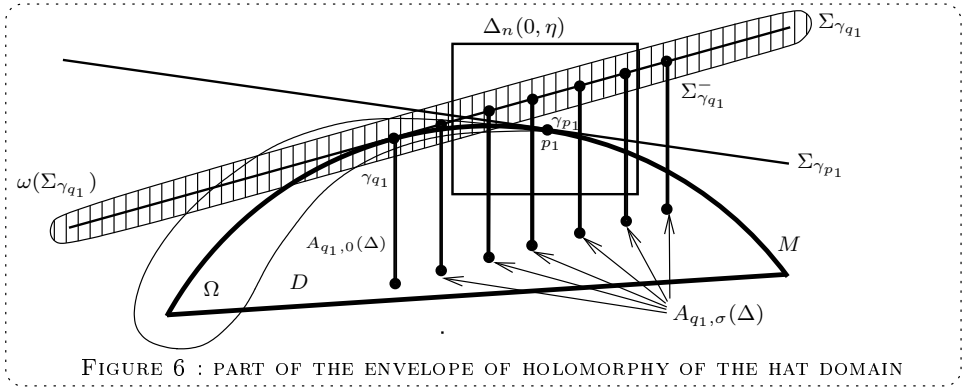
constructions. As explained in \S2.2, we intend to study the envelope
of holomorphy of the union of $D\cup \Omega$ together with an
arbitrary thin neighborhood of a Levi-flat hypersurface
$\Sigma_\gamma$. We need real arcs and analytic discs.

\section*{\S6. Envelopes of holomorphy of domains with Levi-flat hats}

\subsection*{6.1.~A family of real analytic arcs}
To start with, we choose coordinates $t$ and $t'$ as above near $M$
and near $M'$ in which $p_1:=h^{-1}(p_1')$ and $p_1'$ are the origin
and in which the complex equations of $M$ and of $M'$ are given by
$\bar w=\Theta(\bar z, t)$ and $\bar w'=\Theta'(\bar z', t')$
respectively. Geometrically speaking, it is convenient to assume
$T_0M=\{{\rm Im}\, w=0\}$ and $T_0M'=\{{\rm Im}\, w'=0\}$.  We shall
denote the real equations of $M$ and
of $M'$ by $v=\varphi(z,\bar z,u)$ and
$v'=\varphi'(z',\bar z',u')$ respectively. We assume that the power
series defining $\Theta$ and $\Theta'$ converge normally in the
polydisc $\Delta_{2n-1}(0,\rho)$ and $\Delta_{2n-1}(0,\rho')$
respectively. For $q_1'\in M'$ close to the origin in the target
space, we now consider a convenient, sufficiently rich, family of
embedded real analytic arcs $\gamma_{q_1'}'(s')$, depending on
$(2n-1)$ very small real parameters $(z_{q_1'}',u_{q_1'}')\in \C^{n-1}
\times \R$ satisfying $\vert z_{q_1'}'\vert < \varepsilon'$, $\vert
u_{q_1'}'\vert < \varepsilon'$, where $\varepsilon' < < \rho'$, with
the ``time parameter'' $s'$ satisfying $\vert s'\vert \leq \rho'/2$,
which are all transverse to the complex tangential directions of $M'$,
and which are defined as follows:
\def\theequation{6.2}\begin{equation}
\left\{
\aligned
{}
&
\gamma_{q_1'}':=\left\{(x_{1,q_1'}'-{s'}^2-(y_{1,q_1'}'+s')^2-
\vert z_{\sharp q_1'}'\vert^2-{u_{q_1'}'}^2+i[y_{1,q_1'}'+s'],\right.\\
&
\ \ \ \ \ \ \ \ \ \ \ \ \ \ \ \ \ \ \ \ \ \ \ \ \ \ \ \ \ \
\left.
,z_{\sharp q_1'}', \, u_{q_1'}') \in M'\,:
\ s'\in \R, \, \vert s' \vert \leq \rho'/2\right\}.
\endaligned\right.
\end{equation}
Here, in the definition of $\gamma_{q_1'}'$, we identify a point of
$M'$ with its $(2n-1)$ real coordinates $(z',u')=(x_1'+iy_1',
z_\sharp',u')$. We also
recall that $z_\sharp'=(z_2',\dots,z_{n-1}')$ and that ${M_1'}^-$ is given
by~(5.6).  The following figure, in which we have reversed the
curvature of $M_1'$ for easier readability, explains how the
$\gamma_{q_1'}'$ and $M_1'$ are disposed.

\bigskip
\begin{center}
\input figure2.pstex_t
\end{center}
\bigskip

We identify the arcs $\gamma_{q_1'}'$ with the mappings $s'\mapsto
\gamma_{q_1'}'(s')$.  It can be straightforwardly checked that the
following properties hold:
\begin{itemize}
\item[{\bf (1)}]
{\it The mapping $(z_{q_1'}',u_{q_1'}')\mapsto
\gamma_{q_1'}'(0)$ is a real analytic diffeomorphism onto a
neighborhood of $0$ in $M'$. Furthermore, the inverse image of $M_1'$
and of ${M_1'}^-$ simply correspond to the sets $\{x_{1,q_1'}'=0\}$ and
$\{x_{1,q_1'}'<0\}$, respectively.}
\item[{\bf (2)}]
{\it For $x_{1,q_1'}'<0$, we have $\gamma_{q_1'}'\subset \subset
{M_1'}^-$.}
\item[{\bf (3)}]
{\it For $x_{1,q_1'}'=0$, we have 
$\gamma_{q_1'}'\cap M_1'=\{\gamma_{q_1'}'(0)\}$.}
\item[{\bf (4)}]
{\it  For $x_{1,q_1'}'=0$, the order of contact of $\gamma_{q_1'}'$ with 
$M_1'$ at the point $\gamma_{q_1'}'(0)$ equals $2$.}
\item[{\bf (5)}]
For all $\vert z_{q_1'}'\vert, \vert u_{q_1'}'\vert<\varepsilon'$, we have
$\gamma_{q_1'}'([-\rho'/2,-\rho'/4])\subset {M_1'}^-$
and $\gamma_{q_1'}'([\rho'/4,\rho'/2])\subset {M_1'}^-$.
\end{itemize} 

\subsection*{6.3.~Inverse images} Since $h$ is a $\mathcal{C}^\infty$-smooth
CR diffeomorphism, by inverse image, we get in $M$ a family of
$\mathcal{ C}^\infty$-smooth arcs, namely $h^{-1}(\gamma_{q_1'}')$. In
analogy with the notation $\gamma_{q_1'}'(s')$, we shall denote these
arcs by $\gamma_{q_1}(s)$. By the index notation $\cdot_{q_1}$, we
mean that these arcs are parametrized by the point
$q_1:=h^{-1}(q_1')\in M$. Of course, a point $q_1\in M$ can be
identified with its coordinates
$(z_{q_1},u_{q_1})\in\C^{n-1}\times\R$, so the arcs $\gamma_{q_1}$ are
concretely parameterized by $(z_{q_1},u_{q_1})\in\C^{n-1}\times\R$ and
by the ``time'' $s\in\R$. It is convenient to identify the point $p_1$
with the origin (in the coordinate system $t$) and the point $q_1$
close to $p_1$ with its coordinates $(z_{q_1},u_{q_1})$. Of course,
shrinking a bit $\rho$ if necessary, there exists $\varepsilon < <
\rho$ such that the parameters satisfy $\vert z_{q_1}\vert <
\varepsilon$, $\vert u_{q_1}\vert < \varepsilon$ and $\vert s\vert
\leq \rho/2$.  Evidently, the $\mathcal{C}^\infty$-smooth arcs
$\gamma_{q_1}$ satisfy four properties similar to {\bf (1--5)} above
with respect to $M_1$. Let us summarize the geometric properties 
that will be of important use later, 
when envelopes of holomorphy will appear on scene. 

\def\thelemma{6.4}\begin{lemma}
For all small $x_{1,q_1}<0$ and $z_{\sharp,q_1}, u_{q_1}$ arbitrary,
the following two properties holds:
\begin{itemize}
\item[{\bf (1)}]
The center points
$\gamma_{q_1}(0)$ of the smooth arcs
$\gamma_{q_1}$ cover diffeomorphically the left negative one-sided
neighborhood $M_1^-$ of $M_1$ in a neighborhood of $p_1$.
\item[{\bf (2)}]
The arcs $\gamma_{q_1}$ are entirely contained in $M_1^-$ and satisfy,
$\gamma_{q_1}([-\rho/2,-\rho/4])\subset M_1^-$ and
$\gamma_{q_1}([\rho/4,\rho/2])\subset M_1^-$,
{\rm even for small} $x_{1,q_1}\geq 0$.
\end{itemize} 
\end{lemma}

\subsection*{6.5.~Construction of a family of Levi-flat hats}
Next, if $\gamma$ is a $\mathcal{C}^\infty$-smooth arc in $M$
transverse to $T^cM$ at each point, we can construct the union of
Segre varieties attached to the points running in $\gamma$:
$\Sigma_{\gamma}:= \bigcup_{q\in \gamma} S_{\bar q}$. We remind that
the Segre variety $S_{\bar q}$ associated to an arbitrary point $q$
close to the origin is the complex hypersurface of $\Delta_n(0,\rho)$
of equation $w=\overline{\Theta}(z,\bar t_q)$. For various arcs
$\gamma_{q_1}$, we obtain various sets $\Sigma_{\gamma_{q_1}}$ which
are in fact $\mathcal{C}^\infty$-smooth Levi-flat hypersurfaces in a
neighborhood of $\gamma_{q_1}$. The uniformity of the size of such
neighborhoods follows immediately from the smooth dependence with
respect to $(z_{q_1},u_{q_1})$. Shrinking $\rho$ if necessary, the
Levi-flat hypersurface $\Sigma_{\gamma_{q_1}}$ is closed in
$\Delta_n(0,\rho/3)$. What we shall need in the sequel can then be
summarized as follows.

\def\thelemma{6.6}\begin{lemma}
There exists $\varepsilon>0$ with $\varepsilon < < \rho$ such that, if
the parameters of $\gamma_{q_1}$ satisfy $\vert z_{q_1}\vert, \vert
u_{q_1}\vert < \varepsilon$, then the set
$\Sigma_{\gamma_{q_1}}\cap\Delta_n(0,\rho/3)$ is a closed
$\mathcal{C}^\infty$-smooth $($and $\mathcal{C}^\infty$-smoothly
parametrized$)$ \text{\rm Levi-flat} hypersurface of
$\Delta_n(0,\rho/3)$.
\end{lemma}

\subsection*{6.7.~Two families of half-attached analytic discs}
Let us now define inverse images of analytic discs. Complexifying the
real analytic arcs $\gamma_{q_1'}'$, we obtain local transverse
holomorphic discs $(\gamma_{q_1'}')^c$, closed in
$\Delta_n(0,\rho/2)$, of which one half part penetrates inside
$D':=h(D)$. Uniformly smoothing out the corners of such half discs
({\it see}\, the right hand side of {\sc Figure~3}),
using Riemann's conformal mapping theorem and then an automorphism of
$\Delta$, we can easily construct a real analytically parameterized
family of analytic discs $A_{q_1'}': \Delta\to \C^n$ which are
$\mathcal{C}^\infty$-smooth up to the boundary $b\Delta$ and 
are embedding of $\overline{\Delta}$ such that, if
we denote $b^+\Delta:= b\Delta\cap \{\text{\rm Re} \, \zeta \geq 0\}$
(and $b^-\Delta:= b\Delta\cap \{\text{\rm Re} \, \zeta \leq 0\}$),
then we have $A_{q_1'}'(1)= \gamma_{q_1'}'(0)$ and also:
\def\theequation{6.8}\begin{equation}
\gamma_{q_1'}'\cap \Delta_n(0,\rho'/4) \ \subset \
A_{q_1'}'(b^+\Delta) \ \subset \ \gamma_{q_1'}'\cap \Delta_n(0,\rho'/3),
\end{equation}
for all $\vert z_{q_1'}'\vert, \vert u_{q_1'}'\vert < \varepsilon'$
({\it cf.}~{\sc Figures}~2 and~3).  Consequently, the composition with
$h^{-1}$ yields a family of analytic discs $A_{q_1}(\zeta):=
h^{-1}(A_{q_1'}'(\zeta))$ which satisfy similar properties, namely:

\def\thelemma{6.9}\begin{lemma}
The mapping $(q_1,\zeta)\mapsto A_{q_1}(\zeta)$ is
$\mathcal{C}^\infty$-smooth and it provides a uniform
family of $\mathcal{C}^\infty$-smooth embeddings 
of the closed unit disc $\overline{\Delta}$ into $\C^n$.
Furthermore, we have $A_{q_1}(1)=\gamma_{q_1}(0)$ and
\def\theequation{6.10}\begin{equation}
\gamma_{q_1}\cap \Delta_n(0,\rho/4) \ \subset \
A_{q_1}(b^+\Delta) \ \subset \ \gamma_{q_1}\cap \Delta_n(0,\rho/3).
\end{equation}
Finally, we have $A_{q_1}(b^-\Delta) \subset\subset D\cup M_1^-$.
\end{lemma}

This family $A_{q_1}$ will be our starting point to study the envelope
of holomorphy of (a certain subdomain of) the union of $D$ together
with a neighborhood $\Omega$ of $M_1^-$ and with an arbitrarily thin
neighborhood $\omega(\Sigma_{\gamma_{q_1}})$ of $\Sigma_{\gamma_{q_1}}$ ({\it
see} {\sc Figure}~3 and {\sc Figure}~6 below).  At first, we must
include $A_{q_1}$ into a larger family of discs obtained by sliding
the half-attached part inside $\Sigma_{\gamma_{q_1}}$ along the
complex tangential directions of $\Sigma_{\gamma_{q_1}}$.

\subsection*{6.11.~Deformation of half-attached analytic discs}
To this aim, we introduce the equation $v=H_{q_1}(z,u)$ of
$\Sigma_{\gamma_{q_1}}$, where the mapping $(q_1,z,u)\mapsto H_{q_1}(z,u)$
is of course $\mathcal{C}^\infty$-smooth and $\vert\vert
H_{q_1}-H_{p_1}\vert\vert_{C^\infty(z,u)}$ is very small.  Further, we
need some formal notation. We denote $A_{q_1}(\zeta):=
(z_{q_1}(\zeta),u_{q_1}(\zeta))$ and $A_{q_1}(1)=\gamma_{q_1}(0)=:
(z_{q_1}^1,u_{q_1}^1)$. To deform these discs by applying the
classical works on analytic discs and because Banach spaces are
necessary, we shall work in the regularity class
$\mathcal{C}^{k,\alpha}$, where $k\geq 1$ is arbitrary and where
$0<\alpha<1$, which is sufficient for our purposes.  Let $T_1$ denote
the Hilbert transform vanishing at $1$ ({\it see} [Tu1,2,3], [MP1,2]).
By definition, $T_1$ is the unique (bounded, by a classical result)
endomorphism of the Banach space $\mathcal{
C}^{k,\alpha}(b\Delta,\R)$, $0<\alpha<1$, to itself such that
$\phi+iT_1(\phi)$ extends holomorphically to $\Delta$ and $T_1\phi$
vanishes at $1\in b\Delta$, {\it i.e.}  $(T_1(\phi))(1)=0$. Our next
reasoning below is similar to the one in A\u\i rapetyan [A\u\i]: 
we shall ``translate'' a small analytic
disc which is attached to a pair of transverse hypersurfaces. We know
that the disc $A_{q_1}$ has one half of its boundary attached to the
smooth hypersurface $v=H_{q_1}(z,u)$. After a possible linear change
of coordinates, we can assume that the other half is attached to
another real hypersurface $\Lambda_{q_1}$ of equation $v=G_{q_1}(z,u)$
smoothly depending on the parameter $q_1$.  Indeed, since the half
disc is transverse to $\Sigma_{\gamma_{q_1}}$ along $b^+\Delta$ and an
embedding of $\overline{\Delta}$ into $\C^n$, there exist infinitely
many such hypersurfaces $\Lambda_{q_1}$.  Furthermore, we can asssume
that $A_{q_1}$ sends neighborhoods of $i$ and $-i$ in $b\Delta$ into
the intersection of the two hypersurfaces $\Sigma_{\gamma_{q_1}}\cap
\Lambda_{q_1}$.  Let $\varphi^-$ and $\varphi^+$ be two
$\mathcal{C}^\infty$-smooth functions on $b\Delta$ satisfying
$\varphi^-\equiv 0$, $\varphi^+\equiv 1$ on $b^+\Delta$ and
$\varphi^-+\varphi^+=1$ on $b\Delta$. The fact that our disc is half
attached to $\Sigma_{\gamma_{q_1}}$ and half attached to
$\Lambda_{q_1}$ can be expressed by saying that
\def\theequation{6.12}\begin{equation}
v_{q_1}(\zeta)=\varphi^+(\zeta) \, H_{q_1}(
z_{q_1}(\zeta),u_{q_1}(\zeta))+\varphi^-(\zeta) \ G_{q_1}(
z_{q_1}(\zeta),u_{q_1}(\zeta)),
\end{equation}
for all $\zeta\in b\Delta$.  Since the two functions $u_{q_1}$ and
$v_{q_1}$ on $b\Delta$ are harmonic conjugates, the following (Bishop)
equation is satisfied on $b\Delta$ by $u_{q_1}$:
\def\theequation{6.13}\begin{equation}
u_{q_1}(\zeta)=-\left[T_1 \left(\varphi^+ \, H_{q_1}(z_{q_1},\, 
u_{q_1})+\varphi^- \, G_{q_1}(z_{q_1},\, 
u_{q_1})
\right)\right](\zeta).
\end{equation} 
We want to perturb these discs $A_{q_1}$ by ``translating'' them along
the complex tangential directions to
$\Sigma_{\gamma_{q_1}}$. Introducing a new parameter
$\sigma\in\C^{n-1}$ with $\vert \sigma \vert < \varepsilon$, we can
indeed include the discs $A_{q_1}$ into a larger parameterized family
$A_{q_1,\sigma}$ by solving the following perturbed Bishop equation on
$b\Delta$ with parameters $(q_1,\sigma)$:
\def\theequation{6.14}\begin{equation}
u_{q_1,\sigma}(\zeta)=-\left[T_1 \left(\varphi^+ \, 
H_{q_1}(z_{q_1}+\sigma,\, 
u_{q_1,\sigma})
+\varphi^- \, G_{q_1}(z_{q_1}+\sigma,\, 
u_{q_1,\sigma})
\right)\right](\zeta).
\end{equation} 
For instance, the existence and the $\mathcal{C}^{k,\beta}$-smoothness
(with $0<\beta<\alpha$ arbitrary) of a solution $u_{q_1,\sigma}$
to~\thetag{6.14} follows from Tumanov's work [Tu3]. Clearly the
solution disc $A_{q_1,\sigma}$ is half attached to
$\Sigma_{\gamma_{q_1}}$. Differentiating the mapping $\C^{n-1}\times
b^+\Delta\ni(\sigma,\zeta)\,\mapsto\,
(z_{p_1}(\zeta)+\sigma,u_{p_1,\sigma}(\zeta))\in\Sigma_{\gamma_{p_1}}$
at the point $0\times 1$, using the fact that $A_{p_1}(b^+\Delta)$ is
tangent to the $u$-axis at $p_1$ (since $\gamma_{p_1}$ is tangent to
the $u$-axis at $p_1$) which gives $[{d\over
d\theta}(z_{p_1}(e^{i\theta}))]_{\theta=0}=0$, we obtain easily
property {\bf (3)} of the next statement. Notice that since the
discs are $\mathcal{C}^{k,\beta}$ for all $k$, and since the solution
$u_{q_1,\sigma}$ of the modified Bishop equation~(6.14) is the same in
$\mathcal{C}^{k,\beta}$ and in $\mathcal{C}^{l,\beta}$, the discs
$A_{q_1,\sigma}$ are in fact of class $\mathcal{C}^\infty$ with
respect to all the variables.

\def\thelemma{6.15}\begin{lemma}
The $\mathcal{C}^{\infty}$-smooth deformation $($``translation'' type$)$ 
of analytic discs $(q_1,\sigma,\zeta)\mapsto A_{q_1,\sigma}(\zeta)$ is
defined for $\vert q_1 \vert < \varepsilon$ and
for $\vert \sigma \vert < \varepsilon$ and satisfies the following three
properties:
\begin{itemize}
\item[{\bf (1)}]
$A_{q_1,0}\equiv A_{q_1}$.
\item[{\bf (2)}]
$A_{q_1,\sigma}(b^+\Delta)\subset \Sigma_{\gamma_{q_1}}$ for all
$\sigma$.
\item[{\bf (3)}]
The mapping $\C^{n-1}\times b^+\Delta\ni (\sigma,\zeta) \mapsto
A_{p_1,\sigma}(\zeta)\in \Sigma_{\gamma_{p_1}}$ is a local
$\mathcal{C}^\infty$ diffeomorphism from a neighborhood of $0\times 1$
onto a neighborhood of $A_{p_1}(1)=p_1$ in $\Sigma_{\gamma_{p_1}}$.
\end{itemize}
\end{lemma}

\subsection*{6.16.~Preliminary to applying the continuity principle} 
At first, we shall let the parameters $(q_1,\sigma)$ range in certain
new precise subdomains. We choose a positive $\delta < \varepsilon$
with the property that the range of the mapping in {\bf (3)} above, when
restricted to $\{\vert\sigma\vert < \delta\}\times b^+\Delta$, {\it
contains} the intersection of $\Sigma_{\gamma_{p_1}}$ with a small
polydisc $\Delta_n(0,2\eta)$, for some $\eta>0$. Recall that $p_1$
is identified with the origin $0\in\C^n$. Of course, there
exists a constant $c>1$, depending only on the Jacobian matrix of the
mapping in {\bf (3)} at $0\times 1$ such that $c^{-1}\, \delta \leq \eta
\leq c\, \delta$. Let $\Delta(1,\delta)$ denote the disc of radius
$\delta$ centered at $1\in\C$. Furthermore, since the boundary of the
disc $A_{p_1,0}$ is transversal to $T_0^c \Sigma_{\gamma_{p_1}}$, then
after shrinking a bit $\eta$ if necessary, we can assume that the set
$\{A_{p_1,\sigma}(\zeta)\,: \vert \sigma \vert < \delta, \ \zeta\in
\Delta\cap \Delta(1,\delta)\}$ {\it contains and foliates by half
analytic discs the whole lower side} $\Delta_n(0,2\eta)\cap
\Sigma_{\gamma_{p_1}}^-$ ({\it see}~{\sc Figure}~6).  Of course, the
side $\Sigma_{\gamma_{p_1}}^-$ is ``the same side'' as $M^-$, {\it
i.e.} the side of $\Sigma_{\gamma_{p_1}}$ where the greatest portion
of $D$ lies. However, $D$ is in general {\it not} entirely contained
in $\Sigma_{\gamma_{p_1}}^-$, because the Segre varieties $S_{\bar q}$
for $q\in \gamma_{p_1}$ may well intersect $D$.

As presented in \S2, we now fix a neighborhood $\Omega$ of $M_1^-$ in
$\C^n$ to which all the components of the reflection function extend
holomorphically. Such a neighborhood is provided by Lemma~5.10
above. Let $\omega(\Sigma_{\gamma_{q_1}})$ be an arbitrary
neighborhood of $\Sigma_{\gamma_{q_1}}$ in $\C^n$. Our goal is to show
that the envelope of holomorphy of $\Omega\cup D\cup
\omega(\Sigma_{\gamma_{q_1}})$ contains at least the lower side
$\Delta_n(0,\eta) \cap\Sigma_{\gamma_{q_1}}^-$ for all $q_1$ small
enough. We shall apply this
to the components of the reflection function 
in \S7 below.  

By construction, the half parts $A_{q_1,\sigma}(b^+\Delta)$
are all contained in $\Sigma_{\gamma_{q_1}}$. It remains now to
control the half parts $A_{q_1,\sigma}(b^-\Delta)$. Using the last
property of Lemma~6.9, namely $A_{q_1,0}(b^-\Delta)\subset\subset
D\cup M_1^-$, it is clear that, after shrinking $\delta$ if necessary,
then we can insure that $A_{q_1,\sigma}(b^-\Delta) \subset \subset
D\cup \Omega$ for all $\vert q_1\vert < \varepsilon$ and all $\vert
\sigma \vert < \delta$. Of course, this shrinking will result in a
simultaneous shrinking of $\eta$, and we still have the important
inclusion relation: $\{A_{p_1,\sigma}(\zeta)\,: \vert \sigma \vert <
\delta, \ \zeta\in \Delta\cap \Delta(1,\delta)\} \supset
\Delta_n(0,2\eta)\cap \Sigma_{\gamma_{p_1}}^-$.  Finally, shrinking
again $\varepsilon$ if necessary, we then come to a situation that we
may summarize:
\def\thelemma{6.17}\begin{lemma}
{\it For all $\vert q_1\vert < \varepsilon$ and $\vert\sigma\vert<\delta$, 
we have}:
\def\theequation{6.18}\begin{equation}
\left\{
\aligned
{} &
\{A_{q_1,\sigma}(\zeta)\,:
\vert \sigma \vert < \delta, \ \zeta\in \Delta\cap \Delta(1,\delta)\} 
\ \supset \ \Delta_n(0,\eta)\cap \Sigma_{\gamma_{q_1}}^-,\\
&
A_{q_1,\sigma}(b^+\Delta)\subset \Sigma_{\gamma_{q_1}} \ \ {\rm and} \ \ 
A_{q_1,\sigma}(b^-\Delta)\subset \subset D\cup \Omega,
\endaligned\right.
\end{equation}
\end{lemma}
Shrinking $\varepsilon$ if necessary we can also insure that the
intersection of $D$ with $\Delta_n(0,\eta)\cap
\Sigma_{\gamma_{q_1}}^-$ is {\it connected} for all $\vert
q_1\vert < \varepsilon$. Implicitely, 
we assume that $\varepsilon < < \delta$, hence also $\varepsilon < < \eta$.

\subsection*{6.19.~Envelopes of holomorphy}
We are now in position to state and to prove the main assertion of this
paragraph. Especially, the following lemma will be applied to 
each member of the collection $\{\Theta_\beta'(h(t))\}_{\beta\in\N^{n-1}}$
in \S7 below.

\bigskip
\begin{center}
\input figure7.pstex_t
\end{center}
\bigskip

\def\thelemma{6.20}\begin{lemma}
Let $\delta,\, \eta,\, \varepsilon>0$ as above, namely satisfying
$\delta\simeq\eta$, $\varepsilon < < \delta$ and
$\varepsilon < < \eta$. If $\delta>0$ is sufficiently small, 
then the following holds. If a holomorphic
function $\psi\in\mathcal{O}(D\cup \Omega)$ extends holomorphically to
a neighborhood $\omega(\Sigma_{\gamma_{q_1}})$ in $\C^n$, then there exists
a unique holomorphic function $\Psi\in \mathcal{O}(D\cup
[\Delta_n(0,\eta)\cap \Sigma_{\gamma_{q_1}}^-])$ such that
$\Psi\vert_D\equiv \psi$.
\end{lemma}

\proof This is an application of the Behnke-Sommer Kontinuit\"atssatz
({\it see} {\sc Figure}~6). Let $q_1$ with $\vert q_1\vert <
\varepsilon$.  We shall explain later how we choose $\delta>0$
sufficiently small.  Let $\psi\in \mathcal{O}(D\cup\Omega)$. By
assumption, there exists a holomorphic function
$\psi_\omega\in\mathcal{O} (\omega(\Sigma_{\gamma_{q_1}}))$ such that
$\psi_\omega=\psi$ in a neighborhood of
$\gamma_{q_1}([-\rho/2,\rho/2])\cap \Delta_n(0,\rho/3)$ 
in $\C^n$.  First of all, we must
construct a domain $B_{q_1}\subset D\cup \Omega$ sufficiently large
such that $\psi$ and $\psi_\omega$ stick together in a unique
holomorphic function defined in the union $B_{q_1}\cup
\omega(\Sigma_{\gamma_{q_1}})$. To get this extension property, we
need that $B_{q_1}\cap \omega(\Sigma_{\gamma_{q_1}})$ be
connected. For this to hold, we construct (equivalently, we shrink)
the neighborhood $\omega(\Sigma_{\gamma_{q_1}})$ as a union of
polydiscs of very small constant radius centered at points of
$\Sigma_{\gamma_{q_1}}$. Next, we construct in two parts $B_{q_1}$ as
follows.  The first part of $B_{q_1}$ consists of a small neighborhood
of $A_{q_1,0}(\overline{\Delta})$ in $\C^n$, for instance a union of
small polydiscs centered at points of $A_{q_1,0}(\overline{\Delta})$
which are of constant very small radius in order to be contained in
$D\cup\Omega$.  The second part of $B_{q_1}$ consists of three
subparts, namely the union of polydiscs of radius $2\delta$ centered
at points of $A_{p_1}(b^-\Delta)$, at points of
$A_{p_1}(b^+\Delta)\cap \gamma_{p_1}([-\rho/2,-\rho/4])$ and at points
of $A_{p_1}(b^+\Delta)\cap\gamma_{p_1}([\rho/4,\rho/2])$. This part is
the same for all $B_{q_1}$. By Lemma~6.4{\bf(2)}, if $\delta$ is small
enough, the second part of $B_{q_1}$ will be contained in
$D\cup\Omega$. This is how we choose $\delta>0$ small
enough. Moreover, because $A_{p_1}(\overline{\Delta})$ are
non-tangentially half-attached to $\Sigma_{\gamma_{p_1}}$ along
$b^+\Delta$, the intersection $B_{q_1}\cap
\omega(\Sigma_{\gamma_{q_1}})$ is connected. So we get a well defined
semi-local holomorphic extension, again denoted by
$\psi\in\mathcal{O}(B_{q_1}\cup
\omega(\Sigma_{\gamma_{q_1}}))$. Geometrically speaking, this domain
$B_{q_1}\cup \omega(\Sigma_{\gamma_{q_1}})$ is a kind of curved
Hartogs domain. We claim that such a function $\psi$ extends
holomorphically to a neighborhood of the union of disc
$A_{q_1,\sigma}(\overline{\Delta})$ for $\vert \sigma\vert <\delta$.
Indeed, we first observe that for all $\vert q_1\vert < \varepsilon$
and $\vert \sigma\vert<\delta$, the boundaries
$A_{q_1,\sigma}(b\Delta)$ are contained in this domain $B_{q_1}\cup
\omega(\Sigma_{\gamma_{q_1}})$.  This is evident for the half
boundaries $A_{q_1,\sigma}(b^+\Delta)$ which are contained in
$\Sigma_{\gamma_{q_1}}$ by Lemma~6.15{\bf (2)}.  On the other hand,
the boundaries $A_{q_1,\sigma}(b^-\Delta)$ stay within a distance of
order say $3\delta/2$ with respect to the boundary
$A_{p_1}(b^-\Delta)$, by the very construction of the smooth family
$A_{q_1,\sigma}$, which proves the
claim. We remind the notion of {\it analytic
isotopy}\, of analytic discs ([Me2,~Definition~3.1]) which is useful
in applying the continuity principle.  For fixed $q_1$ and for varying
$\sigma$, all the discs $A_{q_1,\sigma}$ are analytically isotopic to
each other with their boundaries lying in
$B_{q_1}\cup\omega(\Sigma_{\gamma_{q_1}})$.  Moreover, for $\sigma=0$,
we obviously see that $A_{q_1,0}$ is analytically isotopic to a point
in the domain $B_{q_1}\cup\omega(\Sigma_{\gamma_{q_1}})$, just by the
trivial isotopy $(r,\zeta)\mapsto A_{q_1,0}(r\,\zeta)$ with values in
the neighborhood
$\omega(A_{q_1,0}(\overline{\Delta}))\subset B_{q_1}$.  By Lemma~3.2
in [Me2], it follows that $\psi$ restricted to a neighborhood of
$A_{q_1,\sigma}(b\Delta)$ extends holomorphically to a neighborhood of
$A_{q_1,\sigma}(\overline{\Delta})$ in $\C^n$, for all $\vert
\sigma\vert <\delta$. Furtermore, thanks to the fact that the map
$(\zeta,\sigma)\mapsto A_{q_1,\sigma}(\zeta)$ is an embedding, we get
a well defined holomorphic extension $\psi_{q_1}$ of $\psi$ to the
union $C_{q_1}:=\cup_{\vert \sigma\vert<\delta}\,
A_{q_1,\sigma}(\Delta\cap\Delta(1,\delta))$.  Of course, this
extension coincides with the old $\psi\in\mathcal{O}(D\cup \Omega)$ in
a neighborhood of the intersection of the half boundary 
$A_{q_1,0}(b^+\Delta)$ with $C_{q_1}$.
Since $C_{q_1}\cap D$ is connected and since $C_{q_1}$ contains
$\Delta_n(0,\eta)\cap \Sigma_{\gamma_{q_1}}^-$ by Lemma~6.17, after
sticking $\psi$ with $\psi_{q_1}$, we get the desired holomorphic
extension $\Psi\in\mathcal{O}(D\cup
[\Delta_n(0,\eta)\cap\Sigma_{\gamma_{q_1}}^-])$.  The proof of
Lemma~6.20 is complete.
\endproof

\section*{\S7. Holomorphic extension to a Levi-flat union of Segre varieties}

\subsection*{7.1.~Straightenings}
For each parameter $q_1'$, we consider the real
analytic arc $\gamma_{q_1'}'$ defined by~\thetag{6.2}.  To
this family of analytic arcs we can associate a family of
straightened coordinates as follows.

\def\thelemma{7.2}\begin{lemma} 
For varying $q_1'\in M'$ with $\vert q_1'\vert <
\varepsilon' < <\rho'$, there exists a real
analytically parameterized family of biholomorphic mappings
$\Phi_{q_1'}'$ of $\Delta_n(0,\rho'/2)$ sending $q_1'$ to the origin
and straightening $\gamma_{q_1'}'$ to the $u'$-axis, such that the
image $M_{q_1'}':=\Phi_{q_1'}(M')$ is a closed real analytic
hypersurface of $\Delta_n(0,\rho'/2)$ close to $M'$ in the real
analytic norm which is given by an equation of the form $\bar
w'=\Theta_{q_1'}'(\bar z',t')$, with $\Theta_{q_1'}'(\bar z',t')$
converging normally in the polydisc $\Delta_{2n-1}(0,\rho'/2)$ and
satisfying $\Theta_{q_1'}'(0,0,w')\equiv w'$ and $\Theta_{q_1'}'(\bar z',t')=
w'+{\rm O}(2)$.
\end{lemma}

\subsection*{7.3.~Different reflection functions}
Let us develope these defining equations in the form:
\def\theequation{7.4}\begin{equation}
\bar w'=\Theta_{q_1'}'(\bar z',t')=\sum_{\beta\in \N^{n-1}}
(\bar z')^\beta \, \Theta_{q_1',\beta}'(t').
\end{equation}
Here, $\Theta_{q_1',0}(0,w')\equiv w'$.  We denote by
$h_{q_1'}=(f_{q_1'},g_{q_1'})$ the mapping in these coordinate systems. To
every such system of coordinates, we associate {\it different
reflection functions}\, by setting:
\def\theequation{7.5}\begin{equation}
\mathcal{R}_{q_1',h_{q_1'}}'(t,\bar\nu')
:=\bar\mu'-\sum_{\beta\in\N^{n-1}}
\bar{\lambda'}^\beta \, \Theta_{q_1',\beta}'
(h_{q_1'}(t)).
\end{equation}

\subsection*{7.6.~Holomorphic extension to a Levi-flat hat} 
Recall from \S6.5 that the Levi-flat hypersurfaces
$\Sigma_{\gamma_{q_1}}$ are defined to be the union of the Segre
varieties $S_{\bar q}$ associated to points $q$ varying in
$\gamma_{q_1}$, intersected with the polydisc $\Delta_n(0,\rho/3)$.
Here, we establish our main crucial observation.

\def\thelemma{7.7}\begin{lemma}
If $q_1$ with $\vert q_1\vert < \varepsilon$ belongs to $M_1^-$, then
all the components $\Theta_{q_1',\beta}'(h(t))$ extend as CR functions of
class $\mathcal{C}^\infty$ over $\Sigma_{\gamma_{q_1}}\cap
\Delta_n(0,\rho/3)$.
\end{lemma}

\proof
Let $\overline{L}_1,\dots,\overline{L}_{n-1}$ be the commuting basis
of $T^{0,1}M$ given by $\overline{L}_j={\partial \over \partial \bar
z_j} +\Theta_{\bar z_j}(\bar z,t){\partial \over \partial \bar w}$,
for $1\leq j\leq n-1$. Clearly, the coefficients of these vectors
fields converge normally in the polydisc $\Delta_{2n-1}(0,\rho)$. By
the diffeomorphism assumption, we have $\hbox{det} (\overline{L}_j\,
\overline{f_{q_1',k}}(0))_{ 1\leq j,k\leq n-1}\neq 0$. At points
$(t,\bar t)$ with $t\in M\cap \Delta_n(0,\rho/3)$, we shall denote
this determinant by:
\def\theequation{7.8}\begin{equation}
{\rm det}\, (\overline{L}_j\, \overline{f_{q_1',k}}(\bar t))_{
1\leq j,k\leq n-1}:=
\underline{\mathcal{D}}(t,\bar t,\{\partial_{\bar t_l}
\overline{f_{q_1',k}}(\bar t) \}_{1\leq l\leq n, \, 1\leq k\leq n-1}).
\end{equation}
Here, by its very definition, the function $\underline{\mathcal{D}}$ is
holomorphic in its variables. Replacing $w$ by $\overline{\Theta}(z,\bar t)$
in $\underline{\mathcal{D}}$, we can write $\underline{\mathcal{D}}$
in the form $\mathcal{D}(z,\bar t,\{\partial_{\bar t_l} \overline{f_{q_1',k}}
(\bar t) \}_{1\leq l\leq n, \, 1\leq k\leq n-1})$, where $\mathcal{D}$ is
holomorphic in its variables.  Shrinking $\rho>0$ if necessary, we may
assume that for all fixed point $\bar t_q
\in M$ with $\vert \bar t_q\vert<\rho/3$, then:
\begin{itemize}
\item[{\bf (1)}]
{\it The polarization $\mathcal{D}(z,\bar t_q,\{\partial_{\bar t_l}
\overline{f_{q_1',k}}(\bar t_q) \}_{1\leq l\leq n, \, 1\leq k\leq n-1})$ is
convergent on the Segre variety $S_{\bar t_q}\cap
\Delta_n(0,\rho/3)=\{(z,w)\in \Delta_n(0,\rho/3):
w=\overline{\Theta}(z,\bar t_q)\}$, {\it i.e.}~it is convergent with
respect to $z\in\C^{n-1}$ for all $\vert z \vert < \rho/3$}.
\item[{\bf (2)}] 
{\it This expression $\mathcal{D}(z,\bar t_q,\{\partial_{\bar t_l}
\overline{f_{q_1',k}}(\bar 
t_q) \}_{1\leq l\leq n, \, 1\leq k\leq n-1})$ does {\rm not}
vanish at any point of the Segre variety $S_{\bar t_q}\cap
\Delta_n(0,\rho/3)$, {\it i.e.} it does not vanish for all $\vert z
\vert < \rho/3$}.
\end{itemize} 
Let us choose $q_1'$ satisfying $\gamma_{q_1'}'\subset {M_1'}^-$, with
$\vert q_1'\vert <\varepsilon'$.  We pick the corresponding parameter
$q_1:=h^{-1}(q_1')$ with $\vert q_1 \vert < \varepsilon$. By the choice of
$\Phi_{q_1'}'$, we then have $f_{q_1'}(\gamma_{q_1}(s))=0$ for all
$s\in \R$ with $\vert s\vert \leq \rho/2$. {\it This property will be
really crucial}. As the mapping $h_{q_1'}$ is of class
$\mathcal{C}^\infty$ over $M$, we can apply the tangential
Cauchy-Riemann derivations $\overline{L}_1^{\beta_1}\cdots
\overline{L}_{n-1}^{\beta_{n-1}}$, $\beta\in \N^{n-1}$ of order
$\vert\beta\vert$ infinitely many
times to the identity: \def\theequation{7.9}\begin{equation}
\overline{g_{q_1'}(t)}=\Theta_{q_1'}' (\overline{f_{q_1'}(t)},
h_{q_1'}(t)).
\end{equation}
which holds for $t\in M\cap \Delta_n(0,\rho)$.  To begin with, we
first apply the CR derivations $\overline{L}_j$ to this
identity~(7.9).  This yields
\def\theequation{7.10}\begin{equation}
\overline{L}_j\,\overline{g_{q_1'}}(\bar t)=
\sum_{k=1}^{n-1}\, 
\overline{L}_j\, \overline{f_{q_1',k}}(\bar t) \
{\partial \Theta_{q_1'}'\over\partial \bar z_k'}\,
(\overline{f_{q_1'}}(\bar t), h_{q_1'}(t)).
\end{equation}
Applying Cramer's rule as in the proofs of Lemmas~3.22 and~4.7, 
we see that there exist holomorphic functions $\mathcal{T}_k$ in their
arguments such that
\def\theequation{7.11}\begin{equation}
{\partial\Theta_{q_1'}'\over\partial \bar z_k'}
(\overline{f_{q_1'}(t)},
h_{q_1'}(t))=
{\mathcal{T}_k(z,\bar t,\{\partial_{\bar t_l}
\overline{h_{q_1',j}}(\bar t)\}_{
1\leq l,j\leq n})
\over
\mathcal{D}(z,\bar t,\{\partial_{\bar t_l}\,
\overline{f_{q_1',k}}(\bar t)\}_{1\leq l\leq n,\, 1\leq k\leq n-1})}.
\end{equation}
By CR differentiating further the identities~(7.11), using Cramer's rule 
at each step and making inductive arguments, it follows that for every 
multi-index $\beta\in\N_*^{n-1}$, there exist holomorphic
functions $\mathcal{T}_\beta$ in their variables such that
\def\theequation{7.12}\begin{equation}
\left\{
\aligned
{}
&
\Theta_{q_1',\beta}'(h_{q_1'}(t))+\sum_{\gamma\in\N_*^{n-1}}
 \ (\overline{f_{q_1'}(t)})^\gamma
\ \Theta_{q_1',\beta+\gamma}'(h_{q_1'}(t))
\ {(\beta+\gamma)!\over \beta!\ \gamma!}=
\\
& 
\ \ \ \ \ \ \ \ \ \ \ \ \ \ \ \ \ \ \ \ \
={\mathcal{T}_\beta(z,\bar t,\{\partial_{\bar t}^\gamma 
\overline{h_{q_1',j}}(\bar t)\}_{1\leq j\leq n, \, \vert
\gamma\vert \leq \vert \beta \vert})\over 
[\mathcal{D}(z,\bar t,\{\partial_{\bar t_l}
\overline{f_{q_1',k}}(\bar t) \}_{1\leq l\leq n, 
\, 1\leq k\leq n-1})]^{2\vert \beta\vert -1}}.
\endaligned\right.
\end{equation}
Precisely, the terms $\mathcal{T}_\beta$ are holomorphic with respect to
$(z,\bar t)$ and relatively polynomial with respect to the jets
$\{\partial_{\bar t}^\gamma \overline{h_{q_1',j}}(\bar t) \}_{1\leq
j\leq n,\,\vert\gamma\vert\leq\vert\beta\vert}$.  Also, the variable
$t$ runs in $M$ in a neighborhood of $\gamma_{q_1}(s)$. Now, we remind
that by Lemma~5.10, all the functions
$t\mapsto\Theta_{q_1',\beta}'(h_{q_1', \beta}(t))$ are already
holomorphically extendable to a neighborhood of $\gamma_{q_1}$ in $\C^n$,
since $\gamma_{q_1}\subset M_1^-$. Let us denote by
$\theta_{q_1',\beta}'(t)$ these holomorphic extensions. We shall
first prove Lemma~7.7 in the simpler case where $M'$ is holomorphically 
nondegenerate, in which case the mapping $h$ in fact extends
holomorphically to a neighborhood of $M_1^-$ in $\C^n$. In this case,
for every point $q\in \gamma_{q_1}$ of the form
$q=\gamma_{q_1}(s)$, 
the terms in the right hand side of~(7.12) extend holomorphically 
to a neighborhood 
$(q,\bar q)$ of 
the complexification $\mathcal{M}$ of $M$, which is
the complex hypersurface in $\C^n\times\C^n$
given by the defining equation $w=\overline{\Theta}(z,\tau)$. 
So, for $(t,\tau)$ close to $(q,\bar q)$, we can complexify~(7.12), 
replacing $\bar t$ by $\tau$ and 
$t$ by $(z,\overline{\Theta}(z,\tau))$, which yields an identity between
holomorphic functions:
\def\theequation{7.13}\begin{equation}
\left\{
\aligned
{}
&
\Theta_{q_1',\beta}'(h_{q_1'}(z,\overline{\Theta}(z,\tau)))+
\sum_{\gamma\in\N_*^{n-1}}
 \ (\overline{f_{q_1'}}(\tau))^\gamma
\ \Theta_{q_1',\beta+\gamma}'(h_{q_1'}(z,\overline{\Theta}(z,\tau)))
\ {(\beta+\gamma)!\over \beta!\ \gamma!}\equiv
\\
& 
\ \ \ \ \ \ \ \ \ \ \ \ \ \ \ \ \ \ \ \ \
\equiv {\mathcal{T}_\beta(z,\tau,\{\partial_{\tau}^\gamma 
\overline{h_{q_1',j}}(\tau)\}_{1\leq j\leq n, \, \vert
\gamma\vert \leq \vert \beta \vert})\over 
[\mathcal{D}(z,\tau,\{\partial_{\tau_l}
\overline{f_{q_1',k}}(\tau) \}_{1\leq l\leq n, 
\, 1\leq k\leq n-1})]^{2\vert \beta\vert -1}}.
\endaligned\right.
\end{equation}
Next, we put $\tau:=\bar t_q= \overline{\gamma_{q_1}(s)}$, whence $t$
belongs to the Segre variety $S_{\bar q}$, namely
$t=(z,\overline{\Theta}(z,\bar t_q))$, where the variable $z$ is
free. From the important fact that $\overline{f_{q_1'}}(\bar t_q)=0$,
because $h(q)$ belongs to $\gamma_{q_1'}'$, it follows that 
{\it the queue sum $\sum_{\gamma\in\N_*^{n-1}}$ in~(7.13) disappears}.
Consequently, we get the following identity on $S_{\bar q}$ for 
$z$ close to $z_q$:
\def\theequation{7.14}\begin{equation}
\Theta_{q_1',\beta}'(h_{q_1'}(z,\overline{\Theta}(z,\bar t_q)))
\equiv {\mathcal{T}_\beta(z,\bar t_q,\{\partial_{\bar t}^\gamma 
\overline{h_{q_1',j}}(\bar t_q)\}_{1\leq j\leq n, \, \vert
\gamma\vert \leq \vert \beta \vert})\over 
[\mathcal{D}(z,\bar t,\{\partial_{\bar t_l}
\overline{f_{q_1',k}}(\bar t_q) \}_{1\leq l\leq n, 
\, 1\leq k\leq n-1})]^{2\vert \beta\vert -1}}.
\end{equation}
The crucial observation now is that the right hand side of~(7.14)
converges over a much longer part of the Segre variety $S_{\bar q}$.
Indeed, by {\bf (1)} after~(7.8), it converges for $\vert z\vert <
\rho/3$.  Furthermore, the right hand side of~(7.14) varies in a
$\mathcal{C}^\infty$ way when $\bar t_q$ varies on
$\gamma_{q_1}$. This proves Lemma~7.7 in the case where $h$ extends
holomorphically to a neighborhood of $M_1^-$ in $\C^n$, which holds
true for instance when $M'$ is holomorphically nondegenerate.

In the general case, it is no longer true that $h$ extends
holomorphically to a neighborhood of $M_1^-$ in $\C^n$, so different
arguments are required. Let $q\in\gamma_{q_1}$ be arbitrary. By
assumption, the components $\Theta_{q_1',\beta}'(h_{q_1'}(t))$ extend
holomorphically to a neighborhood of $q$ in $\C^n$ as holomorphic
functions $\theta_{q_1',\beta}'(t)$ defined, say in the
polydisc $\{\vert t-t_q\vert <
\sigma_q\}$, for small $\sigma_q>0$.
By expanding $h_{q_1'}$ in formal power series at $q$, we get
a series $H_{q_1'}(t_q+(t-t_q))\in\C\dl t-t_q\dr^n$. Also, 
we may expand $\theta_{q_1',\beta}'(t_q+(t-t_q))\in\C\{t-t_q\}$.
Then we have the following formal power series
identities
\def\theequation{7.15}\begin{equation}
\Theta_{q_1',\beta}(H_{q_1'}(t_q+(t-t_q)))\equiv
\theta_{q_1',\beta}(t_q+(t-t_q))
\end{equation} 
in $\C\dl t-t_q\dr$ 
for all $\beta$. Since the Taylor series of $(h_{q_1'},
\overline{h_{q_1'}})$ at $(t_q,\bar t_q)$ induces a formal CR
mapping between the complexifications $\mathcal{M}$ centered at
$(q,\bar q)$ and the complexification $\mathcal{M}'$ centered at
$(q',\bar q')$, it follows that we can write the following formal
power series identities valuable in $\C\dl t-t_q,\tau-\bar t_q\dr$ 
for $(t_q+(t-t_q),\bar t_q+(\tau-\bar
t_q))$ in $\mathcal{M}$:
\def\theequation{7.16}\begin{equation}
\left\{
\aligned
{}
&
\Theta_{q_1',\beta}'(H_{q_1'}(t_q+(t-t_q)))+\\
&
+\sum_{\gamma\in\N_*^{n-1}}\,
\Theta_{q_1',\beta+\gamma}'(H_{q_1'}(t_q+(t-t_q))\,
\overline{F_{q_1'}}(\bar t_q+(\tau-\bar t_q))^\gamma\equiv\\
&
\equiv
{\mathcal{T}_\beta(z_q+(z-z_q),\bar t_q+\tau-\bar t_q,
\{\partial_{\tau}^\gamma 
\overline{H_{q_1',j}}(\bar t_q+(\tau-\bar t_q))\}_{1\leq j\leq n, \, \vert
\gamma\vert \leq \vert \beta \vert})\over 
[\mathcal{D}(z_q+(z-z_q),\bar t_q+(\tau-\bar t_q),\{\partial_{\tau_l}
\overline{F_{q_1',k}}(\bar t_q+(\tau-\bar t_q)) \}_{1\leq l\leq n, 
\, 1\leq k\leq n-1})]^{2\vert \beta\vert -1}}.
\endaligned\right.
\end{equation}
Putting $\tau:=\bar t_q$ in~(7.16), taking~(7.15) into account, and
using the important fact that $\overline{F_{q_1'}}(\bar t_q)=0$, we
get the formal power series identities between two holomorphic functions
which are valuable for $\vert z-z_q\vert<\sigma_q$ in $\C\{z-z_q\}$
and for all $\beta$:
\def\theequation{7.17}\begin{equation}
\left\{
\aligned
{}
&
\theta_{q_1',\beta}'(z_q+(z-z_q),
\overline{\Theta}(z_q+(z-z_q),\bar t_q))\equiv\\
&
\equiv
{\mathcal{T}_\beta(z_q+(z-z_q),\bar t_q,\{\partial_{\bar t}^\gamma 
\overline{h_{q_1',j}}(\bar t_q)\}_{1\leq j\leq n, \, \vert
\gamma\vert \leq \vert \beta \vert})\over 
[\mathcal{D}(z_q+(z-z_q),\bar t_q,\{\partial_{\bar t_l}
\overline{f_{q_1',k}}(\bar t_q) \}_{1\leq l\leq n, 
\, 1\leq k\leq n-1})]^{2\vert \beta\vert -1}}.
\endaligned\right.
\end{equation}
Consequently, we get on $S_{\bar q}$ the following
identities between holomorphic functions of $z$ valuable for $\vert
z-z_q\vert<\sigma_q$ and for all $\beta$: 
\def\theequation{7.18}\begin{equation}
\theta_{q_1',\beta}'(z,\overline{\Theta}(z,\bar t_q))\equiv
{\mathcal{T}_\beta(z,\bar t_q,\{\partial_{\bar t}^\gamma 
\overline{h_{q_1',j}}(\bar t_q)\}_{1\leq j\leq n, \, \vert
\gamma\vert \leq \vert \beta \vert})\over 
[\mathcal{D}(z,\bar t_q,\{\partial_{\bar t_l}
\overline{f_{q_1',k}}(\bar t_q) \}_{1\leq l\leq n, 
\, 1\leq k\leq n-1})]^{2\vert \beta\vert -1}}.
\end{equation}
As in the holomorphically nondegenerate case, we see that the right
hand side of~(7.18) converges for $\vert z\vert <\rho/3$, so the
holomorphic functions $\theta_{q_1',\beta}'(z,\overline{\Theta}(z,\bar
t_q))$ converge in a long piece of the Segre variety $S_{\bar q}$.
The $\mathcal{C}^\infty$-smoothness of the right hand side extension
over $\Sigma_{\gamma_{q_1}}\cap \Delta_n(0,\rho/3)$ yields a CR
extension to $\Sigma_{\gamma_{q_1}}$ which is of class
$\mathcal{C}^\infty$. This completes the proof of Lemma~7.7.
\endproof

\def\thelemma{7.19}\begin{lemma}
If $q_1$ with $\vert q_1\vert <\varepsilon$ belongs to $M_1^-$, then
all the components $\Theta_{q_1',\beta}'(h_{q_1'}(t))$ of the reflection
function $\mathcal{R}_{{q_1'},h_{q_1'}}$ extend
as holomorphic functions to a neighborhood
$\omega(\Sigma_{\gamma_{q_1}})$ of $\Sigma_{\gamma_{q_1}}$ in $\C^n$.
\end{lemma}

\proof
By the hypotheses of Theorem~5.5 and by Lemma~5.10, we remind the
reader that the components $\Theta_{q_1',\beta}'(h_{q_1'}(t))$ already
extend holomorphically to a neighborhood
$\omega(\gamma_{q_1})\subset\Omega$ of $\gamma_{q_1}\subset M_1^-$ in
$\C^n$ as the holomorphic functions $\theta_{q_1',\beta}'(t)$. Thanks to
Lemma~7.7, the statement follows by an application of the following
known propagation result:
\endproof

\def\thelemma{7.20}\begin{lemma}
Let $\Sigma$ be a $\mathcal{C}^\infty$-smooth Levi-flat hypersurface
in $\C^n$ $(n\geq 2)$ foliated by complex hypersurfaces
$\mathcal{F}_\Sigma$. If a continuous CR function $\psi$ defined on
$\Sigma$ extends holomorphically to a neighborhood $\mathcal{U}_p$ of
a point $p$ belonging to a leaf $\mathcal{F}_\Sigma$ of $\Sigma$, then
$\psi$ extends holomorphically to a neighborhood
$\omega(\mathcal{F}_\Sigma)$ of $\mathcal{ F}_\Sigma$ in $\C^n$. The
size of this neighborhood $\omega(\mathcal{ F}_\Sigma)$ depends on the
size of $\mathcal{U}_p$ and is stable under sufficiently small $($even
non-Levi-flat$)$ perturbations of $\Sigma$.
\end{lemma}

\proof
The first part of this statement was first proved by Hanges and Treves
([HaTr]) using microlocal concepts, the {\sc fbi} transform and
controlled deformations of manifolds. Interesting generalizations were
given by Sj\"ostrand and by Tr\'epreau ([Tr2]) in arbitrary
codimension. Another proof using deformations of analytic discs has
been provided by Tumanov ([Tu2]).  Both proofs are constructive and
the second statement about the size of the neighborhoods to which
extension holds follows after a careful inspection of the techniques
therein. Since it is superfluous to repeat the arguments word by word,
we do not enter the details.
\endproof

\section*{\S8. Relative position of the neighbouring Segre varieties}

\subsection*{8.1.~Intersection of Segre varieties}
We are now in position to complete the proof of Theorem~5.5, hence to
achieve the proof of Theorem~1.9.  It remains to show that the
functions $\Theta_{q_1',\beta}'$ extend holomorphically at $p_1$, for
$\gamma_{q_1}$ chosen conveniently. For this choice, we are led to the
following dichotomy: either $S_{\bar p_1}\cap M_1^-=\emptyset$ in a
sufficiently small neighborhood of $p_1$ or there exists a sequence
$(q_k)_{k\in\N}$ of points of $S_{\bar p_1}\cap M_1^-$ tending towards
$p_1$. In the first case, we shall distinguish two sub-cases. Either
$S_{\bar p_1}$ lies below $M_1^-$ or it lies above $M_1^-$. Let us
write this more precisely. We may choose a $\mathcal{C}^\infty$-smooth
hypersurface $H_1$ transverse to $M$ at $p_1$ with $H_1$ satisfying
$H_1\cap M=M_1$ and $H_1^-\cap M=M_1^-$ ({\it see} {\sc
Figure}~7). Thus $H_1$ together with $M$ divides $\C^n$ near $0$ in
four connected parts. More precisely, we say that either $S_{\bar
p_1}\cap H_1^-$ is contained in the lower left quadrant $H_1^-\cap
M^-=H_1^-\cap D$ or it is contained in the upper left quadrant
$H_1^-\cap M^+$. To summarize, we have distinguished three possible
cases:
\begin{itemize}
\item[{\bf Case~I.}] \ 
The half Segre variety $S_{\bar p_1}\cap H_1^-$ cuts $M_1^-$ along an
infinite sequence of points $(q_k)_{k\in\N}$ tending towards $p_1$.
\smallskip
\item[{\bf Case~II.}] \ 
The half Segre variety $S_{\bar p_1}\cap H_1^-$ does not intersect
$M_1^-$ in a neighborhood of $p_1$ and it passes under $M_1^-$, 
namely inside $D$.
\smallskip
\item[{\bf Case~III.}] \ 
The half Segre variety $S_{\bar p_1}\cap H_1^-$ does not intersect
$M_1^-$ in a neighborhood of $p_1$ and it passes over $M_1^-$, namely 
over $D\cup M_1^-$.
\end{itemize}  
In the first two cases, for every point
$q_1$ close enough to $p_1$, the Segre variety
$S_{\bar q_1}$ will intersect $D\cup \Omega$ and the
neighborhoods $\omega(\Sigma_{\gamma_{q_1}})$ constructed in
Lemma~7.20 will always contain the point $p_1$ (we give more
arguments below). The third case could
be {\it a priori}\, the most delicate one. But we can already delineate the
following crucial geometric property, which says that Lemma~6.20 will
apply.

\def\thelemma{8.2}\begin{lemma}
If $S_{\bar p_1} \cap H_1^-$ is contained in $M^+$, then $p_1$ lies in the
lower side $\Sigma_{\gamma_{q_1}}^-$ for every arc
$\gamma_{q_1}\subset M_1^-$ of the family~\thetag{6.2}.
\end{lemma}

\proof In normal coordinates $t$ vanishing at $p_1$,
the real equation of $M$ is given by $v=\varphi(z,\bar z,u)$, where
$\varphi$ is a certain converging real power series satisfying
$\varphi(0)=0$, $d\varphi(0)=0$ and $\varphi(z,0,u)\equiv 0$. We can
assume that $dh(0)={\rm Id}$. We can assume that the ``minus'' side
$D\equiv M^-$ of automatic extension of CR functions is given by
$\{v<\varphi(z,\bar z,u)\}$. Replacing $u$ by $(w+\bar w)/2$ and $v$ by
$(w-\bar w)/2i$, and solving with respect to $w$, we get for $M$ an
equation as above, say $w=\bar w+i \, \overline{\Xi}(z,\bar t)$, with
$\overline{\Xi}(0,\bar t)\equiv 0$. We have $\overline{\Theta} (z,\bar
t)\equiv \bar w+i \, \overline{\Xi}(z,\bar t)$ in our previous
notation.  We claim that every such arc $\gamma_{q_1}\subset
M_1^-$ contains a point $p\in M_1^-$ whose coordinates are of the form
$(z_p,0+i\varphi(z_p,\bar z_p,0))$.  Indeed, by construction, the arcs
$\gamma_{q_1}$ are all elongated along the $u$-coordinate axis, since
it is so for $\gamma_{q_1'}'$ and since $dh(0)=\text{\rm Id}$. In
normal coordinates, the Segre variety $S_{\bar p_1}$
passing through the origin $p_1$
has the simple equation $\{w=0\}$. By assumption, the point
$(z_p,0)\in S_{\bar p_1}$ lies over $M$ in $M^+$, so we have
$\varphi(z_p,\bar z_p,0)<0$. Then the Segre variety $S_{\bar p}$
(which is a leaf of $\Sigma_{\gamma_{q_1}}$), has the equation
$w=-i\varphi(z_p,\bar z_p,0)+i \, \overline{\Xi}(z,\bar
z_p,-i\varphi(z_p,\bar z_p,0))$.  Therefore, the intersection point
$\{z=0\}\cap S_{\bar p}\subset \Sigma_{\gamma_{q_1}}$ has coordinates
equal to $(0,-i\varphi (z_p,\bar z_p,0))$.  This point clearly lies
above the origin $p_1$, so $p_1$ lies in the lower side
$\Sigma_{\gamma_{q_1}}^-$, which completes the proof of Lemma~8.2.
\endproof

\subsection*{8.3.~Extension across $(M,0)$ of the components
$\Theta_{q_1',\beta}'$} We are now prepared to complete the proof of
Theorems~5.5 and~1.9.  We first choose $\delta, \, \eta, \,
\varepsilon$ and various points $\vert q_1 \vert<\varepsilon$ as in
Lemma~6.20 and we consider the two associated arcs $\gamma_{q_1}$ and
$\gamma_{q_1'}'$, the associated mapping $h_{q_1'}$ and the associated
reflection function $\mathcal{R}_{q_1', h_{q_1'}}'$. By Lemma~6.20,
for each such choice of $q_1$, then all the components
$\Theta_{q_1',\beta}'$ extend holomorphically to $D\cup
[\Sigma_{\gamma_{q_1}}^-\cap \Delta_n(0,\eta)]$. Our goal is to show
that for suitably chosen $\gamma_{q_1}$ in Cases I, II and III, then
the components $\Theta_{q_1',\beta}'$ extend holomorphically to a
neighborhood of $p_1$.  Afterwards, thanks to 
Artin's approximation theorem, the Cauchy estimates are
automatic, as explained in Lemma~3.16.

\subsection*{8.4.~Case I}
In Case~I, we choose one of the points 
$q_k\in M_1^-\cap S_{\bar p_1}$ which is
arbitrarily close to $p_1$ and we denote it simply by $q_1$. We can
assume that $\vert q_1\vert < \varepsilon$. Next, we consider the
associated arc $\gamma_{q_1}$. By an application of Lemma~7.19, all
the components $\Theta_{q_1',\beta}'(h_{q_1'}(t))$ of the reflection
function $\mathcal{R}_{q_1',h_{q_1'}}$ extend holomorphically to a
neighborhood $\omega(\Sigma_{\gamma_{q_1}})$ of
$\Sigma_{\gamma_{q_1}}$ in $\C^n$. Of course, this neighborhood
contains the point $p_1\in S_{\bar q_1}
\subset \Sigma_{\gamma_{q_1}}$. However, because of
possible pluridromy, the extension at $p_1$ might well differ from the
extension in the one-sided neighborhood $D$ near $p_1$.  Fortunately,
thanks to Lemma~6.20, all these holomorphic functions extend
holomorphically in a unique way to $D\cup [\Sigma_{\gamma_{q_1}}^-\cap
\Delta_n(0,\eta)]$. The neighborhood $\omega(\Sigma_{\gamma_{q_1}})$
being constructed as a certain union of polydiscs of small radius, it is
geometrically smooth, so its intersection with 
$D\cup [\Sigma_{\gamma_{q_1}}^-\cap \Delta_n(0,\eta)]$ is
connected. In sum, we get unique holomorphic extensions of
the functions $\Theta_{q_1',\beta}'(h_{q_1'}(t))$ to the domain 
\def\theequation{8.5}\begin{equation}
\omega(\Sigma_{\gamma_{q_1}})\ \cup \
D \ \cup \ [\Sigma_{\gamma_{q_1}}^-\cap \Delta_n(0,\eta)],
\end{equation}
which yields the desired holomorphic extensions at $p_1$.
Case~I is achieved.

\subsection*{8.6.~Case II}
Case~II is treated almost the same way as Case~I. Since $S_{\bar
p_1}\cap H_1^-$ is contained in $D$, we can choose a fixed point
$\tilde{q}$ of $S_{\bar p_1}$ which belongs to
$D$. So there exists a radius $\tilde{\rho}>0$ such
that the polydisc $\Delta_n(\tilde{q},\tilde{\rho})$ 
is contained in $D$. For $\vert q_1\vert < \varepsilon$
sufficiently close to $p_1$, there exists a point $\tilde{q}_1\in
S_{\bar q_1}$ sufficiently close to $\tilde{q}$ such that the polydisc
$\Delta_n(\tilde{q}_1,\tilde{\rho}/2)$ is again contained in $D$.
Thanks to Lemma~7.20, if $q_1$ is sufficiently close to $p_1$, the
neighborhood $\omega(\Sigma_{\gamma_{q_1}})$ constructed by
deformations of analytic discs as in [Tu2] will contain the point
$p_1$, since its size along $S_{\bar q_1}$ depends only on the fixed
size of the polydisc $\Delta_n(\tilde{q}_1,\tilde{\rho}/2)$ which is
of radius at least $\tilde{\rho}/2$ uniformly. Finally, as in Case~I, the
monodromy of the extension follows by an application of 
Lemma~6.20.

\bigskip
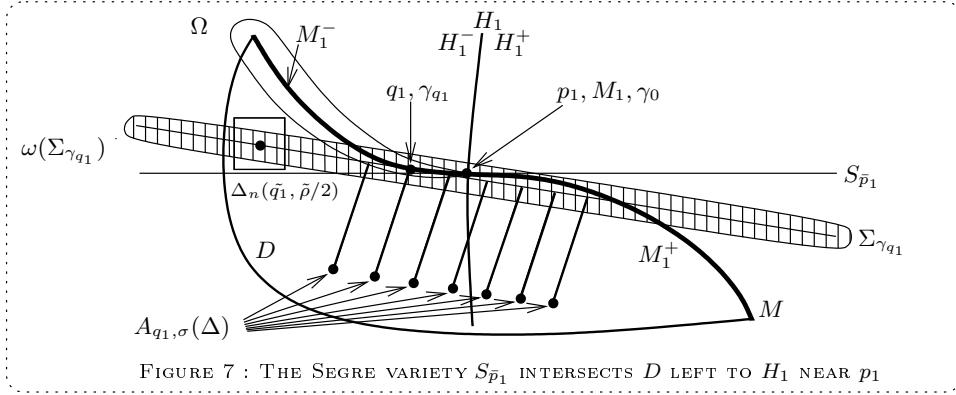
\begin{center}
\input figure8.pstex_t
\end{center}
\bigskip

\subsection*{8.7.~Case III}
For Case~III, thanks to Lemma~8.2, we know already that $p_1$ belongs
to the lower side $\Sigma_{\gamma_{q_1}}^-$. Thus Case~III follows
immediately from the application of Lemma~6.20 summarized in \S8.3
above.  Case~III is achieved. The proofs of Theorems~5.5, 1.9 and~1.2 
are complete.
\endproof

\section*{\S9.~Analyticity of some degenerate 
$\mathcal{C}^\infty$-smooth CR mappings}

\subsection*{9.1.~Presentation of the results}
Theorems~1.9 and~1.14 are concerned with $\mathcal{C}^\infty$-smooth
CR {\it diffeomorphisms}. It is desirable to remove the diffeomorphism
assumption. Taking inspiration from the very deep article of Pinchuk
[P4], we have been successful in establishing the following
statement. We refer the reader to the work of Diederich-Forn{\ae}ss
[DF1] and to the book of D'Angelo [D'A] for fundamentals about complex
curves contained in real analytic hypersurfaces.

\def\thetheorem{9.2}\begin{theorem}
Let $h: M\to M'$ be a $\mathcal{C}^\infty$-smooth CR mapping between
two connected real analytic hypersurfaces in $\C^n$ $(n\geq 2)$.  If
$M$ and $M'$ do not contain any complex curve, then $h$ is
real analytic at {\rm every} point of $M$.
\end{theorem}

At first, we need to recall some known facts about the
local CR geometry of real analytic hypersurfaces.

\begin{itemize}
\item[{\bf (1)}]
If $M$ does not contain complex curves, it is essentially finite.
This is obvious, because the coincidence loci of Segre varieties are
complex analytic subsets which are {\it contained}\, in $M$ ({\it cf.}
[DP1,2]).
\smallskip
\item[{\bf (2)}]
If $M$ is essentially finite at every point, it is locally minimal at
every point, so it consists of a single CR orbit, namely it is
globally minimal. As we have seen in \S3.6 above, CR functions on $M$
(and in particular the components of $h$) extend holomorphically to a
global one-sided neighborhood $D$ of $M$ in $\C^n$.
\smallskip
\item[{\bf (3)}]
If $M$ does not contain complex curves, then $M$ is Levi nondegenerate
at each point of the complement of some proper closed real analytic
subset of $M$. On the contrary, the everywhere Levi degenerate CR manifolds
are locally regularly foliated by complex leaves of dimension equal to
the dimension of the kernel the Levi form, at points where this kernel
is of maximal hence
locally constant dimension. This may happen in the class of essentially 
finite hypersurfaces.
\smallskip
\item[{\bf (4)}]
If $M$ does not contain complex curves, then either $h$ is constant
or it is of real generic rank $(2n-1)$ over an open dense
subset of $M$ and
its holomorphic extension is of complex generic rank $n$ over $D$.
This is easily established by looking at a point where $h$ is 
of maximal, hence locally constant, rank.
\smallskip
\item[{\bf (5)}]
In Theorem~9.2, there exists at least an everywhere dense open subset
$U_M$ of $M$ such that $h$ is real analytic at every point of $U_M$.
\end{itemize}

Based on these observations, Theorem~9.2 will be implied by 
the following more general statement to which the remainder
of \S9 is devoted.

\def\thetheorem{9.3}\begin{theorem}
Let $h: M\to M'$ be a $\mathcal{C}^\infty$-smooth CR mapping between
two connected real analytic hypersurfaces in $\C^n$ $(n\geq 2)$.  If
$M$ and $M'$ are essentially finite at {\rm every} point and if the maximal
generic real rank of $h$ over $M$ is equal to $(2n-1)$, then $h$ is
real analytic at {\rm every} point of $M$.
\end{theorem}

In [BJT], [BR1], [BR2], an apparently similar result is proved. In
these articles, it is always assumed at least that the formal Taylor
series of $h$ at every point of $M$ has Jacobian determinant not
identically zero. It follows that all the results proved in these
papers are superseded by the unification provided in the recent
articles [CPS1,2] and [Da2] expressed in terms of the characteristic
variety~(1.4).  However, the difficult problem would be to treat the
points of $M$ where {\it nothing}\, is {\it a priori}\, known about
the behavior of $h$, for instance points where all the $h_j$ could
vanish to infinite order hence have an identically zero formal Taylor
series. In this case, of course, the characteristic variety is
positive-dimensional. Unless $M$ is strongly pseudoconvex or there
exist local peak functions, it seems impossible to show {\it ab
initio}\, that $h$ is not flat at every point of $M$. Thus the
strategy of working {\it only at one fixed ``center point'' of $M$}\,
might well necessarily fail ({\it cf.}~[BJT], [BR1,2,4],
[BER1,2,3]). On the contrary, a strategy of propagation from nearby
points as developed in [P3,4], [DFY], [DP1,2], [Sha], [V], [PV] (and also in
the previous paragraphs) is really adequate.  Philosophically
speaking, there is no real surprise here, because the propagation
along Segre varieties is a natural generalization of the
weierstrassian conception of analytic continuation.

\subsection*{9.4.~Dense holomorphic extension} Let $D$ be a global 
one-sided neighborhood of $M$ in $\C^n$ to which CR functions extend
holomorphically. It follows from the assumptions of Theorem~9.3 that
the generic complex rank of $h$ in $D$ equals $n$. Recall that the two
everywhere essentially finite hypersurfaces $M$ and $M'$ are of course
holomorphically nondegenerate, namely $\chi_M=n$ and
$\chi_{M'}'=n$. At first, we prove the following lemma. Recall that
the intrinsic exceptional locus $E_M$ defined in \S3.47 is a proper
real analytic subset of $M$. Let $U_M$ denote the open subset
consisting of points $p\in M\backslash E_M$ at which the real rank of
$h$ equals $(2n-1)$.

\def\thelemma{9.5}\begin{lemma}
The open subset $U_M$ is dense in $M$. 
\end{lemma}

\proof
Indeed, suppose on the contrary that $M\backslash U_M$ contains an
open set $V$. Then the rank of $h$ is strictly less than $(2n-1)$ over
$V$. By the principle of analytic continuation and by the boundary
uniqueness theorem, it follows that $h$ is of generic complex rank
strictly less than $n$ in the domain $D$, contradiction.  
\endproof

\def\thelemma{9.6}\begin{lemma}
The mapping $h$ extends holomorphically to a neighborhood of 
every point $p\in U_M$.
\end{lemma}

\proof
Indeed, at such a point $p\in U_M$, $h$ is a local CR diffeomorphism
of class $\mathcal{C}^\infty$. By Lemma~4.3, the image $p':=h(p)$ of
$p$ belongs to $M'\backslash E_{M'}'$. Then Lemma~4.11 applies directly
(with $\chi_{M'}'=n$ 
of course) to show that $h$  extends holomorphically at $p$.
\endproof

\subsection*{9.7.~Holomorphic and formal mappings of 
essentially finite hypersurfaces} Let $h: M\to M'$ be as in the a
hypotheses of
Theorem~9.3. Let $p\in M$ and let $p':=h(p)$. Let $t$ be coordinates
vanishing at $p$ and let as usual a complex equation for the extrinsic
complexification $\mathcal{M}$ of $M$ be of the form $
w=\overline{\Theta}(z,\zeta,\xi)$, where $t=(z,w)\in \C^{n-1}\times\C$ and
$\tau=(\zeta,\xi)\in\C^{n-1}\times \C$.  Similarly, let
$w'=\overline{\Theta}'( z',\zeta',\xi')$ be an equation of
$\mathcal{M}'$. As in the proof of Lemma~4.3, the $\mathcal{C}^\infty$-smooth
CR mapping $h$ induces a formal CR mapping $(H(t),\overline{H}(\tau))$
between $(\mathcal{M},(p,\bar p))$ and $(\mathcal{M}',(p',
\overline{p'}))$. Precisely, this means that the Taylor series
$H_j(t)=\sum_{\gamma\in\N^n}\, H_{j,\gamma}\, t^\gamma$ of $h_j$ at
the origin and there conjugates $\overline{H}(\tau)$ satisfy a formal
power series identity of the form
\def\theequation{9.8}\begin{equation}
G(t)-\overline{\Theta}'(F(t),\overline{F}(\tau),\overline{G}(\tau))\equiv
A(t,\tau)\,\left[w-\overline{\Theta}(z,\zeta,\xi)\right],
\end{equation} 
where we denote $H=(F_1,\dots,F_{n-1},G)$ and where $A(t,\tau)$ is a
formal power series.  Without loss of
generality, we can assume that the coordinates $(z,w)$ and $(z',w')$
are normal, namely the defining functions satisfy
$\overline{\Theta}(0,\zeta,\xi)\equiv \overline{\Theta}(
z,0,\xi)\equiv \xi$ and {\it idem}\, for $\overline{\Theta}'$.  Such
coordinates are not unique, but they specify a certain component
$H_n=G$ of the formal CR mapping $H$ which is called a {\it
transversal component}. In [BR2], two facts about formal CR
mappings between small local pieces of real analytic hypersurfaces
(and even between {\it formal}\, hypersurfaces) are
established. Recall that $M$ and $M'$ are assumed to be essentially
finite at the origin and that the coordinates are normal.
\smallskip
\begin{itemize}
\item[{\bf 1.}]
If the transversal power series $G$ does not vanish identically, then
$H$ is of {\it finite multiplicity}, namely ({\it cf.}~[BR88]), the
ideal generated by the power series
$F_1(z,0),\dots,F_{n-1}(z,0)$ is of finite codimension in 
$\C\dl z\dr$. We denote this codimension by 
${\rm Mult}\,(H,0)$. It is independent of 
normal coordinates.
\smallskip
\item[{\bf 2.}]
If $H$ is of finite multiplicity, then a formal Hopf Lemma holds at
the origin, which tells us that the induced formal mapping
$T_0M/T_0^cM\to T_0M'/T_0^cM'$ represented by $G(0,w)$ is of
formal rank equal to $1$. Equivalently, 
$(\partial G/\partial w)(0)\neq 0$.
\end{itemize}
\smallskip
The multiplicity ${\rm Mult}\,(H,0)$ is independent of normal
coordinates, so it is a meaningful invariant of $h$ at an arbitrary
point of $M$. In normal coordinates, essential finiteness of $M$ at
$p$ is characterized by the finite codimensionality in 
$\C\{t\}$ of the ideal
generated by the $\Theta_\beta(t)$ for all $\beta\in\N^{n-1}$.  
This codimension is independent of coordinates and
denoted by ${\rm Ess\,Type}\,(M,p)$. Recall that $M\backslash E_M$ is
defined to be the set of points $q\in M$ at which the mapping
$t\mapsto (\Theta_\beta(t))_{\beta\in\N^{n-1}}$ is of rank $n$
in coordinates vanishing at $q$. Consequently

\def\thelemma{9.9}\begin{lemma}
For every $q\in M\backslash E_M$, we have 
${\rm Ess\,Type}\,(M,q)=1$.
\end{lemma}

A refinement of the analytic reflection principle proved in [BR1] is
as follows ([BR2, Theorem~6]).

\def\thelemma{9.10}\begin{lemma}
The $\mathcal{C}^\infty$-smooth CR mapping $h$ extends holomorphically
to a neighborhood of a point $q\in \C^n$ provided that in normal
coordinates centered at $q$ and at $q'=h(q)$, the normal component $g$
of $h$ is {\it not}\, flat at the origin, namely its formal power
series $G$ does not vanish identically.
\end{lemma}

Furthermore, in the case where the mapping $h$ extends holomorphically
at one point, four interesting nondegeneracy properties hold:

\def\thelemma{9.11}\begin{lemma}
With the same assumptions as in Theorem~9.3, let $q\in M$, let
$q':=h(q)$ and assume that $h$ extends holomorphically to a
neighborhood of $q$. Then
\begin{itemize}
\item[{\bf (1)}]
The induced differential $dh: T_qM/T_q^cM\to T_{q'}M'/T_{q'}^cM'$ 
is of rank $1$.
\item[{\bf (2)}]
The mapping $h$ is of finite multiplicity $m:={\rm Mult}\,(h,q)<\infty$
and $h$ is a local $m$-to one holomorphic mapping in a neighborhood of
$q$.
\item[{\bf (3)}]
We have the multiplicative relation 
\def\theequation{9.12}\begin{equation}
{\rm Ess\,Type}\,(M,q)={\rm Mult}\, (h,q) \cdot {\rm
Ess\,Type}\,(M',q').
\end{equation}
\item[{\bf (4)}]
If $q\in M\backslash E_M$, then $h$ is a local biholomorphism
at $q$.
\end{itemize}
\end{lemma}

\subsection*{9.13.~Installation of the proof of Theorem~9.3}
Let $E_{\rm na}$ be the closed set of points of $M$ at which the
mapping $h$ is not real analytic. By Lemma~9.6, the complement
$M\backslash E_{\rm na}$ is nonempty and in fact dense in $M$.
If $E_{\rm na}$ is empty, then Theorem~9.3 is proved, 
gratuitously. As in \S2 and \S5 above, we shall assume that
$E_{\rm na}$ is nonempty and we shall construct a contradiction
by showing that there exists in fact a point $p_1$ of $E_{\rm na}$ 
at which $h$ is real analytic. By Lemma~5.4, we are reduced to 
the following statement, which is analogous to Theorem~5.5.

\def\thetheorem{9.14}\begin{theorem}
Let $p_1\in E_{\rm na}$ and assume that there exists a real analytic
one-codimensional submanifold $M_1$ of $M$ with $p_1\in M_1$
which is generic in $\C^n$ such that $E_{\rm na}\backslash \{p_1\}$ is
completely contained in one of the two open sides of $M$ divided by
$M_1$, say in $M_1^+$, and such that $h$ is real analytic at every
point $q\in M\backslash E_{\rm na}$. Then $h$ is real
analytic at $p_1$.
\end{theorem}

To prove this theorem, we shall start as follows. We remind that the
intrinsic exceptional locus $E_M$ of $M$ is of real codimension at
least two in $M$. At each point $q\in M\backslash E_M$, the
hypersurface $M$ is finitely nondegenerate. It follows from
Lemma~9.11{\bf (4)} that at each point $q\in M\backslash 
(E_M\cup E_{\rm na})$, the
mapping $h$ extends as a local biholomorphism from a 
neighborhood of $q$ in $\C^n$ onto a neighborhood of
$h(q)$ in $\C^n$. Consider the relative disposition of the center
point $p_1$ with respect to $E_M$. In principle, there are two cases
to be considered.  Either $p_1\in E_M$ or $p_1\in M\backslash E_M$. In
both cases, we have the following useful existence property.

\def\thelemma{9.15}\begin{lemma}
There exists a small two-dimensional open real analytic manifold 
$K$ passing through $p_1$ and contained in $M$ such that
\begin{itemize}
\item[{\bf (1)}]
$K$  is transversal to $M_1$.
\item[{\bf (2)}]
$K\cap E_M=\{p_1\}$ and the line $T_{p_1}K\cap T_{p_1}M_1$
is not contained in $T_{p_1}^cM$.
\end{itemize}
\end{lemma}

\proof
Indeed, introducing real analytic coordinates on $M$, this follows
from a more general statement. Given a locally defined real analytic
set $E$ in $\R^\nu$ of dimension $1\leq \mu\leq \nu-1$ passing through
the origin, then for almost all $(\nu-\mu)$-dimensional linear planes
$K$ passing through the origin, the intersection of $K$ with $E$
consists of the singleton $\{0\}$ in a neighborhood of the origin.
\endproof

It follows from Lemma~9.15 that the intersection $K\cap M_1$
coincides with a geometrically smooth real analytic arc $\gamma_1$ passing
through $p_1$ which is not complex tangential at $p_1$. By
construction, $\gamma_1\backslash \{p_1\}$ is contained in the locus
$M\backslash E_{\rm na}$ where $h$ is already real analytic.
Moreover, $\gamma_1\backslash \{p_1\}$ is also contained in
$M\backslash E_M$. Its complexification $(\gamma_1)^c$ is a
complex disc transversal to $M$ with $(\gamma_1)^c\cap M=\gamma_1$.
Recall that $h$ already extends holomorphically to a one-sided
neighborhood $D$ of $M$. To fix ideas, we can assume that $D$ is in
the lower side $M^-$ of $M$ in $\C^n$. Moreover, $h$ extends
holomorphically to an open neighborhood $\Omega$ of $M\backslash 
E_{\rm na}$ in $\C^n$.  We
choose {\it normal}\, coordinates $t$ vanishing at $p_1$ in which the
equation of $M$ is of the form $\bar w=\Theta(\bar z,t)$, with
$\Theta(0,t)\equiv w$. Especially, we choose such coordinates 
in order that $\gamma_1$ coincides with a small 
neighborhood of the origin in the $u$-axis in these
normal coordinates, which is possible. 
Also, we choose some arbitrary {\it normal}\, coordinates $t'$
vanishing at $p_1':=h(p_1)$ in which the equation of $M'$ is of the
form $\bar w'=\Theta'(\bar z',t')$, with $\Theta'(0,t')\equiv w'$.
We denote the mapping by $h=(f,g)=(f_1,\ldots,f_{n-1},g)$ in these
coordinates.

Suppose for a while that we have proved that the normal component $g$
of the mapping extends holomorphically to a neighborhood of the point
$p_1$ in the transverse holomorphic disc $(\gamma_1)^c$, which
coincides with a small neighborhood of the origin in the $w$-axis.
Notice that we speak only of holomorphic extension to the single
transverse holomorphic disc passing through $p_1$, because our method
below will not give more.  Then we claim that the proofs of
Theorems~9.14 and~9.3 are achieved. Indeed, it suffices to show that
the holomorphic extension $g(0,w)$ at $w=0$ does not vanish identically, since
then it follows afterwards that the Taylor series $G$ at the origin of the
normal component $g$ does not vanish identically, whence $h$ extends
holomorphically at $p_1$ thanks to Lemma~9.10.  To prove that the
extension $g(0,w)$ is nonzero, we reason as follows.  According to
Lemma~9.11{\bf (1)}, at every point $q\in \gamma_1$ sufficiently close
to $p_1$ and different from $p_1$, the induced differential $dh:
T_qM/T_q^cM\to T_{q'}M'/T_{q'}^cM'$ is of rank one. This entails that
the differential $\partial_w g(0,w)$ is nonzero at $w:=w_q$, which
shows that the holomorphic extension $g(0,w)$ does not vanish
identically, as desired.  In summary, to prove Theorems~9.14 and~9.3,
it remains to establish the following crucial statement.

\def\thelemma{9.16}\begin{lemma} 
The $\mathcal{C}^\infty$-smooth
restrictions $f_1\vert_{\gamma_1},\dots, f_{n-1}\vert_{\gamma_1}$ and
$g\vert_{\gamma_1}$ extend holomorphically to a small neighborhood of
$p_1$ in the complex disc $(\gamma_1)^c$.
\end{lemma}

\subsection*{9.17.~Holomorphic extension to a transverse holomorphic disc}
This subsection is devoted to the proof of Lemma~9.16.  Using the
manifold $K$ of Lemma~9.15, we can include $\gamma_1$ into a
one-parameter family $\gamma_s$ of real analytic arcs, with $s_1<s\leq
1$, contained in $K$ which foliate $K\cap M_1^-$ for $s_1<s<1$.  Since
for $s_1<s<1$, the arcs $\gamma_s$ are contained in $M_1^-$, we have
$\gamma_s \cap E_{\rm na}=\emptyset$. By Lemma~9.15, we also have the
important property $\gamma_s\cap E_M=\emptyset$.  We consider the
complexifications $(\gamma_s)^c$, which are transversal to $M$. One
half of the complex discs $(\gamma_s)^c$ is contained in $D$. The
crucial Lemma~9.19 below is extracted from [P4, Lemma~3.1] and is
particularized to our $\mathcal{C}^\infty$-smooth situation. In the
sequel, it will be applied to the one-dimensional domains of the
complex plane $\C$ defined by
\def\theequation{9.18}\begin{equation}
U_s:=(\gamma_s)^c\ \cap \ D \ \cap \ \{\vert w\vert < r\},
\end{equation}
where $r>0$ is sufficiently small and to certain antiholomorphic
functions to be defined later. First of all, we introduce some
notation. As the complex disc $(\gamma_s)^c$ is transverse to $M$
and almost parallel to the $w$-axis, it
follows that $U_s$ is a small one-dimensional simply connected
domain in $(\gamma_s)^c$
bounded by two real analytic parts which we shall denote by
$\delta_s\subset \gamma_s$ and by $\beta_s\subset \{\vert w\vert
=r\}\cap (\gamma_s)^c\cap D$. 
These two real analytic arcs join together at two
points $q_s^+\in\gamma_s$ and $q_s^-\in\gamma_s$, namely $\{q_s^-,q_s^
+\}=\gamma_s\cap \{w=r\}=\delta_s\cap \beta_s$.  Then the boundaries
$\delta_s$ and $\gamma_s$ depend real analytically on $s$, even in a
neighborhood of $s=1$. We consider the two {\it open}\, real analytic arcs
$\delta_s^\circ:=\delta_s\backslash \{q_s^-,q_s^+\}$ and similarly for
$\beta_s^\circ$. Here is the lemma.

\def\thelemma{9.19}\begin{lemma} Let $U_s\subset \C$ be a
one-parameter family of bounded simply connected domains in $\C$
having piecewise real analytic boundaries with two open pieces
$\delta_s^\circ$ and $\beta_s^\circ$ depending real-analytically on a
real parameter $s_1<s\leq 1$, let $\varphi_s$, $\psi_s$ be
antiholomorphic functions defined in $U_s$ which
depend $\mathcal{C}^\infty$-smoothly on $s$ and set
$\theta_s:=\varphi_s/\psi_s$. Assume that the following four
conditions hold.
\begin{itemize}
\item[{\bf (1)}] For $s<1$, the two functions $\varphi_s$ and $\psi_s$
extend antiholomorphically to a certain neighborhood of
$\overline{U}_s$ in $\C$ and there exists a point $p_1\in
\delta_1^\circ$ so that $\varphi_1$ and $\psi_1$ extend
antiholomorphically to a neighborhood of $\overline{U_1}\backslash
\{p_1\}$ in $\C$ and $\mathcal{C}^\infty$-smoothly up to the open arc
$\delta_1^\circ$.
\item[{\bf (2)}]
The quotient $\theta_1:=\varphi_1/\psi_1$ is of class
$\mathcal{C}^\infty$ over $\delta_1^\circ$.
\item[{\bf (3)}] 
For $s<1$, the function $\psi_s$ does {\rm not}
vanish on $\partial U_s$ and there exists a constant $C>0$ such that
$\vert \theta_s\vert \leq C$ on $\partial U_s$ for all $s_1<s<1$.
\item[{\bf (4)}] 
The function $\psi_1$ does not vanish on
$\overline{U_1}\backslash \{p_1\}$.
\end{itemize}
Then the quotient $\theta_1$ satisfies $\vert \theta_1\vert \leq C$ on
$U_1$ and it extends as an antiholomorphic function to $U_1$ which is of
class $\mathcal{C}^\infty$ up to the open real analytic piece
$\delta_1^\circ$ of the boundary.
\end{lemma}

\proof In view of the nonvanishing of $\psi_s$ in $\partial U_s$, the function
$\psi_s$ has in $U_s$ a certain number $m$ (counting multiplicities)
of zeros which is constant for all $s_1<s<1$. Using a conformal
isomorphism of $U_s$ with the unit disc and an antiholomorphic
Blaschke product, we can construct an antiholomorphic function $b_s$
on $U_s$ extending $\mathcal{C}^\infty$-smoothly 
to the boundary with $\vert b_s\vert
=1$ on $\partial U_s$ such that the $m$ zeros of $b_s$ coincide with
the $m$ zeros of $\psi_s$. Then 
$b_s\theta_s$ is holomorphic in $U_s$ for 
$s_1<s<1$. It follows from the maximum principle that $\vert
b_s\,\theta_s\vert\leq C$ on $\overline{U}_s$ for all $s_1<s<1$.
Since $\psi_1\neq 0$ in $\overline{U_1}\backslash \{p_1\}$, when $s\to
1$, all zeros of the function $\psi_s$ converge to the single point
$p_1\in\partial U_1$. From the form of a Blaschke product,
we observe that $\lim_{s\to 1}\, \vert b_s(z)\vert =1$ for every point
$z\in U_1$. Therefore, for $z\in U_1$, we have
\def\theequation{9.20}\begin{equation}
\vert \theta_1(z)\vert =
\lim_{s\to 1}\,\vert \theta_s(z)\vert =
\lim_{s\to 1}\,\vert b_s(z)\,\theta_s(z)\vert \leq C.
\end{equation}
So the function $\theta_1$ is bounded in $U_1$. Since its boundary
value $\varphi_1/\psi_1$ is of class $\mathcal{C}^\infty$ on
$\delta_1^\circ$, it follows that the antiholomorphic
function $\theta_1$ extends
$\mathcal{C}^\infty$-smoothly up to $\delta_1^\circ\cup
\beta_1^\circ$.  The proof of Lemma~9.19 is complete.
\endproof

We can bow begin the proof of Lemma~9.16. Let $\overline{L}_1,\dots,
\overline{L}_{n-1}$ denote the commuting basis of $T^{0,1}M$ given
by $\overline{L}_j={\partial \over \partial \bar z_j}+
\Theta_{\bar z_j}(\bar z,t)\,{\partial\over \partial \bar w}$, for
$j=1,\dots,n-1$. In a neighborhood of the arc $\gamma_s$ for $s<1$, 
the mapping $h$ extends holomorphically as a local biholomorphism.
It follows that the determinant
\def\theequation{9.21}\begin{equation}
{\rm det}\, (\overline{L}_j\, \overline{f_{k}}(\bar t))_{
1\leq j,k\leq n-1}:=
\mathcal{D}(z,\bar t,\{\partial_{\bar t_l}
\overline{f_{k}}(\bar t) \}_{1\leq l\leq n, \, 1\leq k\leq n-1})
\end{equation}
does not vanish for $t\in M$ in a neighborhood of $\gamma_s$.  Also,
it extends as a certain antiholomorphic function to the domain
$U_s$. Let us denote this extension by $\psi_s$.  In order that the
function $\psi_s$ satisfies the assumption {\bf (4)} of Lemma~9.19, we
first observe that the determinant~(9.21) does not vanish on the part
$\delta_1\backslash\{p_1\}$ of $\partial U_1\backslash \{p_1\}$.
Indeed, since $h$ is real analytic at every point of
$\delta_1\backslash \{p_1\}$ and since $\delta_1\backslash \{p_1\}$ is
contained in $M\backslash E_M$, this follows from Lemma~9.11{\bf(4)}.
For the second part $\beta_1$ of $\partial U_1$, we observe that for
every small $r>0$ as in~\thetag{9.18}, 
the determinant~(9.10), extends as an antiholomorphic
function to $U_1$ and is in fact real analytic in a neighborhood of
$\beta_1$. Since the determinant~(9.21) does not vanish on
$\delta_1\backslash \{p_1\}$, there exist arbitrarily small $r>0$ such
that $\psi_1$ does not vanish over $\beta_1$. Shrinking $s_1$, we can
assume that $\psi_s$ does not vanish on $\beta_s$ for all $s_1<s\leq
1$. Finally, we know already that for $s<1$, the function $\psi_s$
does not vanish on $\delta_s$, thanks to the fact that $\gamma_s\cap
E_M$ is empty. Since $\psi_s$ does not vanish on the boundary
$\partial U_s$ for all $s_1<s<1$, it follows from Rouch\'e's theorem that
the number of zeros of $\psi_s$ in $U_s$ is constant equal to $m$
(counting multiplicities). Therefore, even for $s=1$, the function
$\psi_1$ has in $U_1$ not more than $m$ zeros.  Decreasing $r>0$ once
more, we can assume that $\psi_1$ does not vanish in $U_1$. This shows
that $\psi_s$ satisfies all the assumptions of Lemma~9.19.

Next, as the mapping $h$ is of class $\mathcal{C}^\infty$ over $M$, we
can apply the tangential Cauchy-Riemann derivations
$\overline{L}_1^{\beta_1}\cdots \overline{L}_{n-1}^{\beta_{n-1}}$,
$\beta\in \N^{n-1}$, of order $\vert\beta\vert$ infinitely many times
to the identity
\def\theequation{9.22}\begin{equation}
\overline{g(t)}=\Theta' (\overline{f(t)}, h(t)),
\end{equation}
which holds for all $t\in M$ in a neighborhood of $p_1$. As in \S7 above, 
using the nonvanishing of the determinant 
we get for all $\beta\in\N^{n-1}$ and for all $t\in \gamma_s$ with 
$s<1$ the following identities
\def\theequation{9.23}\begin{equation}
{1\over \beta!}\
{\partial^{\vert \beta\vert} \Theta'\over
\partial (z')^\beta}(\overline{f(t)},h(t))
={\mathcal{T}_\beta(z,\bar t,\{\partial_{\bar t}^\gamma 
\overline{h_j}(\bar t)\}_{1\leq j\leq n, \, \vert
\gamma\vert \leq \vert \beta \vert})\over 
[\mathcal{D}(z,\bar t,\{\partial_{\bar t_l}
\overline{f_{k}}(\bar t) \}_{1\leq l\leq n, 
\, 1\leq k\leq n-1})]^{2\vert \beta\vert -1}}.
\end{equation}
Precisely, the $\mathcal{T}_\beta$'s are holomorphic with respect to
$(z,\bar t)$ and relatively polynomial with respect to the jets
$\{\partial_{\bar t}^\gamma \overline{h_j}(\bar t) \}_{1\leq j\leq
n,\,\vert\gamma\vert\leq\vert\beta\vert}$. It follows that the
numerator $\mathcal{T}_\beta$ extends antiholomorphically to $U_s$ for
every $s_1<s\leq 1$ as a certain function which we shall denote by
$\varphi_{\beta,s}$.  We set
$\psi_{\beta,s}:= [\psi_s]^{ 2\vert\beta\vert-1}$. For $t\in\gamma_s
\subset M$ with $s<1$, let us rewrite~(9.23) as follows:
\def\theequation{9.24}\begin{equation}
(1/\beta!)\,
[\partial^{\vert \beta\vert} \Theta'/
\partial (z')^\beta](\overline{f(t)},h(t))
=\varphi_{\beta,s}(\bar t)/\psi_{\beta,s}(\bar t).
\end{equation}
As the left hand side of~(9.24) is of class $\mathcal{C}^\infty$ on
$\gamma_s$, it follows that the right hand side is of class
$\mathcal{C}^\infty$ on $\delta_s^\circ$, for all $s\leq 1$.  By
construction, for all $\beta\in\N^{n-1}$, the function
$\psi_{\beta,s}$ has no zeros on the boundary $\partial U_s$ for
$s<1$ and it also has no zeros on $\partial U_1\backslash
\{p_1\}$. Furthermore, these two functions 
$\varphi_{\beta,s}$ and $\psi_{\beta,s}$ both
extend antiholomorphically to a
neighborhood of $\overline{U_s}$ in $(\gamma_s)^c$ for $s<1$ and to a
neighborhood of $\overline{U_1}\backslash \{p_1\}$ for $s=1$. Let us
define $\theta_{\beta,s}:=\varphi_{\beta,s}/\psi_{\beta,s}$.  By
Lemma~9.19, for $s=1$, the functions $\theta_{\beta,1}$ extend
antiholomorphically to $U_1$ as bounded functions and
$\mathcal{C}^\infty$-smoothly up to the open real analytic arc
$\delta_1^\circ$. In summary, we have shown that for
all $\beta\in\N^{n-1}$, there exist functions $\theta_{\beta,1}(\bar t)$
defined for $t\in\delta_1$ and extending as antiholomorphic
functions to $U_1$ which are of class $\mathcal{C}^\infty$ up
to $\delta_1^\circ$ such that the following identities hold on 
$\delta_1$:
\def\theequation{9.25}\begin{equation}
(1/\beta!)\, [\partial^{\vert \beta\vert} \Theta'/ \partial
(z')^\beta](\overline{f(t)},h(t)) =\theta_{\beta,1}(\bar t).
\end{equation}

Next, we may derive some polynomial identities in the spirit of [BJT],
[BR1]. By the relation~(9.25) written for $t:=p_1\in \delta_1^\circ$,
we see that $\theta_{\beta,1}(\bar p_1)=0$, because $h(p_1)=p_1'$
sends the origin $p_1$ (in the coordinate system $t$) to the origin
$p_1'$ (in the coordinate system $t'$) and because the coordinates are
normal.  As $M'$ is essentially finite at the origin, there exists an
integer $\kappa\in\N_*$ such that the ideal $(\Theta_\beta'(t'))_{
\vert \beta\vert \leq \kappa}$ is of finite codimension in $\C\{t'\}$.
It follows from~(9.24) and from a classical computation ({\it
cf.}~[BJT], [BR1]) that there exist analytic cooeficients $A_{j,k}$ in
their variables which vanish at the origin and integers $N_j\geq 1$
such that, after possibly shrinking $r>0$, we have
\def\theequation{9.26}\begin{equation}
h_j^{N_j}(t)+
\sum_{k=1}^{N_j}\,
A_{j,k}(\overline{f(t)},\{\theta_{\beta,1}(\bar t)\}_{\vert\beta\vert\leq
\kappa})\, h_j^{N_j-k}(t)
=0,
\end{equation}
for all $t\in \delta_1$.  It follows that these coefficients
$A_{j,k}$, considered as functions of one real variable in $\delta_1$,
extend as antiholomorphic functions to $U_1$. In summary, we have
constructed some polynomial identities for the components of the
mapping $h$ with antiholomorphic coefficients which hold only on the single
transverse half complex disc $U_1=(\gamma_1)^c\cap D$ in a neighborhood of
$p_1$. These polynomial identities are crucial to show that the
mapping $h$ restricted to $(\gamma_1)^c\cap D$ extends holomorphically
to a neighborhood of $p_1$ in $(\gamma_1)^c$.

Indeed, by following the last steps
of the general approach of [BJT], [BR1,\S7], we deduce
that the reflection function (as denoted in 
equation~(8.1) of [BR1,\S8]) extends holomorphically to a neighborhood
of the point $p_1$ in $(\gamma_1)^c$ as a function of one complex
variable $w$ (remember that $(\gamma_1)^c$ is contained in the
$w$-axis). We would like to mention that in the strongly pseudoconvex
case, such a holomorphic extension to a single transverse holomorphic
disc was first derived by Pinchuk in [P4] in the more general case where $h$
is only continuous at $p_1$ and real analytic in $M\backslash E_{\rm
na}$.  Finally, using the real analyticity of the reflection function,
using the $\mathcal{C}^\infty$-smoothness of $h\vert_{\gamma_1}$ and using
Puiseux series as in [BJT], we deduce that $h\vert_{\gamma_1}$
is real analytic at $p_1$. The proof of Theorem~9.3 is complete.
\qed

\smallskip
A careful inspection of the above arguments shows that 
there is no obvious possibility to get an extension 
to the complex discs $(\gamma_s)^c$ with a uniform control
of the size of the domains of extension. Only the extension to the
limit complex disc $(\gamma_1)^c$ can be obtained.

\subsection*{9.27.~Strong uniqueness principle for CR mappings}
We end up this section by an application of Theorem~9.2. A similar
application of Theorem~9.3 may be stated.

\def\thetheorem{9.28}\begin{theorem}
Let $h: M\to M'$ and $h^*: M\to M'$ be two $\mathcal{C}^\infty$-smooth
CR mappings between two connected, real analytic hypersurfaces in
$\C^n$ and let $p\in M$. If $M$ and $M'$ do not contain complex
curves, then there exists an integer $\kappa\in\N_*$ which depends
only on $p$, on $M$ and on $M'$ such that if the two $\kappa$-jets of $h$ and
$h^*$ coincide at $p$, then $h\equiv h^*$ over $M$.
\end{theorem}

\proof
By Theorem~9.2, we can assume that $h$ and $h^*$ are both holomorphic
in a neighborhood of $p$ and nonconstant. By Lemma~9.11, the two
mappings $h$ and $h_*$ satisfies the Hopf Lemma at $p$ and are of
finite multiplicity. It follows from a careful inspection of the
analytic versions of the reflection principle given in [BR1], [BR2]
that if $\kappa$ is large enough, then the two mappings $h$ and $h^*$
coincide in a neighborhood of $p$. In fact, the complete arguments
already appeared in a more general context in [BER3, Theorem~2.5]. 
Then $h\equiv h^*$ all over $M$ by analytic
continuation. For the particular
case of germs, Theorem~9.28 is conjectured in [BER4,~p.~238].
\endproof

\section*{\S10.~Open problems and conjectures}

In the celebrated article [DP2], the following
conjecture stated {\it without pseudoconvexity assumption}, was solved
in the case $n=2$.

\def\theconjecture{10.1}\begin{conjecture} 
Let $h:D\to D'$ be a proper
holomorphic mapping between two bounded domains in $\C^n$ $(n\geq 2)$
having real analytic and geometrically smooth boundaries. Then $h$
extends holomorphically to an open neighborhood of $\overline{D}$ in
$\C^n$.
\end{conjecture}

To the author's knowledge, the conjecture is open for $n\geq 3$.  In
fact, among other conjectures, it has been conjectured for a long time
that every such proper holomorphic mapping
$h: D\to D'$ extends {\it continuously}\, to the boundary $M$ of $D$
and that in this case, $h$ is real analytic at every point of $M$. In
the much easier case where $h$ extends $\mathcal{C}^\infty$-smoothly
up to $M$, Theorem~9.2 above, in which no formal rank assumption is
imposed on the Taylor series of $h$ at points of $M$, provides a
positive answer.  Analogously, in Theorems~1.2 and~1.9, it would be
very desirable to remove the diffeomorphism assumption and also the
$\mathcal{C}^\infty$-smoothness assumption. We have strongly used these
two assumptions in the proof and we have found no way to do without.
Nevertheless, inspired by above conjectures, 
it is natural to suggest the following two open problems.

\def\theconjecture{10.2}\begin{conjecture}
Let $h: M\to M'$ be a {\rm continuous} CR mapping between two globally
minimal real analytic hypersurfaces in $\C^n$ $(n\geq 2)$ and assume
that the holomorphic extension of $h$ to a global one-sided neighborhood $D$
of $M$ in $\C^n$ is of generic rank equal to $n$. Then the reflection
function extends holomorphically to a neighborhood of {\rm every} point
$p\times \overline{h(p)}$ in the graph of $\bar h$.
\end{conjecture}

The rank assumption is really necessary, as shown by the following
trivial example.  Let $M\subset \C^4$ be the product of
$\C_{z_2}^1\times \C_{z_3}^1$ with the unbounded representation of the
$3$-sphere given by the equation $w=\bar w+iz\bar z$, let $M'\subset
\C^4$ be given by $w'=\bar w'+iz_1'\bar z_1'+iz_2'\bar z_3'+ i\bar
z_2'z_3'$, let $h_2(z_1,w)$ be a CR function on $M$ independent of
$(z_2,z_3)$, of class $\mathcal{C}^\infty$, which does not extend
holomorphically to the pseudoconcave side of $M$ at any point. Then
the degenerate mapping $(z_1,z_2,z_3,w) \mapsto (z_1,h_2(z_1,w),0,w)$
maps $M$ into $M'$ but does not extend holomorphically to a
neighborhood of $M$ in $\C^2$. Suppose by contradiction that the
globally defined reflection function $\mathcal{R}_h'(t,\bar\nu')=
\bar\mu'-w-i\bar\lambda_1'z_1-i\bar\lambda_3'h_2(z_1,w)$ extends
holomorphically to a neighborhood of $0\times 0$ in $\C^4\times \C^4$.
Differentiating with respect to $\bar\lambda_3'$, we deduce that
$h_2(w_1,z)$ extends holomorphically at the origin in $\C^4$,
contradiction. In fact, to speak of the extendability of the reflection 
function, one has to choose for $M'$ the {\it minimal for inclusion}\,
real analytic subset containing the image $h(M)$, as argued in [Me5].
In the case where the generic complex rank of $h$ over $D$ equals $n$, 
then $M'$ necessarily is the minimal for inclusion real analytic
set containing $h(M)$. This explains the rank assumption in Conjecture~10.2.

Finally, in the holomorphically nondegenerate case, we expect that $h$
be holomorphically extendable to a neighborhood of $M$.

\def\theconjecture{10.3}\begin{conjecture}
Let $h: M\to M'$ be a {\rm continuous} CR mapping between two globally
minimal real analytic hypersurfaces in $\C^n$ $(n\geq 2)$, assume
that the holomorphic extension of $h$ to a global one-sided neighborhood $D$
of $M$ in $\C^n$ is of generic complex rank equal to $n$ and assume that
$M'$ is holomorphically nondegenerate. Then $h$ is real analytic at
{\rm every} point of $M$.
\end{conjecture}

\vfill
\end{document}

%% file: figure6.pstex_t
\begin{picture}(0,0)%
\epsfig{file=figure6.pstex}%
\end{picture}%
\setlength{\unitlength}{3947sp}%
\begingroup\makeatletter\ifx\SetFigFont\undefined%
\gdef\SetFigFont#1#2#3#4#5{%
  \reset@font\fontsize{#1}{#2pt}%
  \fontfamily{#3}\fontseries{#4}\fontshape{#5}%
  \selectfont}%
\fi\endgroup%
\begin{picture}(6024,2424)(1305,-2769)
\put(2839,-2413){\makebox(0,0)[lb]{\smash{\SetFigFont{9}{10.8}{\familydefault}{\mddefault}{\updefault}$M_1$}}}
\put(6079,-2406){\makebox(0,0)[lb]{\smash{\SetFigFont{9}{10.8}{\familydefault}{\mddefault}{\updefault}$M_1'$}}}
\put(2580,-554){\makebox(0,0)[lb]{\smash{\SetFigFont{9}{10.8}{\familydefault}{\mddefault}{\updefault}$M_1^-$}}}
\put(3347,-550){\makebox(0,0)[lb]{\smash{\SetFigFont{9}{10.8}{\familydefault}{\mddefault}{\updefault}$M_1^+$}}}
\put(5730,-538){\makebox(0,0)[lb]{\smash{\SetFigFont{9}{10.8}{\familydefault}{\mddefault}{\updefault}${M_1'}^-$}}}
\put(6503,-544){\makebox(0,0)[lb]{\smash{\SetFigFont{9}{10.8}{\familydefault}{\mddefault}{\updefault}${M_1'}^+$}}}
\put(1971,-855){\makebox(0,0)[lb]{\smash{\SetFigFont{9}{10.8}{\familydefault}{\mddefault}{\updefault}$\gamma$}}}
\put(5169,-856){\makebox(0,0)[lb]{\smash{\SetFigFont{9}{10.8}{\familydefault}{\mddefault}{\updefault}$\gamma'$}}}
\put(6547,-1797){\makebox(0,0)[lb]{\smash{\SetFigFont{9}{10.8}{\familydefault}{\mddefault}{\updefault}$E_{\rm na}'$}}}
\put(3405,-1811){\makebox(0,0)[lb]{\smash{\SetFigFont{9}{10.8}{\familydefault}{\mddefault}{\updefault}$E_{\rm na}$}}}
\put(5040,-1392){\makebox(0,0)[lb]{\smash{\SetFigFont{9}{10.8}{\familydefault}{\mddefault}{\updefault}$q_1'$}}}
\put(5481,-1385){\makebox(0,0)[lb]{\smash{\SetFigFont{9}{10.8}{\familydefault}{\mddefault}{\updefault}$p_1'$}}}
\put(2280,-1394){\makebox(0,0)[lb]{\smash{\SetFigFont{9}{10.8}{\familydefault}{\mddefault}{\updefault}$p_1$}}}
\put(1875,-1395){\makebox(0,0)[lb]{\smash{\SetFigFont{9}{10.8}{\familydefault}{\mddefault}{\updefault}$q_1$}}}
\put(4193,-791){\makebox(0,0)[lb]{\smash{\SetFigFont{9}{10.8}{\familydefault}{\mddefault}{\updefault}$h$}}}
\put(1955,-2664){\makebox(0,0)[lb]{\smash{\SetFigFont{9}{10.8}{\familydefault}{\mddefault}{\updefault}{\sc Figure 1: Geometric similarity through the CR diffeomorphism $h$}}}}
\end{picture}

%% file: figure3.pstex_t
\begin{picture}(0,0)%
\epsfig{file=figure3.pstex}%
\end{picture}%
\setlength{\unitlength}{3947sp}%
\begingroup\makeatletter\ifx\SetFigFont\undefined%
\gdef\SetFigFont#1#2#3#4#5{%
  \reset@font\fontsize{#1}{#2pt}%
  \fontfamily{#3}\fontseries{#4}\fontshape{#5}%
  \selectfont}%
\fi\endgroup%
\begin{picture}(6024,2424)(1189,-2608)
\put(2597,-1654){\makebox(0,0)[lb]{\smash{\SetFigFont{10}{12.0}{\familydefault}{\mddefault}{\updefault}$A(\Delta)$}}}
\put(2223,-2025){\makebox(0,0)[lb]{\smash{\SetFigFont{10}{12.0}{\familydefault}{\mddefault}{\updefault}$D$}}}
\put(6017,-1807){\makebox(0,0)[lb]{\smash{\SetFigFont{10}{12.0}{\familydefault}{\mddefault}{\updefault}$D'$}}}
\put(4708,-1863){\makebox(0,0)[lb]{\smash{\SetFigFont{10}{12.0}{\familydefault}{\mddefault}{\updefault}$M$}}}
\put(6914,-1323){\makebox(0,0)[lb]{\smash{\SetFigFont{10}{12.0}{\familydefault}{\mddefault}{\updefault}$M'$}}}
\put(4791,-672){\makebox(0,0)[lb]{\smash{\SetFigFont{10}{12.0}{\familydefault}{\mddefault}{\updefault}$\Sigma_\gamma$}}}
\put(3387,-1467){\makebox(0,0)[lb]{\smash{\SetFigFont{10}{12.0}{\familydefault}{\mddefault}{\updefault}$\gamma$}}}
\put(3886,-1056){\makebox(0,0)[lb]{\smash{\SetFigFont{10}{12.0}{\familydefault}{\mddefault}{\updefault}$p_1$}}}
\put(3363,-1146){\makebox(0,0)[lb]{\smash{\SetFigFont{10}{12.0}{\familydefault}{\mddefault}{\updefault}$q$}}}
\put(3201,-829){\makebox(0,0)[lb]{\smash{\SetFigFont{10}{12.0}{\familydefault}{\mddefault}{\updefault}$q_1$}}}
\put(5985,-1494){\makebox(0,0)[lb]{\smash{\SetFigFont{10}{12.0}{\familydefault}{\mddefault}{\updefault}$A'(\Delta)$}}}
\put(6641,-1365){\makebox(0,0)[lb]{\smash{\SetFigFont{10}{12.0}{\familydefault}{\mddefault}{\updefault}$\gamma'$}}}
\put(6227,-858){\makebox(0,0)[lb]{\smash{\SetFigFont{10}{12.0}{\familydefault}{\mddefault}{\updefault}$(\gamma')^c$}}}
\put(1712,-437){\makebox(0,0)[lb]{\smash{\SetFigFont{10}{12.0}{\familydefault}{\mddefault}{\updefault}$S_{\bar q}$}}}
\put(5095,-1349){\makebox(0,0)[lb]{\smash{\SetFigFont{10}{12.0}{\familydefault}{\mddefault}{\updefault}$h$}}}
\put(1796,-2474){\makebox(0,0)[lb]{\smash{\SetFigFont{9}{10.8}{\familydefault}{\mddefault}{\updefault}{\sc Figure 2: The domain and its head covered by a Levi-flat hat}}}}
\end{picture}

%% file: figure4.pstex_t
\begin{picture}(0,0)%
\epsfig{file=figure4.pstex}%
\end{picture}%
\setlength{\unitlength}{3947sp}%
\begingroup\makeatletter\ifx\SetFigFont\undefined%
\gdef\SetFigFont#1#2#3#4#5{%
  \reset@font\fontsize{#1}{#2pt}%
  \fontfamily{#3}\fontseries{#4}\fontshape{#5}%
  \selectfont}%
\fi\endgroup%
\begin{picture}(6024,2424)(112,-1708)
\put(3431,-891){\makebox(0,0)[lb]{\smash{\SetFigFont{10}{12.0}{\familydefault}{\mddefault}{\updefault}$M$}}}
\put(2442,-351){\makebox(0,0)[lb]{\smash{\SetFigFont{10}{12.0}{\familydefault}{\mddefault}{\updefault}$p_1$}}}
\put(1394,-940){\makebox(0,0)[lb]{\smash{\SetFigFont{10}{12.0}{\familydefault}{\mddefault}{\updefault}$A(\Delta)$}}}
\put(1375,-704){\makebox(0,0)[lb]{\smash{\SetFigFont{10}{12.0}{\familydefault}{\mddefault}{\updefault}$D$}}}
\put(372,-975){\makebox(0,0)[lb]{\smash{\SetFigFont{10}{12.0}{\familydefault}{\mddefault}{\updefault}$\Omega$}}}
\put(196,-278){\makebox(0,0)[lb]{\smash{\SetFigFont{10}{12.0}{\familydefault}{\mddefault}{\updefault}$\omega(\Sigma_\gamma)$}}}
\put(3005,-1327){\makebox(0,0)[lb]{\smash{\SetFigFont{10}{12.0}{\familydefault}{\mddefault}{\updefault}$A_\sigma(\Delta)$}}}
\put(1817,439){\makebox(0,0)[lb]{\smash{\SetFigFont{10}{12.0}{\familydefault}{\mddefault}{\updefault}$q_1$}}}
\put(4001,-703){\makebox(0,0)[lb]{\smash{\SetFigFont{10}{12.0}{\familydefault}{\mddefault}{\updefault}$\Sigma_{\gamma'}'$}}}
\put(4172,586){\makebox(0,0)[lb]{\smash{\SetFigFont{7}{8.4}{\familydefault}{\mddefault}{\updefault}{\sc  Smoothing the corners}}}}
\put(5791,-968){\makebox(0,0)[lb]{\smash{\SetFigFont{10}{12.0}{\familydefault}{\mddefault}{\updefault}$M'$}}}
\put(5011,-1013){\makebox(0,0)[lb]{\smash{\SetFigFont{10}{12.0}{\familydefault}{\mddefault}{\updefault}$A'(\Delta)$}}}
\put(4366,-1178){\makebox(0,0)[lb]{\smash{\SetFigFont{10}{12.0}{\familydefault}{\mddefault}{\updefault}$D'$}}}
\put(5441,-698){\makebox(0,0)[lb]{\smash{\SetFigFont{10}{12.0}{\familydefault}{\mddefault}{\updefault}$p_1'$}}}
\put(5596, -6){\makebox(0,0)[lb]{\smash{\SetFigFont{5}{6.0}{\familydefault}{\mddefault}{\updefault}$A'(\Delta)$}}}
\put(3896,-298){\makebox(0,0)[lb]{\smash{\SetFigFont{10}{12.0}{\familydefault}{\mddefault}{\updefault}$h$}}}
\put(5256,162){\makebox(0,0)[lb]{\smash{\SetFigFont{8}{9.6}{\familydefault}{\mddefault}{\updefault}$\gamma'$}}}
\put(5599,406){\makebox(0,0)[lb]{\smash{\SetFigFont{8}{9.6}{\familydefault}{\mddefault}{\updefault}$(\gamma')^c$}}}
\put(3597,514){\makebox(0,0)[lb]{\smash{\SetFigFont{10}{12.0}{\familydefault}{\mddefault}{\updefault}$\Sigma_\gamma$}}}
\put(266,-1609){\makebox(0,0)[lb]{\smash{\SetFigFont{9}{10.8}{\familydefault}{\mddefault}{\updefault}{\sc Figure 3: Envelope of holomorphy of the domain and its Levi-flat hat}}}}
\put(4968,-586){\makebox(0,0)[lb]{\smash{\SetFigFont{10}{12.0}{\familydefault}{\mddefault}{\updefault}$q_1'$}}}
\put(4759,-202){\makebox(0,0)[lb]{\smash{\SetFigFont{8}{9.6}{\familydefault}{\mddefault}{\updefault}$(\gamma')^c$}}}
\put(2907,-91){\makebox(0,0)[lb]{\smash{\SetFigFont{10}{12.0}{\familydefault}{\mddefault}{\updefault}$\Sigma_\gamma^-$}}}
\end{picture}

%% file: figure1.pstex_t
\begin{picture}(0,0)%
\epsfig{file=figure1.pstex}%
\end{picture}%
\setlength{\unitlength}{3947sp}%
\begingroup\makeatletter\ifx\SetFigFont\undefined%
\gdef\SetFigFont#1#2#3#4#5{%
  \reset@font\fontsize{#1}{#2pt}%
  \fontfamily{#3}\fontseries{#4}\fontshape{#5}%
  \selectfont}%
\fi\endgroup%
\begin{picture}(6024,2424)(129,-1663)
\put(3691,-859){\makebox(0,0)[lb]{\smash{\SetFigFont{7}{8.4}{\familydefault}{\mddefault}{\updefault}$M_1'$}}}
\put(4586,131){\makebox(0,0)[lb]{\smash{\SetFigFont{7}{8.4}{\familydefault}{\mddefault}{\updefault}$E_{\rm na}'$}}}
\put(5741,-309){\makebox(0,0)[lb]{\smash{\SetFigFont{7}{8.4}{\familydefault}{\mddefault}{\updefault}$\gamma'$}}}
\put(3966,-524){\makebox(0,0)[lb]{\smash{\SetFigFont{7}{8.4}{\familydefault}{\mddefault}{\updefault}$p_1'$}}}
\put(5176,-354){\makebox(0,0)[lb]{\smash{\SetFigFont{7}{8.4}{\familydefault}{\mddefault}{\updefault}$q'$}}}
\put(326,-279){\makebox(0,0)[lb]{\smash{\SetFigFont{7}{8.4}{\familydefault}{\mddefault}{\updefault}$\gamma'$}}}
\put(442,-1574){\makebox(0,0)[lb]{\smash{\SetFigFont{9}{10.8}{\familydefault}{\mddefault}{\updefault}{\sc Figure 4: Construction of the generic wall by blowing out ellipsoids}}}}
\put(2191, 99){\makebox(0,0)[lb]{\smash{\SetFigFont{7}{8.4}{\familydefault}{\mddefault}{\updefault}$\Upsilon'$}}}
\put(967,561){\makebox(0,0)[lb]{\smash{\SetFigFont{7}{8.4}{\familydefault}{\mddefault}{\updefault}$L'$}}}
\put(3136,  9){\makebox(0,0)[lb]{\smash{\SetFigFont{7}{8.4}{\familydefault}{\mddefault}{\updefault}$Q_{\delta_1}'$}}}
\put(2181,-420){\makebox(0,0)[lb]{\smash{\SetFigFont{7}{8.4}{\familydefault}{\mddefault}{\updefault}$p'$}}}
\end{picture}

%% file: figure2.pstex_t
\begin{picture}(0,0)%
\epsfig{file=figure2.pstex}%
\end{picture}%
\setlength{\unitlength}{3947sp}%
\begingroup\makeatletter\ifx\SetFigFont\undefined%
\gdef\SetFigFont#1#2#3#4#5{%
  \reset@font\fontsize{#1}{#2pt}%
  \fontfamily{#3}\fontseries{#4}\fontshape{#5}%
  \selectfont}%
\fi\endgroup%
\begin{picture}(6024,3024)(741,-3381)
\put(4830,-1628){\makebox(0,0)[lb]{\smash{\SetFigFont{9}{10.8}{\familydefault}{\mddefault}{\updefault}$E_{\rm na}'$}}}
\put(3761,-1067){\makebox(0,0)[lb]{\smash{\SetFigFont{9}{10.8}{\familydefault}{\mddefault}{\updefault}${M_1'}^-$}}}
\put(4586,-762){\makebox(0,0)[lb]{\smash{\SetFigFont{9}{10.8}{\familydefault}{\mddefault}{\updefault}$M_1'$}}}
\put(4851,-1142){\makebox(0,0)[lb]{\smash{\SetFigFont{9}{10.8}{\familydefault}{\mddefault}{\updefault}${M_1'}^+$}}}
\put(3793,-1914){\makebox(0,0)[lb]{\smash{\SetFigFont{9}{10.8}{\familydefault}{\mddefault}{\updefault}$p_1'$}}}
\put(1323,-3259){\makebox(0,0)[lb]{\smash{\SetFigFont{8}{9.6}{\familydefault}{\mddefault}{\updefault}{\sc Figure~5: The family of real analytic arcs on the left side of the wall}}}}
\put(3791,-556){\makebox(0,0)[lb]{\smash{\SetFigFont{9}{10.8}{\familydefault}{\mddefault}{\updefault}$u'$}}}
\put(2846,-3057){\makebox(0,0)[lb]{\smash{\SetFigFont{9}{10.8}{\familydefault}{\mddefault}{\updefault}$y_1', z_*'$}}}
\put(5971,-1922){\makebox(0,0)[lb]{\smash{\SetFigFont{9}{10.8}{\familydefault}{\mddefault}{\updefault}$x_1'$}}}
\put(1352,-745){\makebox(0,0)[lb]{\smash{\SetFigFont{9}{10.8}{\familydefault}{\mddefault}{\updefault}$\gamma_{q_1'}'$}}}
\put(2915,-1932){\makebox(0,0)[lb]{\smash{\SetFigFont{9}{10.8}{\familydefault}{\mddefault}{\updefault}$q_1'$}}}
\put(2436,-1932){\makebox(0,0)[lb]{\smash{\SetFigFont{9}{10.8}{\familydefault}{\mddefault}{\updefault}$q_1'$}}}
\put(3433,-1924){\makebox(0,0)[lb]{\smash{\SetFigFont{9}{10.8}{\familydefault}{\mddefault}{\updefault}$q_1'$}}}
\end{picture}

%% file: figure7.pstex_t
\begin{picture}(0,0)%
\epsfig{file=figure7.pstex}%
\end{picture}%
\setlength{\unitlength}{3947sp}%
\begingroup\makeatletter\ifx\SetFigFont\undefined%
\gdef\SetFigFont#1#2#3#4#5{%
  \reset@font\fontsize{#1}{#2pt}%
  \fontfamily{#3}\fontseries{#4}\fontshape{#5}%
  \selectfont}%
\fi\endgroup%
\begin{picture}(6024,2424)(59,-1643)
\put(4925,-780){\makebox(0,0)[lb]{\smash{\SetFigFont{8}{9.6}{\familydefault}{\mddefault}{\updefault}$M$ }}}
\put(1321,-1063){\makebox(0,0)[lb]{\smash{\SetFigFont{8}{9.6}{\familydefault}{\mddefault}{\updefault}$\Omega$}}}
\put(1860,-1020){\makebox(0,0)[lb]{\smash{\SetFigFont{8}{9.6}{\familydefault}{\mddefault}{\updefault}$D$}}}
\put(3075,592){\makebox(0,0)[lb]{\smash{\SetFigFont{8}{9.6}{\familydefault}{\mddefault}{\updefault}$\Delta_n(0,\eta)$}}}
\put(551,-1559){\makebox(0,0)[lb]{\smash{\SetFigFont{9}{10.8}{\familydefault}{\mddefault}{\updefault}{\sc Figure 6: part of the envelope of holomorphy of the hat domain}}}}
\put(4318,179){\makebox(0,0)[lb]{\smash{\SetFigFont{8}{9.6}{\familydefault}{\mddefault}{\updefault}$\Sigma_{\gamma_{q_1}}^-$}}}
\put(4267,-1216){\makebox(0,0)[lb]{\smash{\SetFigFont{8}{9.6}{\familydefault}{\mddefault}{\updefault}$A_{q_1,\sigma}(\Delta)$}}}
\put(5181,621){\makebox(0,0)[lb]{\smash{\SetFigFont{8}{9.6}{\familydefault}{\mddefault}{\updefault}$\Sigma_{\gamma_{q_1}}$}}}
\put(1961,-760){\makebox(0,0)[lb]{\smash{\SetFigFont{6}{7.2}{\familydefault}{\mddefault}{\updefault}$A_{q_1,0}(\Delta)$}}}
\put(150,-384){\makebox(0,0)[lb]{\smash{\SetFigFont{8}{9.6}{\familydefault}{\mddefault}{\updefault}$\omega(\Sigma_{\gamma_{q_1}})$}}}
\put(2310,-339){\makebox(0,0)[lb]{\smash{\SetFigFont{6}{7.2}{\familydefault}{\mddefault}{\updefault}$\gamma_{q_1}$}}}
\put(5049,-258){\makebox(0,0)[lb]{\smash{\SetFigFont{8}{9.6}{\familydefault}{\mddefault}{\updefault}$\Sigma_{\gamma_{p_1}}$}}}
\put(3504, 19){\makebox(0,0)[lb]{\smash{\SetFigFont{6}{7.2}{\familydefault}{\mddefault}{\updefault}$\gamma_{p_1}$}}}
\put(3439,-126){\makebox(0,0)[lb]{\smash{\SetFigFont{6}{7.2}{\familydefault}{\mddefault}{\updefault}$p_1$}}}
\end{picture}

%% file: figure8.pstex_t
\begin{picture}(0,0)%
\epsfig{file=figure8.pstex}%
\end{picture}%
\setlength{\unitlength}{3947sp}%
\begingroup\makeatletter\ifx\SetFigFont\undefined%
\gdef\SetFigFont#1#2#3#4#5{%
  \reset@font\fontsize{#1}{#2pt}%
  \fontfamily{#3}\fontseries{#4}\fontshape{#5}%
  \selectfont}%
\fi\endgroup%
\begin{picture}(6024,2477)(91,-1733)
\put(4831,-1239){\makebox(0,0)[lb]{\smash{\SetFigFont{9}{10.8}{\familydefault}{\mddefault}{\updefault}$M$}}}
\put(1923,489){\makebox(0,0)[lb]{\smash{\SetFigFont{9}{10.8}{\familydefault}{\mddefault}{\updefault}$M_1^-$}}}
\put(4075,-891){\makebox(0,0)[lb]{\smash{\SetFigFont{9}{10.8}{\familydefault}{\mddefault}{\updefault}$M_1^+$}}}
\put(3168,439){\makebox(0,0)[lb]{\smash{\SetFigFont{9}{10.8}{\familydefault}{\mddefault}{\updefault}$H_1^+$}}}
\put(2817,444){\makebox(0,0)[lb]{\smash{\SetFigFont{9}{10.8}{\familydefault}{\mddefault}{\updefault}$H_1^-$}}}
\put(3038,580){\makebox(0,0)[lb]{\smash{\SetFigFont{9}{10.8}{\familydefault}{\mddefault}{\updefault}$H_1$}}}
\put(3567,141){\makebox(0,0)[lb]{\smash{\SetFigFont{9}{10.8}{\familydefault}{\mddefault}{\updefault}$p_1,M_1,\gamma_0$}}}
\put(5366,-389){\makebox(0,0)[lb]{\smash{\SetFigFont{9}{10.8}{\familydefault}{\mddefault}{\updefault}$S_{\bar p_1}$}}}
\put(938,-1620){\makebox(0,0)[lb]{\smash{\SetFigFont{8}{9.6}{\familydefault}{\mddefault}{\updefault}{\sc Figure 7: The Segre variety $S_{\bar p_1}$ intersects $D$ left to $H_1$ near $p_1$}}}}
\put(5458,-790){\makebox(0,0)[lb]{\smash{\SetFigFont{9}{10.8}{\familydefault}{\mddefault}{\updefault}$\Sigma_{\gamma_{q_1}}$}}}
\put(2492,153){\makebox(0,0)[lb]{\smash{\SetFigFont{9}{10.8}{\familydefault}{\mddefault}{\updefault}$q_1,\gamma_{q_1}$}}}
\put(1265,539){\makebox(0,0)[lb]{\smash{\SetFigFont{9}{10.8}{\familydefault}{\mddefault}{\updefault}$\Omega$}}}
\put(203,-218){\makebox(0,0)[lb]{\smash{\SetFigFont{9}{10.8}{\familydefault}{\mddefault}{\updefault}$\omega(\Sigma_{\gamma_{q_1}})$}}}
\put(908,-1343){\makebox(0,0)[lb]{\smash{\SetFigFont{9}{10.8}{\familydefault}{\mddefault}{\updefault}$A_{q_1,\sigma}(\Delta)$}}}
\put(1670,-865){\makebox(0,0)[lb]{\smash{\SetFigFont{9}{10.8}{\familydefault}{\mddefault}{\updefault}$D$}}}
\put(1516,-477){\makebox(0,0)[lb]{\smash{\SetFigFont{7}{8.4}{\familydefault}{\mddefault}{\updefault}$\Delta_n(\tilde{q_1},\tilde{\rho}/2)$}}}
\end{picture}